\documentclass{article}

\usepackage[numbers,sort&compress]{natbib}

\usepackage{graphics}
\usepackage{graphicx}
\usepackage{color}
\usepackage{subfigure}

\usepackage{lscape}
\usepackage{rotating}

\usepackage{amsfonts}
\usepackage{amsmath}
\usepackage{amsthm}

\usepackage{sublabel}

\usepackage[all]{xypic}
\usepackage{multicol}

\def\R{\mathbb{R}}
\def\Hess#1{\operatorname{Hess}(#1)}
\def\eref#1{(\ref{#1})}
\def\diiv{\operatorname{div}}

\newcommand{\onehalf}{\frac{1}{2}}

\newcommand{\real}{\mathbb{R}}

\newcommand{\calV}{{\cal V}}
\newcommand{\calT}{{\cal T}}

\newcommand{\Div}{\operatorname{div}}
\newcommand{\Grad}{\nabla}

\newcommand{\Iso}{\operatorname{iso}}

\providecommand{\abs}[1]{\lvert#1\rvert}
\providecommand{\norm}[1]{\lVert#1\rVert}

\newcommand{\TV}{\text{TV}}

\newcommand{\Linf}{L^{\infty}}

\newcommand{\Honehalf}{H^{\onehalf}(\partial\Omega)}
\newcommand{\Hnegonehalf}{H^{-\onehalf}(\partial\Omega)}
\newcommand{\LambdaOpNorm}{{\cal L}(\Honehalf,\Hnegonehalf)}

\newcommand{\uhat}{\hat{u}}

\newcommand{\myphi}{\varphi}
\newcommand{\bigO}{{\cal O}}

\newcommand{\sh}{s^h}
\newcommand{\Qh}{Q^h}

\newcommand{\qh}{q^h}

\newcommand{\tr}{\text{tr}\,}

\newcommand{\minimise}{\text{minimise}}
\newcommand{\subjectto}{\text{subject to}}

\newcommand{\cotan}{\cot}

\newcounter{tempc}
\newenvironment{oldequation}[1]
{
\setcounter{tempc}{\value{equation}}

\begin{equation}
}
{
\end{equation}
\setcounter{equation}{\value{tempc}}

}

\swapnumbers
\theoremstyle{definition}
\newtheorem{theorem}{Theorem}[section]
\newtheorem{proposition}[theorem]{Proposition}

\newtheorem{definition}[theorem]{Definition}
\newtheorem{remark}[theorem]{Remark}
\newtheorem{algorithm}[theorem]{Algorithm}

\renewenvironment{itemize}{
  \begin{list}{$\bullet$}{
      \setlength{\itemsep}{0.5ex} \setlength{\parsep}{0.ex}
      \setlength{\partopsep}{0.ex} \setlength{\topsep}{+0.05ex}
      \setlength{\leftmargin}{3ex}
    } }  {
  \end{list}
}

\numberwithin{equation}{section}
\numberwithin{figure}{section}
\numberwithin{table}{section}

\newtheorem{Theorem}{Theorem}[section]

\theoremstyle{remark}

\theoremstyle{definition}

\theoremstyle{definition}

\begin{document} \author {Mathieu Desbrun, Roger D.\ Donaldson, and Houman Owhadi\footnote{California
Institute of Technology, MC 217-50 Pasadena, CA 91125, USA. \newline mathieu@cs.caltech.edu,
rdonald@acm.caltech.edu, owhadi@caltech.edu}}

\title {Discrete Geometric Structures in Homogenization and Inverse Homogenization with application to EIT.}

\date{\today} \maketitle

\begin{abstract}
  We introduce a new geometric approach for the homogenization and inverse
  homogenization of the divergence form elliptic operator with rough
  conductivity coefficients $\sigma(x)$ in dimension two.  We show that
  conductivity coefficients are in one-to-one correspondence with
  divergence-free matrices and convex functions $s(x)$ over the domain $\Omega$.
  Although homogenization is a non-linear and non-injective operator when
  applied directly to conductivity coefficients, homogenization becomes a linear
  interpolation operator over triangulations of $\Omega$ when re-expressed using
  convex functions, and is a volume averaging operator when re-expressed
  with divergence-free matrices.  We explicitly give the transformations which
  map conductivity coefficients into divergence-free matrices and convex
  functions, as well as their respective inverses.  Using optimal weighted
  Delaunay triangulations for linearly interpolating convex functions, we apply
  this geometric framework to obtain an optimally robust homogenization
  algorithm for arbitrary rough coefficients, extending the global optimality of
  Delaunay triangulations with respect to a discrete Dirichlet energy to
  weighted Delaunay triangulations.  Next, we consider inverse homogenization,
  that is, the recovery of the microstructure from macroscopic information, a
  problem which is known to be both non-linear and severly ill-posed.  We show
  how to decompose this reconstruction into a linear ill-posed problem and a
  well-posed non-linear problem.  We apply this new geometric approach to
  Electrical Impedance Tomography (EIT) in dimension two.  It is known that the EIT problem
  admits at most one isotropic solution.  If an isotropic solution exists, we
  show how to compute it from any conductivity having the same boundary
  Dirichlet-to-Neumann map. This is of practical
  importance since the EIT problem always admits a unique
  solution in the space of divergence-free matrices and is stable
 with respect to $G$-convergence in that space (this property fails  for isotropic matrices).  As such,
  we suggest that the space of convex functions is the natural space to use to parameterize solutions of the EIT problem.
\end{abstract}

\setcounter{tocdepth}{2} \tableofcontents

\section{Introduction}

In this paper, we introduce a new geometric framework of the homogenization and inverse homogenization
of the divergence-form elliptic operator \begin{equation}\label{kjewehkjhwkhewe}
    u \rightarrow -\Div(\sigma \Grad u)
\end{equation} where $\sigma$ is symmetric, uniformly elliptic and with entries $\sigma_{ij}\in
L^\infty$.  Owing to its physical interpretation, we refer to $\sigma$ as the conductivity.

The classical theory of homogenization is based on abstract operator convergence and deals with the
asymptotic limit of a sequence of operators of the form \eref{kjewehkjhwkhewe} parameterized by a
small parameter $\epsilon$.  We refer to $G$-convergence for symmetric operators, $H$-convergence for
non-symmetric operators and $\Gamma$-convergence for variational problems \cite{Mur78, Gio75,
  MR630747, MR0477444, MR0240443, MR506997, MR1968440}. We also refer to \cite{BeLiPa78} for the original formulation based on asymptotic analysis and \cite{JiKoOl91} for a review.

Instead of considering the homogenized limit of an $\epsilon$-family of operators of the form
\eref{kjewehkjhwkhewe}, we will construct in this paper a sequence of finite dimensional and low rank
operators approximating~\eqref{kjewehkjhwkhewe} with arbitrary bounded $\sigma(x)$. More precisely,
since no small parameter $\epsilon$ is introduced in the formulation, one has to understand
homogenization in the context of finite dimensional approximation using a parameter $h$ that
represents a computational scale determined by the available computational power and the desired
precision. We are motivated by the fact that in most engineering problems, one has to deal with a
given medium and not with an $\epsilon$-family of media.

 This observation gave rise to methods such as special
finite element methods, metric based upscaling and harmonic change of coordinates considered
in~\cite{BabOsb83, BabCalOsb94, owha2005, owha2007, owha2006, OwDes09, OwBer09}. This point of view
recovers not only results from classical homogenization with periodic or ergodic coefficients but also
allows for homogenization of a given medium with arbitrary rough coefficients.  In particular we need
not make assumptions of ergodicity, scale separation and we do not need to introduce small parameter
$\epsilon$.

Our formalism, in not relying on small parameter $\epsilon$, is closely related to numerical
homogenization which deals with coarse scale numerical approximations of solutions of~\eqref{e:basic}
below.  Here we refer to the subspace projection formalism~\cite{MR2399542}, the multiscale finite
element method~\cite{MR1455261}, the mixed multiscale finite element method~\cite{MR2231859}, the
heterogeneous multiscale method~\cite{MR2314852,
  EnSou08}, sparse chaos approximations~\cite{MR2317004, MR2399150}; finite
difference approximations based on viscosity solutions~\cite{MR2361302}, operator splitting
methods~\cite{MR2342991} and generalized finite element methods~\cite{MR2283892}. We refer to
\cite{MR2322432, MR2281625} for an numerical implementation of the idea of a global change of harmonic
coordinates for porous media and reservoir modeling.

\paragraph{Contributions.} In
this paper, we focus on the intrinsic geometric framework underlying homogenization. First we show
that conductivities $\sigma$ can be put into one-to-one correspondence with, that is, can be
parameterized by, symmetric definite positive divergence free matrices $Q$, and by convex functions
$s$ as well (Section~\ref{khslkjdhdhklkjhe}).  While the transformation which maps $\sigma$ into
effective conductivities $q^h$ per coarse edge element is a highly non-linear transformation
(Section~\ref{jhskdgkjsghdkj}), we show that homogenization in the space of symmetric definite
positive divergence free matrices $Q$ acts as volume averaging, and hence is linear, while
homogenization in the space of convex functions $s$ acts as a linear interpolation operator
(Section~\ref{lksldhslkhsdlkhjej}). Moreover, we show that homogenization as it is formulated here is
self-consistent and satisfies a semi-group property (Section~\ref{ksjdlkshdhksjdh}).

\begin{figure}
  \begin{center}
 \xymatrix{
*--------------{} & \txt{\footnotesize \it Physical \\\footnotesize \it conductivity space} &&
\txt{\footnotesize \it Divergence-free \\ \footnotesize \it matrix space } && \txt{\footnotesize \it
Convex functions\\ \footnotesize \it space }\\
 *-------{\txt{\Large $\Omega$}} & *+[F-,]{\txt{{\huge $\sigma$}}} \ar@{~>} @<1ex>
 [rr]^{\txt{\footnotesize non-linear (injection if $\sigma$ scalar)\\ \footnotesize  $\sigma$-harm.
 coord. \eref{e:harmonic},\eref{edgeconddfdsuct} }} \ar@{~>}[rrddd]_{\txt{\footnotesize  non-linear,
 \\ \footnotesize  non-injective \\ \footnotesize  \eref{edgeconduct}, \eref{jashgajsggsjhw}}}  & &
*+[F-,]{\txt{ {\huge $Q$}}} \ar@{~>}@<1ex>[ll]^{\txt{ \footnotesize $\frac{Q}{\sqrt{\det(Q)}}$-har.
coor. \eref{eQharmonic}, \eref{jwhgdjhwjh}} } \ar@{->}@<1ex>[rr]^{\txt{\footnotesize linear
bijection\\\footnotesize rotation+integration\eref{kjshwjhgkjwhg}}} \ar@{->}[ddd]|-{\txt{\footnotesize
volume-averaging\\ \footnotesize linear, non-injective\\ \eref{klyiuyoiyu}}} & & *+[F-,]{\txt{{\huge
$s$}}}\ar@{->}@<1ex>[ll]^{\txt{\footnotesize differentiation + rotation
\eref{ksdjskhdhskdhjd}}}\ar@{->}[ddd]|-{\txt{\footnotesize  linear-interpolation\\ \footnotesize
linear, non-injective\\ \eref{klkglkhg}}}\\ *--------------{} & && &&\\ *--------------{} & && &&\\
*--------------{\txt{\Large $\Omega_h$}} & && *+[F-,]{\txt{ {\huge $q^h$}}}
\ar@{->}@<1ex>[rr]^{\txt{\footnotesize linear bijection\\ \footnotesize rotation+integration
\eref{lkglkjhgkjg}}} && *+[F-,]{\txt{ {\huge $s^h$}}}\ar@{->}@<1ex>[ll]^{\txt{\footnotesize
differentiation+rotation \eref{lkglkjhgkjg}}}\\ *--------------{} & && \txt{\footnotesize \it
effective edge\\ \footnotesize \it conductivities}&&\\
  }
$\quad$

\end{center} \caption[Relationships between parameterizations]{Relationships between
  parameterizations of conductivity.  Straight and wavy lines represent linear
  and non-linear relationships, respectively. \label{f:relationships}}
\end{figure}
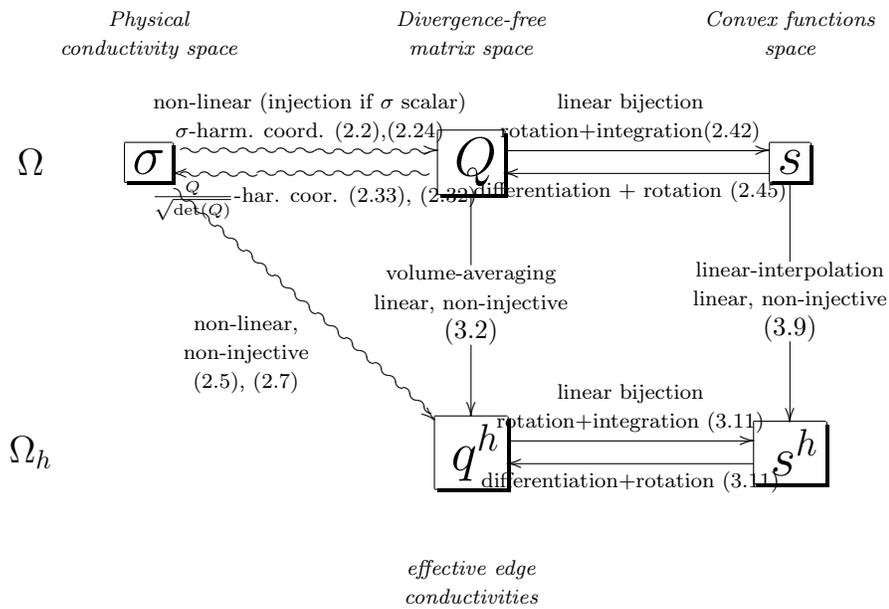

Hence, once formulated in the proper space, homogenization is a linear interpolation operator acting
on convex functions. We apply this observation to construct optimally robust and stable algorithms for
homogenizing divergence form equations with arbitrary rough coefficients by using weighted Delaunay
triangulations for linearly interpolating convex functions (Section~\ref{s:triangulate}).
Figure~\ref{f:relationships} summarizes relationships between the different parameterizations for
conductivity we study.

We use this new geometric framework for reducing the complexity of an inverse homogenization problem
(Section~\ref{lshkdgshdgjjhdsg}).  Inverse homogenization deals with the recovery of the physical
conductivity $\sigma$ from coarse scale information, for instance, from the set of effective
conductivities.  This problem is ill-posed insofar as it has no unique solution, and the space of
solutions is a highly nonlinear manifold.  We use this new geometric framework to re-cast inverse
homogenization into an optimization problem within a linear space.

We apply this result to Electrical Impedance Tomography (EIT), the problem of computing $\sigma$ from
Dirichlet and Neumann data measured on the boundary of our domain. First, we provide a new method for
solving EIT problems through parameterization via convex functions. Next we use this new geometric
framework to obtain new theoretical results on the EIT problem (Section~\ref{klcolhdckhcdeh}).
Although the EIT problem admits at most one isotropic solution, this isotropic solution may not exist
if the boundary data have been measured on an anisotropic medium.
We show that the EIT
problem admits a unique solution in the space of divergence-free matrices.  The uniqueness property has also been obtained in \cite{MR2308972}. When conductivities are endowed with the topology of $G$-convergence the inverse conductivity problem is discontinuous when restricted to isotropic matrices \cite{MR773707, MR2308972} and continuous when restricted to divergence-free matrices \cite{MR2308972}.

If an isotropic solution exists we
show here how to compute it for any conductivity having the same boundary data.
 This is of practical
importance since the medium to be recovered in a real application may not be isotropic and the
associated EIT problem may not admit an isotropic solution but if such an isotropic solution exists it can be computed from the divergence-free solution by solving  PDE \eref{eQharmsdsonic}.

As such, we suggest that the space
of divergence-free matrices parametrized by the space of convex functions is the natural space to look into for solutions of the EIT
problem.

\section{Homogenization and parametrization of the conductivity space.} \label{s:setup}

To illustrate our new approach, we will consider, as a first example, the homogenization of the
Dirichlet problem \begin{equation}
  \left\{
  \begin{aligned}
    -\Div\left(\sigma \Grad u\right) &= f, \enspace &x\in\Omega, \\
    u &= 0, \enspace &x\in\partial\Omega.
  \end{aligned}
  \right.
  \label{e:basic}
\end{equation} $\Omega$ is a bounded convex subset of $\real^d$ with a $C^2$ boundary, and $f\in
\Linf(\Omega)$.  The condition on $f$ can be relaxed to $f\in L^2(\Omega)$, but for the sake of
simplicity we will restrict our presentation to $f\in \Linf(\Omega)$.

Let $F:\Omega \rightarrow \Omega$ denote the {\em harmonic coordinates} associated with
\eref{e:basic}.  That is, $F(x)=\big(F_1(x),\ldots,F_d(x)\big)$ is a $d$-dimensional vector field
whose coordinates satisfy \begin{equation}\label{e:harmonic}
  \left\{
    \begin{aligned}
      \Div\left(\sigma \nabla F_i\right)&=0, \enspace &x\in\Omega,\\
      F_i(x)&=x_i, \enspace &x\in\partial\Omega.
    \end{aligned}
  \right.
\end{equation} In dimension $d=2$ it is known that $F$ is a homeomorphism from $\Omega$ onto $\Omega$
and $\det(\nabla F)>0$ a.e. \cite{An99, MR2001070, MR1892102}. For $d\geq 3$, $F$ may be non-injective, even
if $\sigma$ is smooth~\cite{An99, MR1892102,
  MR2073507}. We will restrict our presentation to $d=2$.

For a given symmetric matrix $M$, we denote by $\lambda_{\max}(M)$ and $\lambda_{\min}(M)$ its maximal
and minimal eigenvalues. Define
\begin{equation}\label{sjdsdsddsgkjdksjhdg}
  \mu:=\Big\|\frac{\lambda_{\max}((\nabla F)^T \nabla F )}{\lambda_{\min}((\nabla F)^T \nabla F
  )}\Big\|_{L^\infty(\Omega)}.
\end{equation}
We will call condition \eref{sjdsgkjdksjhdg} the following non-degeneracy condition on the anisotropy
of $(\nabla F)^T \nabla F$
\begin{equation}\label{sjdsgkjdksjhdg}
  \mu<\infty.
\end{equation}
For $d=2$, \eref{sjdsgkjdksjhdg} is always  satisfied if $\sigma$ is smooth
\cite{MR2001070} or even H\"{o}lder continuous \cite{MR1892102}.

\subsection{Homogenization as a non-linear operator} \label{jhskdgkjsghdkj}

Let $\Omega_h$ be a regular triangulation of $\Omega$ having resolution $h$. Let $X_h$ be the set of
piecewise linear functions on $\Omega_h$ with Dirichlet boundary conditions.  Let $\mathcal{N}_h$ be
the set of interior nodes of $\Omega_h$. For each node $i\in \mathcal{N}_h$, denote $\varphi_i$ the
piecewise linear nodal basis functions equal to $1$ on the node $i$ and $0$ on the other nodes. Let
$\mathcal{E}_h$ be the set of interior edges of $\Omega_h$, hence if $e\in \mathcal{E}_h$ then
$e=(i,j)$ where $i$ and $j$ are distinct interior nodes and share the edge of two triangles of
$\Omega_h$.

\begin{definition}[Effective edge conductivities]
  Let $q^h$ be the mapping from $\mathcal{E}_h$ onto $\R$, such that
  for $(i,j) \in \mathcal{E}_h$
  \begin{equation}\label{edgeconduct}
    q_{ij}^h:=-\int_{\Omega} (\nabla (\varphi_i \circ F))^T \sigma(x) \nabla
    (\varphi_j \circ F)\,dx.
  \end{equation}
  Observe that $q_{ij}^h=q_{ji}^h$, hence $q^h$ is only a function of {\em
    undirected} edges $(i,j)$. We refer to $q_{ij}^h$ as the {\em effective
    conductivity} of the edge $(i,j)$.
\end{definition}

Let $\mathcal{M}$ be the space of $2 \times 2$ uniformly elliptic, bounded and symmetric matrix fields
on $\Omega$.  Let $T_{q^h,\sigma}$ be
the operator mapping $\sigma$ onto $q^h$ defined by \eref{edgeconduct}.  Let $\mathcal{Q}_h$ be the
image of $T_{q^h,\sigma}$. \begin{equation}\label{jklghskjhdgjsh} \begin{split}
  T_{q^h,\sigma}\,:\,\mathcal{M} & \longrightarrow \mathcal{Q}_h\\
  \sigma & \longrightarrow T_{q^h,\sigma}[\sigma]:=q^h.
\end{split} \end{equation} Observe that $T_{q^h,\sigma}$ is both non-linear and non-injective.

Let $j\sim i$ be the set of interior nodes $j$, distinct from $i$, that share an edge with $i$.
\begin{definition}[Homogenized problem] Consider the vector
  $(u^h_i)_{i\in \mathcal{N}_h}$ of $\R^{\mathcal{N}_h}$ such that for all
  $i\in\mathcal{N}_h$,
  \begin{equation}\label{jashgajsggsjhw}
    \sum_{j\sim i} q_{ij}^h(u^h_i-u^h_j)= \int_{\Omega} f(x) \varphi_i\circ F(x)\,dx.
  \end{equation}
  We refer to this finite difference problem for $(u^h_i)_{i\in \mathcal{N}_h}$
  associated to $q^h$ as the {\em homogenized problem.}
\end{definition}

The identification of effective edge conductivities and the homogenized problem is motivated by the
following theorem: \begin{theorem}\label{lshgdsljdhghdsgs}
  The homogenized problem \eref{jashgajsggsjhw} has a solution $(u^h_i)_{i\in
    \mathcal{N}_h}$ and it is unique.  Moreover, let $u$ be the solution of
  \eref{e:basic} and define
  \begin{equation}\label{kdsjdhgsjhgsjgh}
    u_h:=\sum_{i\in \mathcal{N}_h}  u_i^h \varphi_i \circ F.
  \end{equation}
  If condition \eref{sjdsgkjdksjhdg} holds, then
  \begin{equation}\label{kjwhlkwhel}
    \|u-u_h\|_{H^1_0(\Omega)}\leq C h \|f\|_{L^\infty(\Omega)}.
  \end{equation}
\end{theorem} \begin{remark} We refer to \cite{owha2005} and \cite{owha2007} for numerical results
associated with theorem \ref{lshgdsljdhghdsgs}.
\end{remark}

\begin{remark}
  The constant $C$ depends on $\|1/\lambda_{\min}(\sigma)\|_{L^\infty(\Omega)}$,
  $\|\lambda_{\max}(\sigma)\|_{L^\infty(\Omega)}$, $\Omega$, and
  $\mu$. Replacing
  $\|f\|_{L^\infty(\Omega)}$ by $\|f\|_{L^2(\Omega)}$ in \eref{kjwhlkwhel} adds
  a dependence of $C$ on $\big\|(\det (\nabla
  F))^{-1}\big\|_{L^\infty(\Omega)}$.
\end{remark}

\begin{remark}
Although the proof of the theorem shows a dependence of $C$ on $\mu$ associated with condition \eref{sjdsgkjdksjhdg}, numerical results in dimension two indicate that $C$ is mainly correlated with the contrast (minimal and maximal eigenvalues) of $a$. This is why we believe that there should be a way of proving \eref{lshgdsljdhghdsgs} without condition \eref{sjdsgkjdksjhdg}. We refer to sections 2 and 3 of \cite{MR2001070} for the detailed analysis of a similar condition (definition 2.1 of \cite{MR2001070}).
\end{remark}

\begin{remark}
  Problem~\eref{jashgajsggsjhw} and Theorem~\ref{lshgdsljdhghdsgs} represent an
  generalization of method I of \cite{BabCalOsb94} to non-laminar media (see
  also~\cite{owha2005}).
\end{remark}

\begin{remark}
  It is also proven in~\cite{owha2005} (proof of Theorem 1.14) that if $f\in
  L^\infty(\Omega)$ then there exist constants $C, \alpha>0$ such that $u\circ
  F^{-1} \in C^{1,\alpha}(\Omega)$ and
  \begin{equation}\label{kjwhlkwhelgt}
    \|\nabla (u\circ F^{-1})\|_{C^\alpha(\Omega)}\leq C \|f\|_{L^\infty(\Omega)},
  \end{equation}
  where constants $C$ and $\alpha$ depend on $\Omega$,
  $\|1/\lambda_{\min}(\sigma)\|_{L^\infty(\Omega)}$,
  $\|\lambda_{\max}(\sigma)\|_{L^\infty(\Omega)}$, and
  $\mu$. We also refer to \cite{BabCalOsb94} (for quasi-laminar media) and \cite{MR2001070} for similar observations (on connections with quasi-regular and quasi-conformal mappings)
\end{remark}

\begin{remark}
  Unlike a canonical finite element treatment, where we consider only
  approximation of the solution, here we are also considering approximation of
  the operator.  This consideration is important, for example, in multi-grid
  solvers which rely on a set of operators which are self-consistent over a
  range of scales.
\end{remark}

\begin{proof}[Proof of Theorem~\ref{lshgdsljdhghdsgs}]
  The proof is similar to the proofs of Theorems 1.16 and 1.23
  of~\cite{owha2005}; we also refer to~\cite{OwBer09}. For the sake of
  completeness we will recall its main lines.

  Write $Q$ the matrix~\eref{edgeconddfdsuct}.  Replacing $u$ by $\hat{u}\circ
  F$ in~\eref{e:basic} we obtain after differentiation and change of variables
  that $\hat{u}:=u\circ F^{-1}$ satisfies
\begin{equation}\label{jsdkjsgdgsdjg}
  \left\{
  \begin{aligned}
    -\sum_{i,j}Q_{ij} \partial_i \partial_j \hat{u}&= \frac{f}{\det(\nabla F)}\circ F^{-1}, \enspace
    &x\in\Omega, \\
    u &= 0, \enspace &x\in\partial\Omega.
  \end{aligned}
  \right.
\end{equation} Similarly, multiplying~\eref{e:basic} by test functions $\varphi \circ F$ (with
$\varphi$ satisfying a Dirichlet boundary condition), integrating by parts and using the change
variables $y=F(x)$ we obtain that \begin{equation}\label{kjwhlsdsddsdhel} \int_{\Omega}(\nabla
\varphi)^T Q
  \nabla \hat{u} =\int_{\Omega} \varphi \frac{f}{\det(\nabla F)}\circ F^{-1}.
\end{equation}

Observing that $\hat{u}$ satisfies the non-divergence form equation, we obtain, using Theorems 1.2.1
and 1.2.3 of \cite{MPG00}, that if $Q$ is uniformly elliptic and bounded, then $\hat{u}\in
W^{2,2}(\Omega)$ with \begin{equation}\label{kjxddddsdhel}
  \|\hat{u}\|_{W^{2,2}(\Omega)}\leq C \|f\|_{L^\infty(\Omega)}.
\end{equation} The constant $C$ depends on $\Omega$ and bounds on the minimal and maximal eigenvalues
of $Q$.  We have used the fact that the Cordes-type condition on $Q$ required by~\cite{MPG00}
simplifies for $d=2$.

Next, let $V_h$ be the linear space defined by $\varphi \circ F$ for $\varphi \in X_h$. Write $u_h$
the finite element solution of~\eref{e:basic} in $V_h$. Writing $u_h$ as in~\eref{kdsjdhgsjhgsjgh} we
obtain that the resulting finite-element linear system can be written as
\begin{equation}\label{jashgazdaassaswzsdjhw} -q_{ii}^h  u^h_i-\sum_{j\sim i} q_{ij}^h  u^h_j=
\int_{\Omega} f(x) \varphi_i\circ F(x)\,dx, \end{equation} for $i\in\mathcal{N}_h$.  We use
definition~\eqref{edgeconduct} for $q_{ij}^h$. Using the change of variables $y=F(x)$ we obtain that
$q_{ij}^h$ can be written \begin{equation}\label{edgesddsnzdduct}
  q_{ij}^h:=-\int_{\Omega} (\nabla \varphi_i )^T Q(x) \nabla \varphi_j \,dx.
\end{equation} Decomposing the constant function $1$ over the basis $\varphi_j$ we obtain that
\begin{equation}\label{edgesddsnzdssdssduiuct}
  -\int_{\Omega} (\nabla \varphi_i )^T Q(x).\nabla (\sum_{j} \varphi_j) \,dx=0,
\end{equation} from which we deduce that \begin{equation}\label{edgesddsdedsnzdsdssduct}
  q_{ii}^h +\sum_{j\sim i} q_{ij}^h =0.
\end{equation} Combining~\eref{edgesddsdedsnzdsdssduct} with~\eref{jashgazdaassaswzsdjhw} we obtain
that the vector $(u^h_i)_{i\in \mathcal{N}_h}$ satisfies equation~\eref{jashgajsggsjhw}.

Using the change of variables $y=F(x)$ in \begin{equation}\label{edgesddsnzdssdssduct} \int_{\Omega}
(\nabla (\varphi_i\circ F) )^T \sigma(x).\nabla u_h \,dx=\int_{\Omega} \varphi_i\circ F\, f,
\end{equation} we obtain that $\hat{u}_h:=u_h \circ F^{-1}$ satisfies
\begin{equation}\label{zdselsdsddsdhel} \int_{\Omega}(\nabla \varphi_i)^T Q \nabla \hat{u}_h
=\int_{\Omega} \varphi_i \frac{f}{\det(\nabla F)}\circ F^{-1}. \end{equation} Hence $\hat{u}_h$ is the
finite element approximation of $\hat{u}$. Using the notation $\sigma[v]:=\int_{\Omega} \nabla v^T
\sigma \nabla v$ we obtain through the change of variables $y=F(x)$ that $\sigma[v]=Q[v\circ F^{-1}]$.
It follows that \begin{equation}\label{zdssdjyuwyedsddsdhel} \sigma[u-u_h]= Q[\hat{u}-\hat{u}_h].
\end{equation} Since $\hat{u}_h$ minimizes $Q[\hat{u}-v]$ over $v\in X_h$ we obtain equation
\eref{kjwhlkwhel} from the $W^{2,2}$-regularity of $\hat{u}$ \eref{kjxddddsdhel}. \end{proof}

The fact that $q^h$, as a quadratic form on $\R^{\mathcal{N}_h}$, is positive definite can be obtained
from the following proposition: \begin{proposition}
  For all vectors $(v_i)_{i\in \mathcal{N}_h}\in\real^{\mathcal{N}_h}$,
\begin{equation} \sum_{i\sim j} v_i q_{ij}^hv_j=\int_{\Omega} (\nabla (v\circ F))^T \sigma \nabla
(v\circ F), \end{equation} where $v:=\sum_{i\in \mathcal{N}_h} v_i \; \varphi_i$. \end{proposition}

\begin{proof}
  The proof follows from first observing that
  \begin{equation}
    \sum_{i\sim j} v_i q_{ij}^hv_j= \int_{\Omega} (\nabla v)^T(y) Q(y) (\nabla v)(y)\,dy,
  \end{equation}
  then applying the change of variables $y=F(x)$.
\end{proof}

\begin{remark}
  Despite the fact that positivity holds for any triangulation $\Omega_h$, as we
  shall examine in Section~\ref{s:triangulate}, we can take advantage of the
  freedom to choose $\Omega_h$ to produce $q_{ij}^h$ which give linear systems
  representing homogenized problems~\eqref{jashgajsggsjhw} having optimal
  conditioning properties.
\end{remark}

\subsection{Parametrization of the conductivity space}\label{khslkjdhdhklkjhe}

We now take advantage of special properties of $\sigma$ when transformed by its harmonic coordinates
$Q$ to parameterize the space of conductivities.

\begin{definition}[Space of divergence-free matrices] We say that a matrix field
  $M$ on $\Omega$ is {\em divergence-free} if its columns are divergence-free
  vector fields.  That is, $M$ is divergence-free if for all $v\in
  C^\infty_0$ and $l\in \R^2$
\begin{equation}
  \int_\Omega (\nabla v)^T\,M.l=0
\end{equation} \end{definition}

\begin{definition}[Divergence-free conductivity] Given the conductivity $\sigma$
  associated to~\eref{kjewehkjhwkhewe} and a domain $\Omega$, define $Q$ to be
  the symmetric $2\times 2$ matrix given by
  \begin{equation}\label{edgeconddfdsuct}
    Q = F_{\ast}\sigma := \frac{(\nabla F)^T \sigma \nabla F}{\det(\nabla F)}\circ F^{-1}.
  \end{equation}
  Again, $F:\Omega \rightarrow \Omega$ are the harmonic
  coordinates~\eref{e:harmonic} associated to $\sigma$.
\end{definition}

\begin{proposition}[Properties of $Q$] \label{kljssdlhkshdlk}  $Q$ satisfies the following properties:
\begin{itemize}
\item $Q$ is
  positive-definite, symmetric and divergence-free.

\item $Q\in (L^1(\Omega))^{d\times d}$.

\item $\det(Q)$ is uniformly bounded away from $0$ and $\infty$.

\item $Q$ is bounded and uniformly
  elliptic if and only if $\sigma$ satisfies the non-degeneracy condition~\eref{sjdsgkjdksjhdg}.
  \end{itemize}
\end{proposition} \begin{proof} Equations \eref{jsdkjsgdgsdjg} and \eref{kjwhlsdsddsdhel} imply that
for all $\hat{u}\in H^1_0\cap H^2(\Omega)$ and all $\varphi\in C^\infty_0(\Omega)$
\begin{equation}\label{jsdkjsgdfrdgsdjg}
   \int_{\Omega}(\nabla \varphi)^T Q
  \nabla \hat{u}=  - \int_{\Omega} \varphi \sum_{i,j}Q_{ij} \partial_i \partial_j \hat{u}.
\end{equation} Let $l\in \R^d$, choosing $\hat{u}$ such that $\nabla \hat{u}=l$ on the support of
$\varphi$ we obtain that for all $l\in \R^d$ \begin{equation}\label{jsdkjsgdfrdddgsdjg}
   \int_{\Omega}(\nabla \varphi)^T Q\cdot l=  0
\end{equation} It follows by integration by parts that $\diiv(Q\cdot l)=0$ in the weak sense and hence
$Q$ is divergence-free (its columns are divergence-free vector fields, this has also been obtained in
  \cite{owha2005}).
The second and third part of the Proposition can be obtained from
  \begin{equation}\label{edgec1onddfdssdsuct}
    \det(Q)=\det\big(\sigma\circ F^{-1}\big),
  \end{equation}
  and
  \begin{equation}\label{edgeceo2nddfdssuct}
    \int_\Omega Q=\int_{\Omega}(\nabla F)^T \sigma \nabla F
  \end{equation}
  The last part of the Proposition can be obtained from the
  following inequalities (valid for $d=2$).  For $x\in \Omega$ a.e.,
  \begin{equation}\label{sjdsdkjd8ksjhdg}
    \lambda_{\max}(Q)\leq \lambda_{\max}(\sigma) \sqrt{\frac{\lambda_{\max}((\nabla F)^T \nabla F
    )}{\lambda_{\min}((\nabla F)^T \nabla F )}}
  \end{equation}
  \begin{equation}\label{sjdoiosgkjdk0sjhdg}
    \lambda_{\min}(Q)\geq \lambda_{\min}(\sigma) \sqrt{\frac{\lambda_{\min}((\nabla F)^T \nabla F
    )}{\lambda_{\max}((\nabla F)^T \nabla F )}}
  \end{equation}
\end{proof}

Write $T_{Q,\sigma}$ the operator mapping $\sigma$
onto $Q$ through equation~\eref{edgeconddfdsuct}.  That is, \begin{equation}
  \begin{split}
    T_{Q,\sigma}\,:\,\mathcal{M} & \longrightarrow \mathcal{M}_{\diiv}\\
    M & \longrightarrow T_{Q,\sigma}[M]:= \frac{(\nabla F_M)^T M \nabla F_M}{\det(\nabla F_M)}\circ
    F^{-1}_M,
  \end{split}
\end{equation} where $F_M$ are the harmonic coordinates associated to $M$ through
equation~\eref{e:harmonic} (for $\sigma\equiv M$) and $\mathcal{M}_{\diiv}$ is the image of $\mathcal{M}$ under the operator $T_{Q,\sigma}$. Observe (from Proposition \ref{kljssdlhkshdlk}) that $\mathcal{M}_{\diiv}$
is a space of $2 \times 2$ of symmetric, positive and divergence-free matrix fields on $\Omega$, with entries in $L^1(\Omega)$ and with determinants uniformly bounded away from $0$ and infinity.

Since for all $M \in \mathcal{M}_{\diiv}$,
$T_{Q,\sigma}[M]=M$ ($T_{Q,\sigma}$ is a non-linear projection onto $\mathcal{M}_{\diiv}$) it follows
that $T_{Q,\sigma}$ is a non-injective operator from $\mathcal{M}$ onto $\mathcal{M}_{\diiv}$. Now
denote $\mathcal{M}_{\Iso}$ the space of $2 \times 2$ isotropic, uniformly elliptic, bounded and
symmetric matrix fields on $\Omega$.
Hence matrices in $\mathcal{M}_{\Iso}$ are of the form $\sigma(x)I_d$ where $I_d$ is the $d\times d$
identity matrix and $\sigma(x)$ is a scalar function.

\begin{theorem}\label{kjlaehlwhelkjwh}
  The following statements hold in dimension $d=2$:
  \begin{enumerate}
  \item The operator $T_{Q,\sigma}$ is an injection from $\mathcal{M}_{\Iso}$
    onto $\mathcal{M}_{\diiv}$.
  \item
    \begin{equation}\label{jwhgdjhwjh}
      T_{Q,\sigma}^{-1}[Q]=\sqrt{\det(Q)\circ G^{-1}}\, I_d,
    \end{equation}
    where $G$ are the harmonic coordinates associated to
    $\frac{Q}{\sqrt{\det(Q)}}$.  That is, $G_i(x),i=1,2$ is the unique solution of
    \begin{equation}\label{eQharmonic}
      \left\{
        \begin{aligned}
          \Div \left(\frac{Q}{\sqrt{\det(Q)}} \nabla G_i\right)&=0 \enspace
          &x\in\Omega, \\
          G_i(x)&=x_i \enspace &x\in\partial\Omega.
        \end{aligned}
      \right.
    \end{equation}
  \item $G=F^{-1}$ where $G$ is the transformation defined by~\eqref{eQharmonic},
    and $F$ are the harmonic coordinates associated to
    $\sigma:=T_{Q,\sigma}^{-1}[Q]$ by~\eref{e:harmonic}.
  \end{enumerate}
\end{theorem}

\begin{remark}
 Observe that the non
  degeneracy condition~\eref{sjdsgkjdksjhdg} is not necessary for the validity of this theorem.
\end{remark} \begin{remark}
  $T_{Q,\sigma}$ is not surjective from $\mathcal{M}_{\Iso}$ onto
  $\mathcal{M}_{\diiv}$. This can be proven by contradiction by assuming $Q$ to
  be a non-isotropic constant matrix. Constant $Q$ is trivially divergence-free,
  yet it follows that $\sigma=\sqrt{\det(Q)}I_d$, $F(x)=x$ and $Q$ is isotropic,
  which is a contradiction.
\end{remark} \begin{remark}
  $T_{Q,\sigma}$ is not an injection from $\mathcal{M}$ onto
  $\mathcal{M}_{\diiv}$ .  However it is known~\cite{MR1899805} that for each
  $\sigma\in \mathcal{M}$ there exists a sequence $\sigma_\epsilon$ in
  $\mathcal{M}_{\Iso}$ $H$-converging towards $\sigma$.  (Moreover, this
  sequence can be chosen to be of the form $a(x,x/\epsilon)$, where $a(x,y)$ is
  periodic in $y$.)  Since $\mathcal{M}_{\Iso}$ is dense in $\mathcal{M}$
  with respect to the topology induced by $H$-convergence, and since
  $T_{Q,\sigma}$ is an injection from $\mathcal{M}_{\Iso}$, the scope of
  applications associated with the existence of $T_{Q,\sigma}^{-1}$ would not
  suffer from a restriction from $\mathcal{M}$ to $\mathcal{M}_{\Iso}$.
\end{remark} \begin{proof}[Proof of Theorem~\ref{kjlaehlwhelkjwh}]
  First observe that if $\sigma$ is scalar then we obtain from
  equation~\eref{edgeconddfdsuct} that
  \begin{equation}\label{edgeconddfdssuct}
    \det(Q)=\big(\sigma\circ F^{-1}\big)^2,
  \end{equation}
  and hence
  \begin{equation}\label{edgsdeeconddfdssuct}
    \sigma=\sqrt{\det(Q)\circ F}.
  \end{equation}

  Consider again equation~\eref{edgeconddfdsuct}.  Let $R$ be the $2\times 2$,
  $\frac{\pi}{2}-$rotation matrix in $\R^2$, that is,
  \begin{equation}\label{jshdkjhewkehj}
    R = \begin{pmatrix} 0 & -1 \\ 1 & 0 \end{pmatrix}.
  \end{equation}
  Observe that for a $2\times 2$ matrix $A$,
  \begin{equation}\label{klejwkheke}
    (A^{-1})^T=\frac{1}{\det (A)}R\,A\,R^T.
  \end{equation}
  Write $G:=F^{-1}$. Recall that
  \begin{equation}\label{lhlkhelkhwep}
    \nabla G=(\nabla F)^{-1}\circ F^{-1}.
  \end{equation}
  Applying~\eref{lhlkhelkhwep} to~\eref{edgeconddfdsuct} gives
  \begin{equation}\label{kljkgwljegwh}
    Q \nabla G= \det(\nabla G)((\nabla G)^{-1})^T \sigma\circ G.
  \end{equation}
  Using $\sqrt{\det (Q)}=\sigma\circ G$ and applying equation~\eref{klejwkheke}
  to $((\nabla G)^{-1})^T$ we obtain that
  \begin{equation}\label{kljksdsdgwljegwh}
    \frac{Q}{\sqrt{\det(Q)}} \nabla G= R \; \nabla G \; R^T.
  \end{equation}
  Observing that in dimension two, for all functions $v\in H^1$,
  $\diiv(R\nabla v)=0$, and we obtain from~\eref{kljksdsdgwljegwh} that $G$
  satisfies equation~\eref{eQharmonic}. The boundary condition comes from the
  fact that $G=F^{-1}$, where $F$ is a diffeomorphism and $F(x)=x$ on $\partial
  \Omega$. Let us now show that equation~\eref{eQharmonic} admits a unique solution. If $G_i'$ is another solution of equation~\eref{eQharmonic} then
    \begin{equation}\label{kljdfdksdsdgwljegwh}
   \nabla (G_i-G_i') \frac{Q}{\sqrt{\det(Q)}} \nabla (G_i-G_i')=0
  \end{equation}
  Since $Q$ is positive with $L^1$ entries and $\det(Q)$ is uniformly bounded away from zero and infinity it follows that $\frac{Q}{\sqrt{\det(Q)}}$ is positive and its minimal  eigenvalue is bounded away from infinity almost everywhere in $\Omega$. It follows that $\nabla (G_i-G_i')=0$ almost everywhere in $\Omega$ and we conclude from the boundary condition on $G_i$ and $G_i'$ that
  $G_i=G_i'$ almost everywhere in $\Omega$.
\end{proof}

\begin{definition}[The space of convex functions]
  Consider the space of $W^{2,1}(\Omega)$ convex functions on
  $\Omega$ whose discriminants (determinant of the Hessian) are uniformly bounded away from zero and infinity. Write $\mathcal{S}$ the quotient set on that space defined by the
  equivalence relation: $s\sim s'$ if $s-s'$ is an affine function. Let $R$ be
  the rotation matrix \eref{jshdkjhewkehj}.
\end{definition}

\begin{theorem}[Scalar parameterization of conductivity in $\real^2$]\label{kjshgdkdjsgdhjg}
  For each $Q\in \mathcal{M}_{\diiv}$ there exists a unique $s\in \mathcal{S}$
  such that
  \begin{equation}\label{kjshwjhgkjwhg}
    \Hess{s}=R^T Q R,
  \end{equation}
  where $\Hess{s}$ is the Hessian of $s$.
\end{theorem} \begin{remark}
  Since $Q$ is positive-definite one concludes that $\Hess{s}$ is
  positive-definite, and thus, $s(x)$ is convex. Furthermore, the principal
  curvature directions of $s(x)$ are the eigenvectors of $Q$, rotated by
  $\pi/2$. Note that this geometric characteristic will be crucial later when we
  approximate $s(x)$ by piecewise-linear polynomials, which are not everywhere
  differentiable---but for which the notion of convexity is still well defined.
\end{remark} \begin{proof}
  In $\R^2$, the symmetry and divergence-free constraints on $Q$ reduce the
  number of degrees of freedom of $Q(x)$ to a single one.  This remaining degree
  of freedom is $s(x)$, our scalar convex parameterizing function.  To construct
  $s(x)$, observe that as a consequence of the Hodge decomposition, there exist
  functions $h, k \in W^{1,1}(\Omega)$ such that
  \begin{equation}
    Q = \begin{pmatrix} a & b \\ b & c \end{pmatrix}
    = \begin{pmatrix} h_y & k_y \\ -h_x & -k_x \end{pmatrix}
  \end{equation}
  where $a,b,c$ are scalar functions.  These choices ensure that the
  divergence-free condition is satisfied, namely that $a_x + b_y = b_x + c_y =
  0$.  Another application of the Hodge decomposition gives the existence of
  $s\in W^{2,1}(\Omega)$ such that $\Grad s = \left(-k,h\right)^T$.  This
  choice ensures that $b = -h_x = k_y = -s_{xy}$, the symmetry condition.  The
  functions $h$ and $k$ are unique up to the addition of arbitrary constants, so
  $s$ is unique up to the addition of affine functions of the type $\alpha x +
  \beta y + \gamma$, where $\alpha, \beta, \gamma \in \real$ are arbitrary
  constants.
\end{proof}

We write $T_{s,Q}$ the operator from $\mathcal{M}_{\diiv}$ onto $\mathcal{S}$ mapping $Q$ onto $s$.
Observe that \begin{equation}
  \begin{split}
    T_{s,Q}\,:\, \mathcal{M}_{\diiv} &\longrightarrow \mathcal{S}\\
    Q &\longrightarrow T_{s,Q}[Q]= s
  \end{split}
\end{equation} is a bijection and \begin{equation}\label{ksdjskhdhskdhjd}
  T_{s,Q}^{-1}[s]= R\Hess{s}R^T.
\end{equation} Refer to Figure~\ref{f:param2} for a summary of the relationships between $\sigma$, $Q$
and $s$. Figures~\ref{f:anis4sigma} and~\ref{f:anis4fine} show an example conductivity in each of the
three spaces. \begin{figure}
  \begin{center}
    \xymatrix{
      *-------{} & \txt{\footnotesize \it Physical \\\footnotesize \it conductivity space\\
      \footnotesize \it $\mathcal{M}$ (tensor) \\ \footnotesize \it $\mathcal{M}_{\Iso}$ (scalar)
      } && \txt{\footnotesize \it Divergence-free \\ \footnotesize \it matrix space \\ \footnotesize
      \it $\mathcal{M}_{\operatorname{div}}$} && \txt{\footnotesize \it Convex functions\\
      \footnotesize \it  space  $\mathcal{S}$ }\\
      *-------{\txt{\Large $\Omega$}} & *+[F-,]{\txt{{\huge $\sigma$}}} \ar@{~>} @<1ex>
      [rr]^{\txt{\footnotesize  $Q=\frac{(\nabla F)^T \sigma \nabla F}{\det(\nabla F)}\circ F^{-1}$
      \eref{edgeconddfdsuct} \\ \footnotesize  $F$: $\sigma$-harmonic coordinates \eref{e:harmonic} \\
      non linear}}  & &
      *+[F-,]{\txt{ {\huge $Q$}}} \ar@{~>}@<1ex>[ll]^{\txt{ \footnotesize $\sigma=\sqrt{\det(Q)}\circ
      G^{-1}$ \eref{jwhgdjhwjh}\\ \footnotesize $G$:  $\frac{Q}{\sqrt{\det(Q)}}$-harm. coord.
      \eref{eQharmonic}} } \ar@{->}@<1ex>[rr]^{\txt{\footnotesize $\Hess{s}=R^T Q R$
      \eref{kjshwjhgkjwhg}\\\footnotesize $R$: $\frac{\pi}{2}$-rot. matrix \eref{jshdkjhewkehj}
      \\linear}}  & & *+[F-,]{\txt{{\huge $s$}}}\ar@{->}@<1ex>[ll]^{\txt{\footnotesize $Q=R \Hess{s}
      R^T$ \eref{ksdjskhdhskdhjd}}}  }
  \end{center}
  \caption[Parameterizations of conductivity]{The three parameterizations of
    conductivity, and the spaces to which each belongs. \label{f:param2}}
\end{figure}
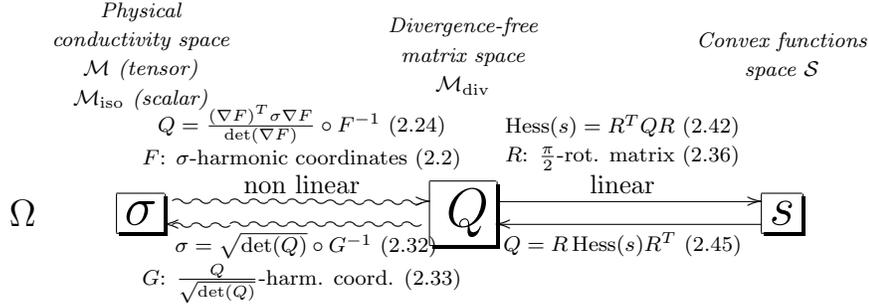 \begin{figure}
  \begin{center}
    \includegraphics[width=0.34\textwidth]{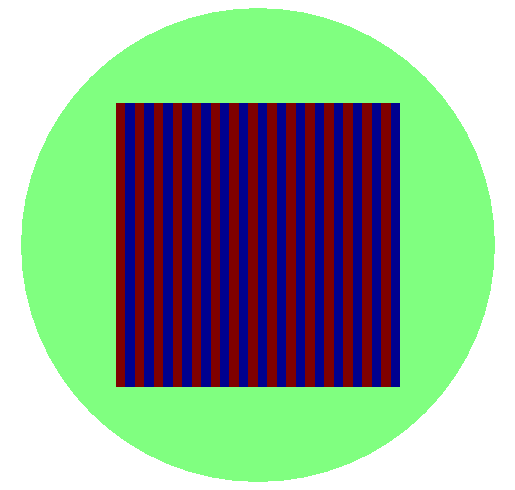}
    \includegraphics[width=0.35\textwidth]{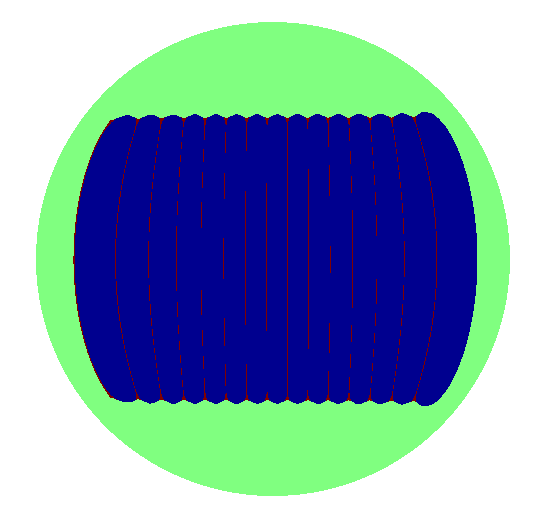}
  \end{center}
  \caption[Example: a laminated conductivity] {The left-hand image shows the original
    scalar conductivity $\sigma(x) = a(x) \operatorname{Id}$.  In blue regions
    a=0.05, in red regions a=1.95, and in green regions, a = 1.0.  The
    right-hand image gives $\sqrt{\det (Q)} = \sigma \circ F^{-1}$, showing how
    harmonic coordinates distort $\sigma(x)$.} \label{f:anis4sigma}
\end{figure} \begin{figure}
  \begin{center}
    \includegraphics[width=0.49\textwidth]{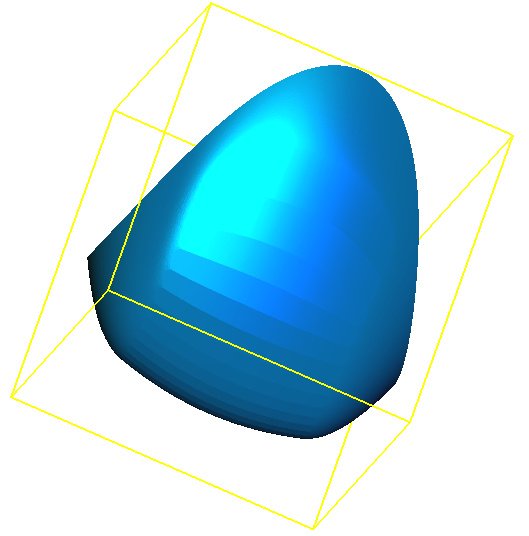}
    \includegraphics[width=0.49\textwidth]{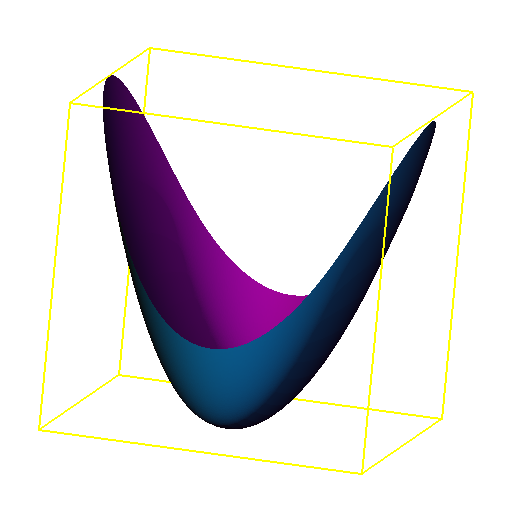}
  \end{center}
  \caption[Fine-scale $s(x)$ for the laminated conductivity] {Two views of the
    fine-scale function $s(x)$ represented as a height field surface for the
    laminated conductivity of Figure~\ref{f:anis4sigma}.  The left-hand view
    shows the fine-scale pattern in $\sigma(x)$, and the right-hand view
    highlights the coarse-scale anisotropy in the
    curvature.} \label{f:anis4fine}
\end{figure}

\section{Discrete geometric homogenization} \label{lksldhslkhsdlkhjej}

We now apply the results of Section~\ref{s:setup} to show that in our framework, homogenization can be
represented either as volume averaging, or as interpolation.  Thus, unlike direct homogenization of
$\sigma \in \mathcal{M}$, homogenization in $\mathcal{M}_{\Div}$ or $\mathcal{S}$ is a linear
operation. Moreover, in this framework, homogenization inherits the semi-group property enjoyed by
volume averaging and interpolation, giving a self-consistency to homogenization in our setting.

\subsection{Homogenization by volume averaging}

The operator $T_{q^h,\sigma}$ defined in~\eref{jklghskjhdgjsh} is a non-linear operator on
$\mathcal{M}$.  However, its restriction to $\mathcal{M}_{\diiv}$, which is a subset of $\mathcal{M}$,
is linear and equivalent to volume averaging as shown by Theorem~\ref{th:volavg}, below.   Using the
notation of Section~\ref{jhskdgkjsghdkj}, we introduce the operator
\begin{equation}\label{jklghskjhddsdsgjsh} \begin{split}
  T_{q^h,Q}\,:\,\mathcal{M}_{\diiv} & \longrightarrow \mathcal{Q}_h\\
  Q & \longrightarrow T_{q^h,\sigma}[Q],
\end{split} \end{equation} where for $Q\in \mathcal{M}_{\diiv}$ and $(i,j)\in \mathcal{E}_h$, one has
\begin{equation}\label{klyiuyoiyu}
  \big(T_{q^h,Q}[Q]\big)_{ij}= -\int_{\Omega}(\nabla \varphi_i)^T Q \nabla \varphi_j.
\end{equation} Observe that $T_{q^h,Q}$ is a volume averaging operator.

\begin{theorem}[Homogenization by volume averaging]
  \label{th:volavg}
  $T_{q^h,Q}$  is a linear volume averaging operator on $\mathcal{M}_{\diiv}$. Moreover:
  \begin{enumerate}
  \item  For $Q\in \mathcal{M}_{\diiv}$  one has
    \begin{equation}\label{kslkdjhjddd}
      T_{q^h,\sigma}[Q]=T_{q^h,Q}[Q].
    \end{equation}
  \item For  $\sigma\in \mathcal{M}$
    \begin{equation}\label{kjdkehsdslkdsejher}
      T_{q^h,\sigma}[\sigma]=T_{q^h,Q} \circ T_{Q,\sigma}[\sigma].
    \end{equation}
  \item Writing $x_j$ the locations of the nodes of $\Omega_h$, for all $l\in \R^2$
    \begin{equation}\label{edgesddsdsdesnzdsdssduct}
      q_{ii}^h (l.x_i)+\sum_{j\sim i} q_{ij}^h (l.x_j)=0.
    \end{equation}
  \end{enumerate}
\end{theorem}

\begin{remark}
  Equation~\eref{kslkdjhjddd} states that $T_{q^h,Q}$ is the restriction of the
  operator $T_{q^h,\sigma}$ to the space of divergence-free matrices
  $\mathcal{M}_{\diiv}$.  It follows from~\eref{kjdkehsdslkdsejher}
  that the homogenization operator $T_{q^h,\sigma}$ is equal to the composition
  of the linear non-injective operator $T_{q^h,Q}$, which acts on
  divergence-free matrices, with the non-linear operator $T_{Q,\sigma}$, which
  projects into the space of divergence-free matrices. Observe also that
  $T_{Q,\sigma}$ is injective as an operator from $\mathcal{M}_{\Iso}$, the
  space of scalar conductivities, onto $\mathcal{M}_{\diiv}$.
\end{remark} \begin{remark}
  Equation~\eref{edgesddsdsdesnzdsdssduct} is essentially stating that $q^h$ is
  divergence free at a discrete level, see~\cite[Section 2.1]{dona2008} for
  details.
\end{remark} \begin{proof}
  Using the change of coordinates $y=F(x)$ we obtain that
  \begin{equation}\label{edgeconduct1}
    \int_{\Omega} (\nabla (\varphi_i \circ F))^T \sigma(x) \nabla (\varphi_j \circ
    F)\,dx=\int_{\Omega} (\nabla \varphi_i)^T Q \nabla \varphi_j
  \end{equation}
  which implies~\eref{kjdkehsdslkdsejher}. One obtains
  equation~\eref{kslkdjhjddd} by observing that since $Q$ is divergence-free its
  associated harmonic coordinates are just linear functions and
  $T_{\sigma,Q}[Q]=Q$.  Since $Q$ is divergence-free, we have, for constant
  $l\in\real^2$,
  \begin{equation}\label{edgesddsnzdsdssduct}
    \int_{\Omega} (\nabla \varphi_i )^T Q(x).l\,dx = 0.
  \end{equation}
  Now, set $\mathcal{V}_h$ the set off all nodes in the triangulation $\Omega_h$
  and set $x_j$ the location of node $j \in \mathcal{V}_h$.  The function $z(x)
  := \sum_{j\in \mathcal{V}_h} x_j\varphi_j(x)$ is the identity map on
  $\Omega_h$, so we can write $l=\nabla \left(\sum_{j\in\mathcal{V}_h}
  (l.x_j)\varphi_j(x) \right)$.  Combining this with~\eref{edgesddsnzdsdssduct}
  gives~\eref{edgesddsdsdesnzdsdssduct}.
\end{proof}

\subsection{Homogenization by linear interpolation}

Write $T_{s^h,s}$ the linear interpolation operator over $\Omega_h$.  Hence for $s\in \mathcal{S}$ and
$s^h:=T_{s^h,s}[s]$, we have, for $x\in \Omega$, \begin{equation}\label{klalkshdklhj}
  s^h(x)=\sum_{i} s(x_i) \varphi_i(x),
\end{equation} where the sum in~\eref{klalkshdklhj} is taken over all nodes of $\Omega_h$ and $x_i$ is
the location of node $i$. Write $\mathcal{S}_h$ the space of linear interpolations of elements of
$\mathcal{S}$ on $\Omega_h$. Hence, \begin{equation}\label{klkglkhg}
  \begin{split}
    T_{s^h,s}\,:\,\mathcal{S} & \longrightarrow \mathcal{S}_h\\
    s & \longrightarrow T_{s^h,s}[s]:= \sum_{i} s(x_i) \varphi_i(x).
  \end{split}
\end{equation}

For $(i,j)\in \mathcal{E}_h$ write $\delta_{(i,j)}(x)$ the uniform Lebesgue (Dirac) measure on the
edge $(i,j)$ (as a subset of $\R^2$). Let $R$ be the rotation matrix already introduced in
\eref{jshdkjhewkehj}. For $s^h \in \mathcal{S}_h$ observe that $R\Hess{s^h}R^T$ is a Dirac measure on
edges of $\Omega_h$. For $(i,j)\in \mathcal{E}_h$ define $\big(T_{q^h,s^h}[s^h]\big)_{ij}$ as the
curvature of $s^h$ in the direction orthogonal to the edge $(i,j)$, hence
\begin{equation}\label{jjhgkjgjkhhjk}
  \begin{split}
    T_{q^h,s^h}\,:\, \mathcal{S}_h &\longrightarrow \mathcal{Q}_h\\
    s^h &\longrightarrow T_{q^h,s^h}[s^h]
  \end{split}
\end{equation} with \begin{equation}\label{lkglkjhgkjg}
  \sum_{(i,j)\in \mathcal{E}_h} \big(T_{q^h,s^h}[s^h]\big)_{ij} \delta_{(i,j)}= R\Hess{s^h}R^T.
\end{equation}

Write $s_i=s(x_i)$ where $s_i$ is the location of the node $i$. Referring to Figure~\ref{f:cotans},
\begin{figure}
 \begin{center}
\begin{picture}(0,0)%
\includegraphics{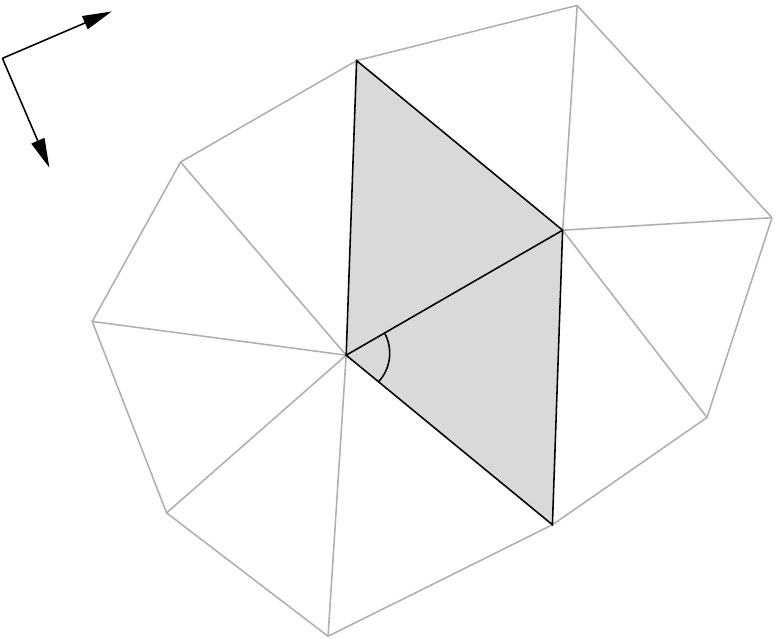}%
\end{picture}%
\setlength{\unitlength}{3947sp}%
\begingroup\makeatletter\ifx\SetFigFont\undefined%
\gdef\SetFigFont#1#2#3#4#5{%
  \reset@font\fontsize{#1}{#2pt}%
  \fontfamily{#3}\fontseries{#4}\fontshape{#5}%
  \selectfont}%
\fi\endgroup%
\begin{picture}(3718,3052)(1339,-3522)
\put(3976,-3136){\makebox(0,0)[lb]{\smash{{\SetFigFont{10}{12.0}{\familydefault}{\mddefault}{\updefault}{\color[rgb]{0,0,0}$k$}%
}}}}
\put(3001,-661){\makebox(0,0)[lb]{\smash{{\SetFigFont{10}{12.0}{\familydefault}{\mddefault}{\updefault}{\color[rgb]{0,0,0}$l$}%
}}}}
\put(1876,-661){\makebox(0,0)[lb]{\smash{{\SetFigFont{10}{12.0}{\familydefault}{\mddefault}{\updefault}{\color[rgb]{0,0,0}$\eta$}%
}}}}
\put(1576,-1336){\makebox(0,0)[lb]{\smash{{\SetFigFont{10}{12.0}{\familydefault}{\mddefault}{\updefault}{\color[rgb]{0,0,0}$\xi$}%
}}}}
\put(4051,-1486){\makebox(0,0)[lb]{\smash{{\SetFigFont{10}{12.0}{\familydefault}{\mddefault}{\updefault}{\color[rgb]{0,0,0}$i$}%
}}}}
\put(2926,-2386){\makebox(0,0)[lb]{\smash{{\SetFigFont{10}{12.0}{\familydefault}{\mddefault}{\updefault}{\color[rgb]{0,0,0}$j$}%
}}}}
\put(3346,-1786){\makebox(0,0)[lb]{\smash{{\SetFigFont{10}{12.0}{\familydefault}{\mddefault}{\updefault}{\color[rgb]{0,0,0}$e_{ij}$}%
}}}}
\put(3634,-2399){\makebox(0,0)[lb]{\smash{{\SetFigFont{10}{12.0}{\familydefault}{\mddefault}{\updefault}{\color[rgb]{0,0,0}$t_{ijk}$}%
}}}}
\put(3226,-2204){\makebox(0,0)[lb]{\smash{{\SetFigFont{10}{12.0}{\familydefault}{\mddefault}{\updefault}{\color[rgb]{0,0,0}$\theta_{ijk}$}%
}}}}
\end{picture}%
 \end{center}
 \caption{Notation for computing $q_{ij}^h$ from $s_i,s_j,s_k$ and $s_l$. \label{f:cotans}}
\end{figure}
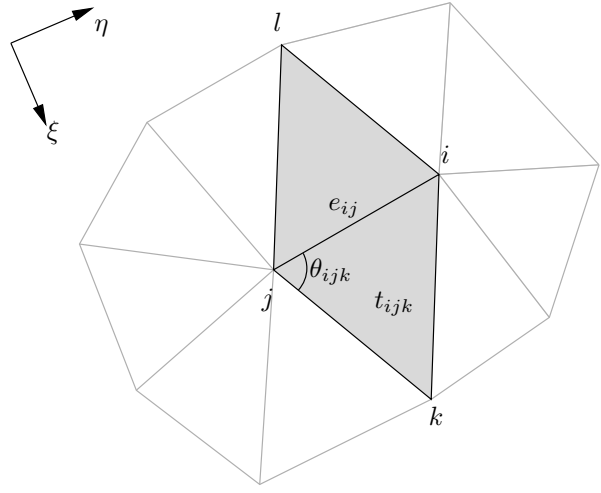 $T_{q^h,s^h}[s^h]$ defined as the trace of the $R$-rotated Hessian of $s^h$ on the edge
$(i,j)$ is \begin{equation} \begin{split}
  \big(T_{q^h,s^h}[s^h]\big)_{ij} =& -\dfrac{1}{\abs{e_{ij}}^2}\left(\cot\theta_{ijk} +
    \cot\theta_{ijl}\right) s_i \\
    & -\dfrac{1}{\abs{e_{ij}}^2}\left(\cot\theta_{jik} +
    \cot\theta_{jil}\right) s_j \\
    & + \dfrac{1}{2 \abs{t_{ijk}}} s_k
    + \dfrac{1}{2 \abs{t_{ijl}}} s_l,
    \label{e:hinge}
\end{split} \end{equation} while diagonal elements $\big(T_{q^h,s^h}[s^h]\big)_{ii}$  are expressed as
\[\big(T_{q^h,s^h}[s^h]\big)_{ii} = -\sum_{j\sim i} \big(T_{q^h,s^h}[s^h]\big)_{ij},\] where $j\sim i$
is the set of vertices distinct from $i$ and sharing an edge with vertex $i$, $\abs{e_{ij}}$ is the
length of edge $(i,j)$, $\abs{t_{ijk}}$ is the area of the triangle with vertices $(i,j,k)$, and
$\theta_{ijk}$ is the interior angle of triangle $t_{ijk}$ at vertex $j$ (see Figure~\ref{f:cotans}).

Note that~\eqref{e:hinge} is valid only for interior edges. Because of our choice to interpolate
$s(x)$ by piecewise linear functions, we have concentrated all of the curvature of $s(x)$ on the edges
of the mesh, and we need a complete hinge, an edge with two incident triangles, in order to
approximate this curvature.  Without values for $s(x)$ outside of $\Omega$ and hence exterior to the
mesh, we do not have a complete hinge on boundary edges. This will become important where we apply our
method to solve the inverse homogenization problem in EIT.  However, for the homogenization problem,
our homogeneous boundary conditions make irrelevant the values of $q_{ij}^h$ on boundary edges.

$T_{q^h,s^h}$ defined through~\eqref{e:hinge} has several nice properties.  For example, direct
calculation shows that $T_{q^h,s^h}[s^h]$ computed using~\eqref{e:hinge} is divergence-free in the
discrete sense given by~\eref{edgesddsdsdesnzdsdssduct} for any values $s_i$.  This fact allows us to
parameterize the space of edge conductivities $q^h$ satisfying the discrete divergence-free
condition~\eref{edgesddsdsdesnzdsdssduct} by linear interpolations of convex functions.
\begin{proposition}[Discrete divergence-free parameterization of conductivity] $T_{q^h,s^h}$ defined
using~\eqref{e:hinge} has the following properties:
  \begin{enumerate}
  \item Affine functions are exactly the nullspace of $T_{q^h,s^h}$; in
    particular, $q^h:=T_{q^h,s^h}[s^h]$ is divergence-free in the discrete sense
    of~\eref{edgesddsdsdesnzdsdssduct}.
  \item The dimension of the range of $T_{q^h,s^h}$ is equal to the
    number of edges in the triangulation, minus the discrete divergence-free
    constraints~\eref{edgesddsdsdesnzdsdssduct}.
  \item $T_{q^h,s^h}$ defines a bijection from $\mathcal{S}_h$ onto
    $\mathcal{Q}_h$ and for $s^h \in \mathcal{S}_h$
    \begin{equation}\label{kljdhlkjhkhedj}
      T_{s^h,q^h}^{-1}[s^h]=T_{s,Q}^{-1}[s^h].
    \end{equation}
  \end{enumerate}
\end{proposition}

\begin{proof}
  These properties can be confirmed in both volume-averaged and interpolation
  spaces:
  \begin{enumerate}
  \item The first property can be verified directly from the hinge
    formula~\eqref{e:hinge}.

  \item For the Dirichlet problem in finite elements, the number of degrees of
    freedom in a stiffness matrix which is not necessarily divergence-free
    equals the number of interior edges on the triangle mesh.  The
    divergence-free constraint imposes two constraints---one for each of the
    $x-$ and $y-$directions---at each interior vertex such that the left term
    of~\eqref{jashgajsggsjhw}, namely $\sum_{j\sim i} q_{ij}^h(v_i-v_j)$, is
    zero for affine functions.  Thus, the divergence-free stiffness matrix has
    \begin{equation}
      E_I - 2 V_I
    \end{equation}
    degrees of freedom, where $E_I$ is the number of interior edges, and $V_I$
    is the number of interior vertices.

    The piecewise linear interpolation of $s(x)$ has $V-3$ degrees of freedom,
    where there are $V$ vertices in the mesh.  The restriction of 3 degrees of
    freedom corresponds to the arbitrary addition of affine functions to $s(x)$
    bearing no change to $Q$.

    Our triangulation $\Omega_h$ tessallates our domain $\Omega$, and so
    $\Omega_h$ is a simply connected domain of trivial topology.  For this
    topology, it can be shown that the number of edges $E$ is
    \begin{equation}
      E = 2V + V_I - 3.
    \end{equation}
    recalling that the number of boundary edges equals the number of boundary
    vertices, we have
    \begin{equation}
      E_I - 2 V_I = V - 3,
    \end{equation}
    In fact $s(x)$ and $Q(x)$ represented on the same mesh have the same degrees
    of freedom when $Q(x)$ is divergence-free.
  \end{enumerate}
  \item This property can be easily checked from the previous ones.
\end{proof}

\begin{theorem}
  $T_{s^h,s}$, a linear interpolation operator on $\mathcal{S}$, has the
  following properties:
  \begin{enumerate}
  \item  For $Q\in \mathcal{M}_{\diiv}$,
    \begin{equation}\label{kslkdjhdjddd}
      T_{q^h,Q}[Q]=T_{q^h,s^h}\circ T_{s^h,s} \circ T_{s,Q}[Q].
    \end{equation}
  \item For $\sigma\in \mathcal{M}$,
    \begin{equation}\label{kjdkehsdsdslkdsejher}
      T_{q^h,\sigma}[\sigma]=T_{q^h,s^h}\circ T_{s^h,s} \circ T_{s,Q}  \circ
      T_{Q,\sigma}[\sigma].
    \end{equation}
  \end{enumerate}
\end{theorem}

\begin{remark}
  It follows from equations~\eref{kslkdjhdjddd} and~\eref{kjdkehsdsdslkdsejher}
  that homogenization is a linear interpolation operator acting on convex
  functions.  Observe that $T_{q^h,s^h}$, $T_{s^h,s}$ and $T_{s,Q}$ are all
  linear operators.  Hence, the non-linearity of the homogenization operator is
  confined to the non-linear projection operator $T_{Q,\sigma}$
  in~\eref{kjdkehsdsdslkdsejher} whereas if $\sigma$ is scalar its
  non-injectivity is confined to the linear interpolation operator $T_{s^h,s}$.
  Equation~\eref{kljdhlkjhkhedj} is understood in terms of measures on edges of
  $\Omega_h$ and implies that the bijective operator mapping $q^h$ onto $s^h$ is
  a restriction of the bijective operator mapping $Q$ onto $s$ to the spaces
  $\mathcal{Q}_h$ and $\mathcal{S}_h$.
\end{remark} \begin{remark}
  Provided that the $s_i$ interpolate a convex function $s(x)$, the $q_{ij}^h =
  \big(T_{q^h,s^h}[s^h]\big)_{ij}$ form a positive semi-definite stiffness
  matrix even if not all $q_{ij}^h$ are strictly positive.  We discuss this
  further in the next section, where we show that even with this flexibility in
  the sign of the $q^h_{ij}$, it is always possible to triangulate a domain such
  that $q^h_{ij}>0$.
\end{remark}

\begin{proof}
  Define a coordinate system $\xi$-$\eta$ such that edge $ij$ is parallel to the
  $\eta$-axis as illustrated in
  Figure~\ref{f:cotans}. Using~\eref{edgesddsnzdduct} to rewrite $T_{q^h,Q}\circ
  T_{Q,s}[s]$ in this rotated coordinate system yields
  \begin{equation}
    q_{ij}^h = -\int_{\Omega} (\Grad \varphi_i)^T
    \begin{pmatrix}
      s_{\eta\eta} & -s_{\xi\eta} \\
      -s_{\xi\eta} & s_{\xi\xi}
    \end{pmatrix}
    \Grad \varphi_j.
    \label{e:stiffness2}
  \end{equation}
  A change of variables confirms that integral~\eqref{e:stiffness2} is invariant
  under rotation and translation.  We abuse notation in that the second
  derivatives are understood here in the sense of measures: we are about to
  interpolate $s(x)$ by piecewise linear functions, which do not have pointwise
  second derivatives everywhere. We are concerned with the values of $s(x)$
  interpolated at $i,j,k,$ and $l$, as these are associated to only the
  corresponding hat basis functions sharing support with those at $i$ and $j$.
  The second derivatives of $\varphi$ are non-zero only on edges, and due to the
  support of the gradients of the $\varphi$, contributions of the second
  derivatives at edges $e_{ik}$, $e_{jk}$, $e_{il}$, and $e_{jl}$ are also zero.
  Finally, the $\partial_{\xi\eta} \varphi$ and $\partial_{\eta\eta} \varphi$ are
  zero along $ij$, so the only contributions of $s(x)$ to $T_{q^h,Q}\circ
  T_{Q,s}[s]$ defined through the integral are its second derivatives with
  respect to $\xi$ along edge $e_{ij}$.  The contributions of four integrals
  remain, and by symmetry, we have only two integrals to compute.  Noting that
  the singularities in the first and second derivatives are not coincident, from
  direct computation of the gradients of the basis functions and integration by
  parts we have
  \begin{align}
    \int_{t_{ijk} \cup t_{ijl}} \partial_{\eta} \varphi_i \partial_{\xi\xi}
    \varphi_i \partial_{\eta} \varphi_j &=
    \frac{1}{\abs{e_{ij}}^2} \left( \cot\theta_{ijk} + \cot\theta_{ijl}
    \right),
    \\
    \int_{t_{ijk} \cup t_{ijl}} \partial_{\eta} \varphi_i \partial_{\xi\xi}
    \varphi_k \partial_{\eta} \varphi_j &=
    -\frac{1}{2 \abs{t_{ijk}}},
  \end{align}
  where $\abs{e_{ij}}$ is the length of the edge with vertices $(i,j)$, and
  $\abs{t_{ijk}}$ is the area of the triangle with vertices
  $(i,j,k)$. $\theta_{ijk}$ is the interior angle of triangle $ijk$ at vertex
  $j$ (see Figure~\ref{f:cotans}).  The only contribution to these integrals is
  in the neighborhood of edge $e_{ij}$. Combining these results, we have that
  the elements of the stiffness matrix are given by formula~\eqref{e:hinge}.
\end{proof}

\begin{figure} \begin{center} \xymatrix{ *-------{} & \txt{\footnotesize \it Physical \\\footnotesize
\it conductivity space\\ \footnotesize \it $\mathcal{M}$ (tensor) \\ \footnotesize \it
$\mathcal{M}_{\Iso}$ (scalar) } && \txt{\footnotesize \it Divergence-free \\ \footnotesize \it matrix
space \\ \footnotesize $\mathcal{M}_{\operatorname{div}}$} && \txt{\footnotesize \it Convex
functions\\ \footnotesize \it space  $\mathcal{S}$ }\\
 *-------{\txt{\Large $\Omega$}} & *+[F-,]{\txt{{\huge $\sigma$}}} \ar@{~>} @<1ex>
 [rr]^{\txt{\footnotesize  $\sigma$-harm. coord. \eref{e:harmonic},\eref{edgeconddfdsuct} }}
 \ar@{~>}[rrddd]_{\txt{\footnotesize  $q^h_{ij}=\int_\Omega \nabla (\varphi_i\circ F)\sigma \nabla
 (\varphi_j\circ F)$\\ \eref{edgeconduct}, \eref{jashgajsggsjhw}}}  & &
*+[F-,]{\txt{ {\huge $Q$}}} \ar@{~>}@<1ex>[ll]^{\txt{ \footnotesize $\frac{Q}{\sqrt{\det(Q)}}$-har.
coor. \eref{eQharmonic}, \eref{jwhgdjhwjh}} } \ar@{->}@<1ex>[rr]^{\txt{\footnotesize $\Hess{s}=R^T Q
R$ \eref{kjshwjhgkjwhg}}} \ar@{->}[ddd]|-{\txt{\footnotesize  $q^h_{ij}=\int_\Omega \nabla \varphi_i
Q\nabla \varphi_j$\\ \eref{klyiuyoiyu}}} & & *+[F-,]{\txt{{\huge
$s$}}}\ar@{->}@<1ex>[ll]^{\txt{\footnotesize $Q=R \Hess{s} R^T$
\eref{ksdjskhdhskdhjd}}}\ar@{->}[ddd]|-{\txt{\footnotesize  $s^h=\sum_{i}s(x_i)\varphi_i(x)$\\
\eref{klkglkhg}}}\\ *-------{} & && &&\\ *-------{} & && &&\\ *-------{\txt{\Large $\Omega_h$}} & &&
*+[F-,]{\txt{ {\huge $q^h$}}} \ar@{->}@<1ex>[rr]^{\txt{\footnotesize $\Hess{s^h}=R^T q^h R$
\eref{lkglkjhgkjg}}} && *+[F-,]{\txt{ {\huge $s^h$}}}\ar@{->}@<1ex>[ll]^{\txt{\footnotesize $q^h=R
\Hess{s^h} R^T$ \eref{lkglkjhgkjg}}}\\ *-------{} & && \txt{\it effective edge\\ \it
conductivities}&&\\} \end{center} \caption[Summary of homognization spaces]{Summary of discrete
homogenization,
  showing the relationships between the discrete spaces which approximate the
  spaces introduced in Section~\ref{s:setup}. \label{f:discretesummary}}
\end{figure}
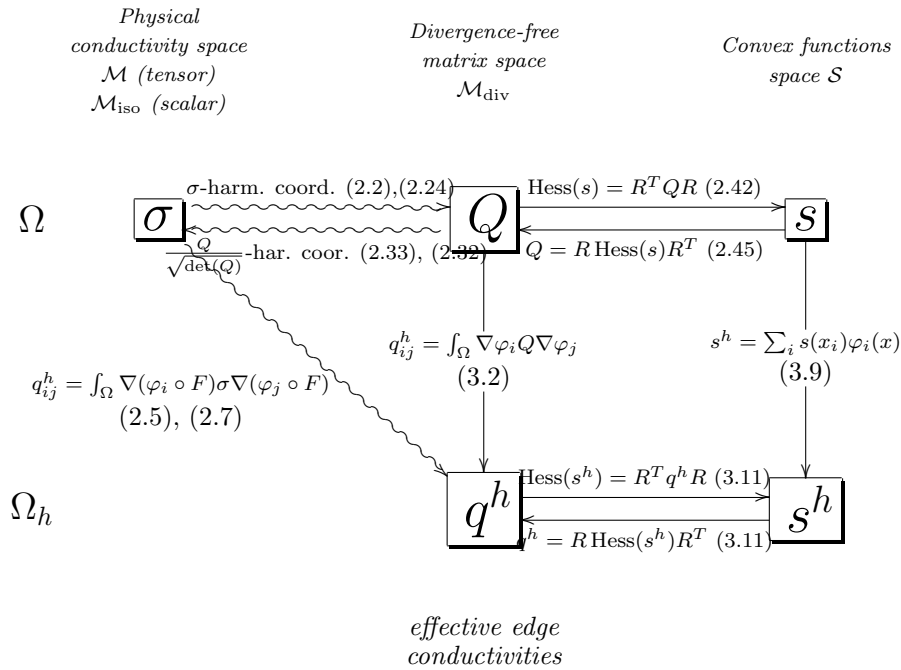

We refer to Figure~\ref{f:discretesummary} for a summary of the results of this subsection.

\subsection{Semi-group properties in geometric homogenization} \label{ksjdlkshdhksjdh}

Consider, for example, three approximation scales $0<h_1<h_2<h_3$.  We now show that homogenization
from $h_1$ to $h_3$ is identical to homogenization from $h_1$ to $h_2$, then from $h_2$ to $h_3$.  We
identify this as a semi-group property.

Let $\Omega_{h_{C}}$ be a coarse triangulation of $\Omega$, and $\Omega_{h_{F}}$ be a finer,
sub-triangulation of $\Omega_{h_{F}}$. Let $\varphi_i^C, \varphi_i^F$ be the piecewise linear nodal
basis functions centered on the interior nodes of $\Omega_{h_{C}}$ and $\Omega_{h_{F}}$. Observe that
for each interior node of the coarse triangulation $i\in \mathcal{N}_{h_C}$, $\varphi_i^C$ can be
written as a linear combination of $\varphi_k^F$, we write $\phi_{ik}$ the coefficients of that linear
combination. Hence \begin{equation}
  \varphi_i^C=\sum_{k\in \mathcal{N}_{h_F}} \phi_{ik} \varphi_k^F.
\end{equation} Define $T_{q^{h_C},q^{h_F}}$ as the operator mapping the effective conductivities of
the edges of fine triangulation onto the effective conductivities of the edges of the coarse
triangulation. Hence \begin{equation}
  \begin{split}
    T_{q^{h_C},q^{h_F}}\,:\, \mathcal{Q}_{h_F} &\longrightarrow \mathcal{Q}_{h_C}\\
    q^{h_F} &\longrightarrow T_{q^{h_C},q^{h_F}}[q^{h_F}]
  \end{split}
\end{equation} with, for $(i,j)\in \mathcal{E}_{h_C}$, \begin{equation}
  \big(T_{q^{h_C},q^{h_F}}[q^{h_F}]\big)_{ij}= \sum_{l,k\in \mathcal{N}_{h_F}\,:\, (l,k)\in
  \mathcal{E}_{h_F}}
  \phi_{ik}\phi_{jl} q^{h_F}_{kl}.
\end{equation}

Let $T_{s^{h_C},s^{h_F}}$ be the linear interpolation operator mapping piecewise linear functions on
$\Omega_{h_F}$ onto piecewise linear functions on $\Omega_{h_C}$.  Hence \begin{equation}
\begin{split}
  T_{s^{h_C},s^{h_F}}\,:\,\mathcal{S}_{h_F} & \longrightarrow \mathcal{S}_{h_C}\\
  s^{h_F}& \longrightarrow T_{s^{h_C},s^{h_F}}[s^{h_F}],
\end{split} \end{equation} and as in \eref{klalkshdklhj}, we have for $x\in \Omega$
\begin{equation}\label{klalkshdxcdklhj}
  T_{s^{h_C},s^{h_F}}[s^{h_F}](x)=\sum_{i\in \mathcal{N}_{h_C}} s^{h_F}(x_i) \varphi_i^C(x).
\end{equation}

\begin{figure} \begin{center} \xymatrix{ *-------{} & \txt{\footnotesize \it Physical
\\\footnotesize \it conductivity space} && \txt{\footnotesize \it Divergence-free \\ \footnotesize \it
matrix space } && \txt{\footnotesize \it Convex functions\\ \footnotesize \it space }\\
 *-------{\txt{\Large $\Omega$}} & *+[F-,]{\txt{{\huge $\sigma$}}} \ar@{~>} @<1ex> [rr] \ar@{~>}[rrd]
 \ar@{~>}[rrdd] \ar@{~>}[rrddd] & &
*+[F-,]{\txt{ {\huge $Q$}}} \ar@{~>}@<1ex>[ll]\ar@{->}@<1ex>[rr] \ar@{->}[d] & & *+[F-,]{\txt{{\huge
$s$}}}\ar@{->}@<1ex>[ll]\ar@{->}[d]\\ *-------{\txt{\Large $\Omega_{h_1}$}} & && *+[F-,]{\txt{ {\huge
$q^{h_1}$}}} \ar@{->}@<1ex>[rr]\ar@{->}[d] && *+[F-,]{\txt{ {\huge
$s^{h_1}$}}}\ar@{->}@<1ex>[ll]\ar@{->}[d] \\ *-------{\txt{\Large $\Omega_{h_2}$}} & && *+[F-,]{\txt{
{\huge $q^{h_2}$}}} \ar@{->}@<1ex>[rr]\ar@{->}[d] && *+[F-,]{\txt{ {\huge
$s^{h_2}$}}}\ar@{->}@<1ex>[ll]\ar@{->}[d] \\ *-------{\txt{\Large $\Omega_{h_3}$}} & && *+[F-,]{\txt{
{\huge $q^{h_3}$}}} \ar@{->}@<1ex>[rr] && *+[F-,]{\txt{ {\huge $s^{h_3}$}}}\ar@{->}@<1ex>[ll]
  }
\end{center} \caption[Sequence of scales in the semi-group properties]{Discrete geometric
  homogenization showing the sequence of scales referred to by the
  semi-group properties.  \label{f:semigroup}}
\end{figure}
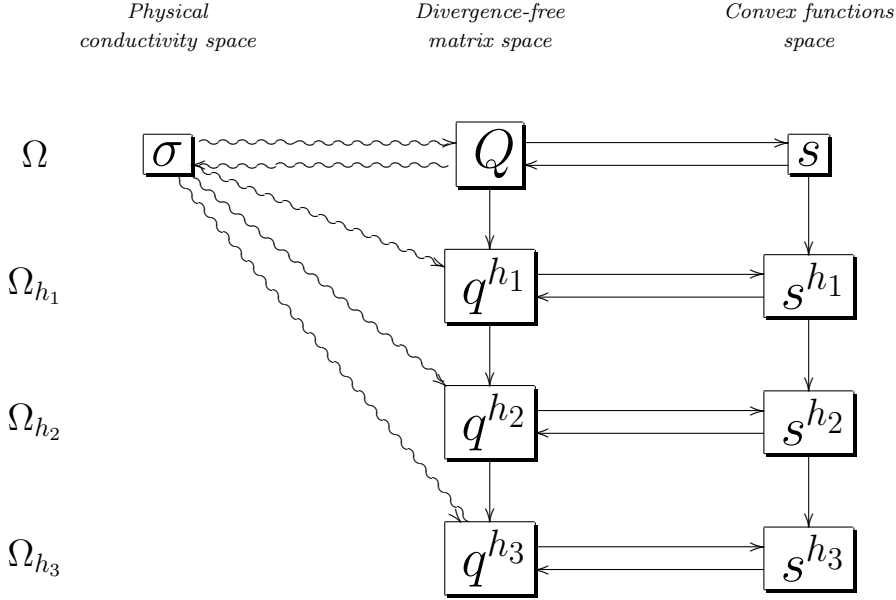

\begin{theorem}[Semi-group properties in geometric homogenization]
  \label{klwekjwhge}
  The linear operators $T_{q^{h_C},q^{h_F}}$ and $T_{s^{h_C},s^{h_F}}$ satisfy the
  following properties:
  \begin{enumerate}
  \item $T_{s^{h_C},s^{h_F}}$ is the restriction of the interpolation operator
    $T_{s^{h_C},s}$ to piecewise linear functions on $\Omega_{h_F}$.  That is,
    for $s^{h_F}\in \mathcal{S}_{h_F}$
    \begin{equation}\label{kslkdjhsdssdsddsdsdsdsddjddd}
      T_{s^{h_C},s^{h_F}}[s^{h_F}]=T_{s^{h_C},s}[s^{h_F}].
    \end{equation}
  \item For $Q\in \mathcal{M}_{\diiv}$
    \begin{equation}\label{kjdkehsddioi8sdcdslkdsejher}
      T_{q^{h_C},Q}[Q]=T_{q^{h_C},q^{h_F}}\circ T_{q^{h_F},Q}[Q].
    \end{equation}
  \item For $s\in \mathcal{S}$
    \begin{equation}\label{kjdkehsddsdcdslkdsejher}
      T_{s^{h_C},s}[s]=T_{s^{h_C},s^{h_F}}\circ T_{s^{h_F},s}[s].
    \end{equation}
  \item For $\sigma\in \mathcal{M}$
    \begin{equation}\label{kjdkslkdsejher}
      T_{q^{h_C},\sigma}[\sigma]=T_{q^{h_C},q^{h_F}}\circ T_{q^{h_F},\sigma}[\sigma].
    \end{equation}
  \item For $q^{h_F} \in \mathcal{Q}_{h_F}$
    \begin{equation}\label{kjdxkehsdher}
      T_{q^{h_C},q^{h_F}}[q^{h_F}]=T_{q^{h_C},s^{h_C}}\circ T_{s^{h_C},s^{h_F}} \circ
      T_{q^{h_F},s^{h_F}}^{-1}[q^{h_F}].
    \end{equation}
  \item For $h_1<h_2<h_3$
    \begin{equation}\label{daekjk}
      T_{s^{h_3},s^{h_1}}=T_{s^{h_3},s^{h_2}}\circ T_{s^{h_2},s^{h_1}}.
    \end{equation}
  \item For $h_1<h_2<h_3$
    \begin{equation}\label{dae}
      T_{q^{h_3},q^{h_1}}=T_{q^{h_3},q^{h_2}}\circ T_{q^{h_2},q^{h_1}}.
    \end{equation}
  \end{enumerate}
\end{theorem} \begin{remark}
  As we will see below, if the triangulation $\Omega_{h_F}$ is not chosen
  properly $s^{h_F}=T_{s^{h_F},s}[s]$ may not be convex. In that situation
  $T_{s^{h_C},s}$ in~\eref{kslkdjhsdssdsddsdsdsdsddjddd}, when acting on
  $s^{h_F}$, has to be interpreted as a linear interpolation operator over
  $\Omega_{h_C}$ acting on continuous functions of $\Omega$. We will show in the
  next section how to choose the triangulation $\Omega_{h_F}$ (resp.\ $\Omega_{h_C}$)
  to ensure the convexity of $s^{h_F}$ (resp.\ $s^{h_C}$).
\end{remark} \begin{remark}
  The semi-group properties obtained in Theorem~\ref{klwekjwhge} are essential
  to the self-consistency of any homogenization theory. The fact that
  homogenizing directly from scale $h_1$ to scale $h_3$ is equivalent to
  homogenizing from scale $h_1$ to scale $h_2$ then from $h_2$ onto $h_3$ is a
  property that is in general not satisfied by most numerical homogenization
  methods found in the literature when applied to PDEs with arbitrary
  coefficients, such as non-periodic or non-ergodic conductivities.
  Figure~\ref{f:semigroup} illustrates the sequence of scales referred to by
  these semi-group properties.
\end{remark} 

\section{Optimal meshes based on convex functions} \label{klsjsldkhklsj} \label{s:triangulate}

In this section, we use the convex function parameterization $s\in\mathcal{S}$ to construct
triangulations of $\Omega$ which give matrices approximating the elliptic operator with optimal
conditioning.  In particular, we show that we can triangulate a given set of vertices such that the
off-diagonal elements of the stiffness matrix $q^h_{ij}$ are always non-positive.  In turn, this
minimizes the radii of the Gershgorin disks containing the eigenvalues of $q^h_{ij}$.  The argument
directly uses the geometry of $s(x)$, constructing the triangulation from the convex hull of points
projected up to $s(x)$.  We show that this procedure, a general case of the convex hull projection
method for producing the Delaunay triangulation from a paraboloid, produces a {\em weighted} Delaunay
triangulation.  That is, we provide a geometric interpretation of the weighted Delaunay triangulation
as well as an efficient method for producing optimal $Q$-adapted meshes.

 We remind the reader that throughout this paper,
our triangulation $\Omega_h$ is a tessellation of the compact and simply-connected domain $\Omega$,
and, as such, is itself simply connected and of trivial topology. Also, since $\Omega\subset \real^2$,
we shall identify the arguments of scalar functions as in $s(x), x\in\real^2$, or $s(x,y),
x,y\in\real$ interchangeably without further comment.

\subsection{Construction of positive Dirichlet weights}

The constant $C$ in~\eref{kjwhlkwhel} can be minimized by choosing the triangulation in a manner that
ensures the positivity of the effective edge conductivities $q^h_{ij}$. The reason behind this
observation lies in the fact that the discrete Dirichlet energy associated to the homogenized problem
\eref{jashgajsggsjhw} is \begin{equation}
  E_Q(u) = \onehalf \sum_{i\sim j} q^h_{ij} \left(u_i - u_j\right)^2
  \label{e:discreteDirichlet}
\end{equation} where $i\sim j$ are the edges of the triangulation, and $u_i$ interpolate $u(x)$ at
vertices.

We now show that for $Q$ divergence-free, we can use a parameterization $s(x)$ to build a
triangulation such that $q^h_{ij} \geq 0$.  $q^h_{ij}$, identified here as {\em Dirichlet weights} are
typically computed as elements of the stiffness matrix, where $Q$ is known exactly.  In this paper, we
have introduced the parameterization $s(x)$ for divergence-free conductivities, and if we interpolate
$s(x)$ by piecewise-linear functions, $q^h_{ij}$ is given by the hinge formula~\eqref{e:hinge}.

In the special case where $Q$ is the identity, it is well-known~\cite{pp-cdms-93} that
\begin{equation}
  q^h_{ij} = \onehalf \left(\cotan \theta_{ikj} + \cotan \theta_{ilj}\right),
\end{equation} and in such case, all $q^h_{ij} \geq 0$ when the vertices are connected by a Delaunay
triangulation.  Moreover, the Delaunay triangulation can be constructed geometrically.  Starting with
a set of vertices, the vertices are projected to the surface of any regular parabaloid
\begin{equation}
  p(x,y) = a \left(x^2 + y^2\right),
\end{equation} where $a>0$ is constant.  The convex hull of these points forms a triangulation over
the surface of $p(x,y)$, and the projection of this triangulation back to the $xy$-plane is Delaunay.
See~\cite{rour1998}, for example.  Our observation is that the correspondance \begin{equation}
  Q = \mbox{identity} \Rightarrow s(x,y) = \onehalf \left(x^2 + y^2\right),
\end{equation} can be extended to all positive-definite and divergence-free $Q$.  By constructing our
triangulation as the projection of the convex hull of a set of points projected on to {\em any} convex
$s(x,y)$, we have the following:

\begin{theorem} \label{th:posq}
  Given a set of points $\calV$, there exists a triangulation of those points
  such that all $q^h_{ij}\geq 0$.  We refer to this triangulation as a
  {\em $Q$-adapted triangulation}.  If there exists no edge for which
  $q^h_{ij}=0$, this triangulation is unique.
\end{theorem}

\begin{remark}
  The points $\calV \supset \calV_h$, where $\calV_h$ is the set of nodes in the
  resulting triangulation $\Omega_h$.  That is, some points in $\calV$ may be
  decimated by the triangulation.  See also the remark following
  Proposition~\ref{p:weightedDelaunay}.
\end{remark} \begin{remark}
  $q^h_{ij}>0$ does not hold for arbitrary triangulations, as each
  triangulation, according to its connectivity, admits a different set of
  piecewise linear basis functions $\myphi_i$.
\end{remark} \begin{remark}
  While $s(x,y)$ may be convex, an arbitrary piecewise linear interpolation may
  not be.  Figure~\ref{f:flips} illustrates two interpolations of $s(x,y)$, one
  of which gives a $q^h_{ij}> 0$, and the other of which does not.  Moreover, we
  note that as long as the function $s(x,y)$ giving our interpolants $s_i$ is
  convex, the discrete Dirichlet operator is positive semi-definite, even if
  some individual elements $q^h_{ij}<0$.  Figure~\ref{f:flips} also illustrates
  how a $Q$-adapted triangulation can be non-unique: if four interpolants
  forming a hinge are co-planar, both diagonals give $q^h_{ij}=0$.
  \begin{figure}
    \begin{center}
\begin{picture}(0,0)%
\includegraphics{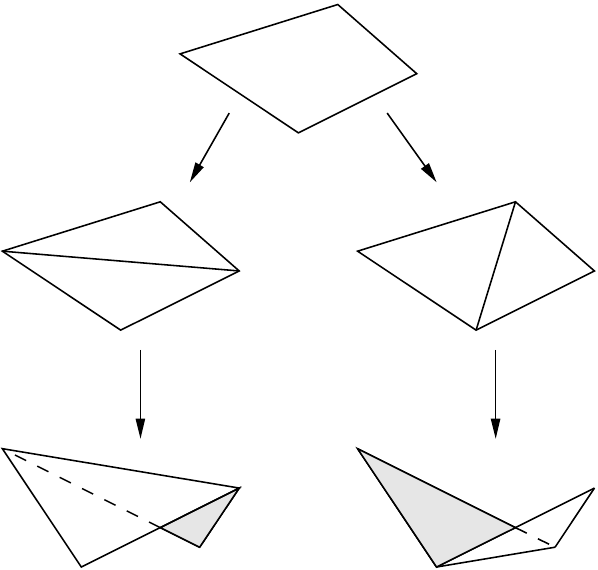}%
\end{picture}%
\setlength{\unitlength}{3947sp}%
\begingroup\makeatletter\ifx\SetFigFont\undefined%
\gdef\SetFigFont#1#2#3#4#5{%
  \reset@font\fontsize{#1}{#2pt}%
  \fontfamily{#3}\fontseries{#4}\fontshape{#5}%
  \selectfont}%
\fi\endgroup%
\begin{picture}(2866,2724)(2547,-4573)
\end{picture}%
    \end{center}
    \caption[Edge flips altering convexity] {Edge flips can replace non-convex
      edges, where $\Rightarrow q_{ab}<0$, with convex edges without changing
      the interpolated values $s_i$.  For the given hinge, the diagonal giving a
      negative edge is on the left; a positive edge is on the right.}
    \label{f:flips}
  \end{figure}
\end{remark}

\begin{proof}[Proof of Theorem~\ref{th:posq}]
  We proceed by constructing the triangulation as follows. Given $\calV$, we
  orthogonally project each 2D point onto the surface $s(x)$ corresponding to
  $Q$.  Take the convex hull of these points in 3D.  Orient each convex hull
  normal so that it faces outward from the convex hull.  Discard polyhedral
  faces of the convex hull with normals having positive $z$-components.
  Arbitrarily triangulate polyhedra on the convex hull which are not already
  triangles. The resulting triangulation, once projected back orthogonally onto
  the plane, is the $Q$-adapted triangulation.  Indeed, it is simple to show by
  direct calculation that hinge formula~\eqref{e:hinge} is invariant under the
  transformation $\{s_i \rightarrow s_i + ax_i + by_i + c\}$, where $a, b, c \in
  \real$ are constants independent of $i$.  This is consistent with the
  invariance of $Q$ under the addition of affine functions to $s(x)$.

  Now consider edge $ij$, referring to Figure~\ref{f:cotans}.  Due to the
  invariance under affine addition, we can add the affine function which results
  in $s_i = s_j = s_k = 0$.  Thus,
  \begin{equation}
    q^h_{ij} = \frac{1}{2\abs{t_{ijl}}} s_l.
  \end{equation}
  $\abs{t_{ijl}}$ is the unsigned triangle area, and so the sign of $q_{ij}$
  equals the sign of $s_l$.  That is, when $s_l$ lies above the $xy$-plane,
  $q_{ij}>0$, showing that the hinge is convex if and only if $q_{ij}>0$, and
  the hinge is flat if and only if $q_{ij}=0$.  All hinges on the convex hull of
  the interpolation of $s(x,y)$ are convex or flat, so all $q_{ij}\geq 0$, as
  expected. Moreover, $q^h_{ij}=0$ corresponds to a flat hinge, which in turn
  corresponds to an arbitrary triangulation of a polyhedron having four or more
  sides.  This is the only manner in which the $Q$-adapted triangulation can be
  non-unique.
\end{proof}

\subsection{Weighted Delaunay and $Q$-adapted triangulations}

There is a connection between $s(x)$ and weighted Delaunay triangulations, the dual graphs of ``power
diagrams.''  Glickenstein~\cite{glic2007} studies the discrete Dirichlet energy in context of weighted
Delaunay triangulations.  In the notation of~\eqref{e:hinge} and Figure~\ref{f:cotans}, Glickenstein
shows that for weights $w_i$, the coefficients of the discrete Dirichlet energy are \begin{equation}
  \begin{split}
    q^h_{ij} &=
    \onehalf \left(\cot\theta_{ikj} + \cot\theta_{ilj}\right) \\
    & \quad + \dfrac{1}{2\abs{e_{ij}}^2}\left(\cot\theta_{ijk} +
      \cot\theta_{ijl}\right) w_i \\
    & \quad + \dfrac{1}{2\abs{e_{ij}}^2}\left(\cot\theta_{jik} +
      \cot\theta_{jil}\right) w_j \\
    & \quad - \dfrac{1}{4 \abs{t_{ijk}}} w_k
    - \dfrac{1}{4 \abs{t_{ijl}}} w_l.
  \end{split}
\end{equation} Comparison of this formula with~\eqref{e:hinge} indicates that this is the
discretization of \begin{equation}
  q^h_{ij} = - \int_{\Omega} \Grad \myphi_i^T\left(I_2 - \onehalf Q_w\right) \Grad \myphi_j,
\end{equation} where $I_2$ is the $2\times 2$ identity matrix, and \begin{equation}
  Q_w = \begin{pmatrix} w_{yy} & -w_{xy} \\ -w_{xy} & w_{xx} \end{pmatrix}.
\end{equation} So, modulo addition of an arbitrary affine function, the interpolants \begin{equation}
  s_i = \onehalf \left(x_i^2 + y_i^2\right) - \onehalf w_i
\end{equation} can be used to compute Delaunay weights {\em from interpolants of $s(x)$}.

Thus, we have demonstrated the following connection between weighted Delaunay triangulations and
$Q$-adapted triangulations: \begin{proposition}\label{p:weightedDelaunay}
  Given a set of points $\calV$, the weighted Delaunay triangulation of those
  points having weights
  \begin{equation}
    w_i = x_i^2 + y_i^2 - 2s_i
    \label{e:weights}
  \end{equation}
  gives the same triangulation as that obtained by projecting the convex hull of
  points $\left(x_i, y_i, s_i\right)$ onto the $xy$-plane, where $s_i = s(x_i,
  y_i)$ are interpolants of the convex interpolation function $s(x)$.
\end{proposition}

\begin{remark}
  Weighted Delaunay can be efficiently computed by current computational
  geometry tools, see for instance~\cite{cgal}. Thus, we use such a weighted
  Delaunay algorithm instead of the convex hull construction to generate
  $Q$-adapted triangulations in our numerical tests below.
\end{remark} \begin{remark}
  In contrast to Delaunay meshes, weighted Delaunay triangulations do not
  necessarily contain all of the original points $\calV$.  The ``hidden'' points
  correspond to values $s_i$ that lie \emph{inside} the convex hull of the
  other interpolants of $s(x,y)$. In our setting, as long as we construct
  $w_i$ from $s_i$ interpolating a \emph{convex} function $s(x,y)$ (that is,
  weights representing a positive-definite $Q$), our weighted Delaunay
  triangulations do contain all the points in $\calV$.
\end{remark} \begin{remark}
  The triangulation is specific to $Q$, not to $s(x,y)$.  The addition of an affine
  function to $s(x,y)$ does \emph{not} alter the $q_{ij}$ given by the hinge formula, a
  fact which can be confirmed by direct calculation.  This is consistent with
  the observation that modifying the weights by the addition of an affine
  function $\{ w_i \rightarrow w_i + ax_i + by_i + c\}$, $a, b, c \in \real$ are
  constants independent of $i$, does not change the weighted Delaunay
  triangulation.  This can be seen by considering the dual graph determined by
  the points and their Delaunay weights, whereby adding an affine function to
  each of the weights only translates the dual graph in space, thereby leaving
  the triangulation unchanged.
\end{remark} \begin{remark}[Global energy minimum]
  The convex hull construction of a weighted Delaunay triangulation gives the
  global energy minimum result which is an extension of the result for the
  Delaunay triangulation.  That is, the discrete Dirichlet
  energy~\eqref{e:discreteDirichlet} with $q_{ij}$ computed using hinge
  formula~\eqref{e:hinge}, where $s_i$ interpolate a convex $s(x)$, gives the
  minimum energy for any given function $u_i$ provided the $q^h_{ij}$ are computed
  over the weighted Delaunay triangulation determined by
  weights~\eqref{e:weights}.

  To see this, consider the set of all triangulations of a fixed set of points.
  Each element of this set can be reached from every other element by performing
  a finite sequence of edge-flips.  The local result is that if two
  triangulations differ only in a single flip of an edge, and the triangulation
  is weighted Delaunay after the flip, then the latter triangulation gives the
  smaller Dirichlet energy.

  A global result is not possible for general weighted Delaunay triangulations
  because the choice of weights can give points with non-positive dual areas,
  whereupon these points do not appear in the final triangulation.  However, if
  the weights are computed from interpolation of a convex function, none of the
  points disappear, and the local result can be applied to arrive at the
  triangulation giving the global minimum of the Dirichlet norm.

  Similarly, if an arbitrary set of weights is used to construct interpolants
  $s_i$ using~\eqref{e:weights}, taking the convex hull of these points removes
  exactly those points which give non-positive dual areas.
  See comments in~\cite{glic2007} for further discussion of this global minimum
  result.
\end{remark}

\subsection{Computing optimal meshes}

Using the connection that we established between $s(x)$ and weighted Delaunay triangulations, we can
design a numerical procedure to produce high quality Q-adapted meshes. Although limited to
two-dimension, we extend the variational approach to isotropic meshing presented in~\cite{alli2005} to
anisotropic meshes. In our case, we seek a mesh that produces a matrix associated to the homogenized
problem~\eref{jashgajsggsjhw} having a small condition number, while still providing good
interpolations of the solution.

The variational approach in~\cite{alli2005} proceeds by moving points on a domain so as to improve
triangulation quality. At each step, the strategy is to adjust points to minimize, for the current
connectivity of the mesh, the cost function \begin{equation}
  E_p = \int_{\Omega} \abs{p(x,y) - p^h(x,y)},
\end{equation} where $p(x,y) = \onehalf (x^2 + y^2)$ and $p^h(x,y)$ is the piecewise linear
interpolation of $p(x,y)$ at each of the points.  That is, $p^h(x,y)$ inscribes $p(x,y)$, and $E_p$
represents the $L^1$ norm between the paraboloid and its piecewise linear interpolation based on the
current point positions and connectivity. The variational approach proceeds by using the critical
point of $E_p$ to update point locations iteratively~\cite{alli2005}.

Our extension consists of replacing the paraboloid $p(x,y)$ with the conductivity parameterization
$s(x,y)$. Computing the critical point of \begin{equation}
  E_s = \int_{\Omega} \abs{s(x,y) - \sh(x,y)}
\end{equation} with respect to point locations is found by solving \begin{multline}
    \Hess{s} (x_i^\ast,y_i^\ast) = \\
    \Hess{s} (x_i,y_i)
    - \dfrac{1}{\abs{K_i}} \sum_{t_j \in K_i}
    \left( \Grad_{(x_i,y_i)}\abs{t_j} \Biggl[
      \sum_{k \in t_j} s(x_k - x_i, y_k - y_i) \Biggr] \right)
  \label{e:simpleupdate}
\end{multline}

for the new position $(x_i^\ast,y_i^\ast)$. $\Hess{s}$ is the Hessian of $s(x,y)$, $K_i$ is the set of
trinagles adjacent to point $i$, $t_j$ is a triangle that belongs to $K_i$, and $\abs{t_j}$ is the
unsigned area of $t_j$.  Once the point positions have been updated in this fashion, we then recompute
a new tessellation based on these points and the weights $s_i$ through a weighted Delaunay algorithm
as detailed in the previous section.

\begin{algorithm}[Computing a $Q$-optimal mesh] Following~\cite{alli2005}, our
  algorithm for producing triangulations that lead to well conditioned stiffness
  matrices for the homogenized problem~\eqref{jashgajsggsjhw} is as follows:
  \begin{center}
    \begin{tabular}{l}
      \hline
      Read the interpolation function $s(x)$ \\
      Generate initial vertex positions $(x_i, y_i)$ inside $\Omega$ \\
      Do \\
      \hspace{0.5in} Compute triangulation weights using~\eqref{e:weights} \\
      \hspace{0.5in} Construct weighted Delaunay triangulation of the points \\
      \hspace{0.5in} Move points to their optimal positions using~\eqref{e:simpleupdate} \\
      Until (convergence or max iteration) \\
      \hline
    \end{tabular}
  \end{center}
\end{algorithm}

\begin{figure}
  \begin{center}
    \includegraphics[width=0.38\textwidth]{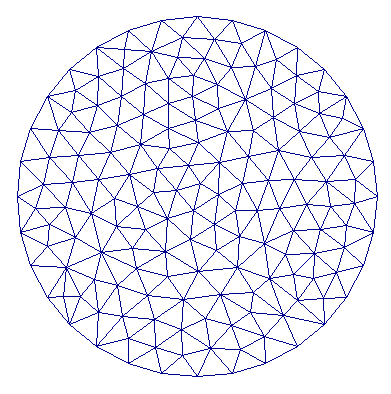}
    \hspace{0.1\textwidth}
    \includegraphics[width=0.38\textwidth]{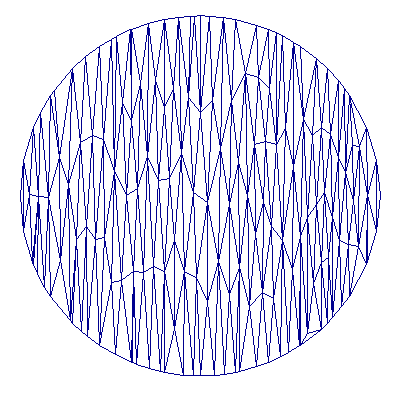}
  \end{center}
  \caption[Comparison of an isotropic and anisotropic mesh] {Comparison of an
    isotropic and an anisotropic mesh.  The figure on the left shows the lack of
    directional bias expected for a mesh suitable for the isotropic problem,
    while the figure on the right is suitable for the case where the
    conductivity is greater in the $y$-direction then in the
    $x$-direction.}
  \label{f:aniscompare}
\end{figure}

\begin{figure}
  \begin{center}
\setlength{\unitlength}{0.240900pt}
\ifx\plotpoint\undefined\newsavebox{\plotpoint}\fi
\begin{picture}(1500,900)(0,0)
\sbox{\plotpoint}{\rule[-0.200pt]{0.400pt}{0.400pt}}%
\put(241.0,123.0){\rule[-0.200pt]{4.818pt}{0.400pt}}
\put(221,123){\makebox(0,0)[r]{ 1}}
\put(1419.0,123.0){\rule[-0.200pt]{4.818pt}{0.400pt}}
\put(241.0,167.0){\rule[-0.200pt]{2.409pt}{0.400pt}}
\put(1429.0,167.0){\rule[-0.200pt]{2.409pt}{0.400pt}}
\put(241.0,226.0){\rule[-0.200pt]{2.409pt}{0.400pt}}
\put(1429.0,226.0){\rule[-0.200pt]{2.409pt}{0.400pt}}
\put(241.0,256.0){\rule[-0.200pt]{2.409pt}{0.400pt}}
\put(1429.0,256.0){\rule[-0.200pt]{2.409pt}{0.400pt}}
\put(241.0,270.0){\rule[-0.200pt]{4.818pt}{0.400pt}}
\put(221,270){\makebox(0,0)[r]{ 10}}
\put(1419.0,270.0){\rule[-0.200pt]{4.818pt}{0.400pt}}
\put(241.0,315.0){\rule[-0.200pt]{2.409pt}{0.400pt}}
\put(1429.0,315.0){\rule[-0.200pt]{2.409pt}{0.400pt}}
\put(241.0,373.0){\rule[-0.200pt]{2.409pt}{0.400pt}}
\put(1429.0,373.0){\rule[-0.200pt]{2.409pt}{0.400pt}}
\put(241.0,404.0){\rule[-0.200pt]{2.409pt}{0.400pt}}
\put(1429.0,404.0){\rule[-0.200pt]{2.409pt}{0.400pt}}
\put(241.0,418.0){\rule[-0.200pt]{4.818pt}{0.400pt}}
\put(221,418){\makebox(0,0)[r]{ 100}}
\put(1419.0,418.0){\rule[-0.200pt]{4.818pt}{0.400pt}}
\put(241.0,462.0){\rule[-0.200pt]{2.409pt}{0.400pt}}
\put(1429.0,462.0){\rule[-0.200pt]{2.409pt}{0.400pt}}
\put(241.0,521.0){\rule[-0.200pt]{2.409pt}{0.400pt}}
\put(1429.0,521.0){\rule[-0.200pt]{2.409pt}{0.400pt}}
\put(241.0,551.0){\rule[-0.200pt]{2.409pt}{0.400pt}}
\put(1429.0,551.0){\rule[-0.200pt]{2.409pt}{0.400pt}}
\put(241.0,565.0){\rule[-0.200pt]{4.818pt}{0.400pt}}
\put(221,565){\makebox(0,0)[r]{ 1000}}
\put(1419.0,565.0){\rule[-0.200pt]{4.818pt}{0.400pt}}
\put(241.0,610.0){\rule[-0.200pt]{2.409pt}{0.400pt}}
\put(1429.0,610.0){\rule[-0.200pt]{2.409pt}{0.400pt}}
\put(241.0,668.0){\rule[-0.200pt]{2.409pt}{0.400pt}}
\put(1429.0,668.0){\rule[-0.200pt]{2.409pt}{0.400pt}}
\put(241.0,698.0){\rule[-0.200pt]{2.409pt}{0.400pt}}
\put(1429.0,698.0){\rule[-0.200pt]{2.409pt}{0.400pt}}
\put(241.0,713.0){\rule[-0.200pt]{4.818pt}{0.400pt}}
\put(221,713){\makebox(0,0)[r]{ 10000}}
\put(1419.0,713.0){\rule[-0.200pt]{4.818pt}{0.400pt}}
\put(241.0,757.0){\rule[-0.200pt]{2.409pt}{0.400pt}}
\put(1429.0,757.0){\rule[-0.200pt]{2.409pt}{0.400pt}}
\put(241.0,816.0){\rule[-0.200pt]{2.409pt}{0.400pt}}
\put(1429.0,816.0){\rule[-0.200pt]{2.409pt}{0.400pt}}
\put(241.0,846.0){\rule[-0.200pt]{2.409pt}{0.400pt}}
\put(1429.0,846.0){\rule[-0.200pt]{2.409pt}{0.400pt}}
\put(241.0,860.0){\rule[-0.200pt]{4.818pt}{0.400pt}}
\put(221,860){\makebox(0,0)[r]{ 100000}}
\put(1419.0,860.0){\rule[-0.200pt]{4.818pt}{0.400pt}}
\put(241.0,123.0){\rule[-0.200pt]{0.400pt}{4.818pt}}
\put(241,82){\makebox(0,0){ 10}}
\put(241.0,840.0){\rule[-0.200pt]{0.400pt}{4.818pt}}
\put(331.0,123.0){\rule[-0.200pt]{0.400pt}{2.409pt}}
\put(331.0,850.0){\rule[-0.200pt]{0.400pt}{2.409pt}}
\put(384.0,123.0){\rule[-0.200pt]{0.400pt}{2.409pt}}
\put(384.0,850.0){\rule[-0.200pt]{0.400pt}{2.409pt}}
\put(421.0,123.0){\rule[-0.200pt]{0.400pt}{2.409pt}}
\put(421.0,850.0){\rule[-0.200pt]{0.400pt}{2.409pt}}
\put(450.0,123.0){\rule[-0.200pt]{0.400pt}{2.409pt}}
\put(450.0,850.0){\rule[-0.200pt]{0.400pt}{2.409pt}}
\put(474.0,123.0){\rule[-0.200pt]{0.400pt}{2.409pt}}
\put(474.0,850.0){\rule[-0.200pt]{0.400pt}{2.409pt}}
\put(494.0,123.0){\rule[-0.200pt]{0.400pt}{2.409pt}}
\put(494.0,850.0){\rule[-0.200pt]{0.400pt}{2.409pt}}
\put(511.0,123.0){\rule[-0.200pt]{0.400pt}{2.409pt}}
\put(511.0,850.0){\rule[-0.200pt]{0.400pt}{2.409pt}}
\put(527.0,123.0){\rule[-0.200pt]{0.400pt}{2.409pt}}
\put(527.0,850.0){\rule[-0.200pt]{0.400pt}{2.409pt}}
\put(541.0,123.0){\rule[-0.200pt]{0.400pt}{4.818pt}}
\put(541,82){\makebox(0,0){ 100}}
\put(541.0,840.0){\rule[-0.200pt]{0.400pt}{4.818pt}}
\put(631.0,123.0){\rule[-0.200pt]{0.400pt}{2.409pt}}
\put(631.0,850.0){\rule[-0.200pt]{0.400pt}{2.409pt}}
\put(683.0,123.0){\rule[-0.200pt]{0.400pt}{2.409pt}}
\put(683.0,850.0){\rule[-0.200pt]{0.400pt}{2.409pt}}
\put(721.0,123.0){\rule[-0.200pt]{0.400pt}{2.409pt}}
\put(721.0,850.0){\rule[-0.200pt]{0.400pt}{2.409pt}}
\put(750.0,123.0){\rule[-0.200pt]{0.400pt}{2.409pt}}
\put(750.0,850.0){\rule[-0.200pt]{0.400pt}{2.409pt}}
\put(774.0,123.0){\rule[-0.200pt]{0.400pt}{2.409pt}}
\put(774.0,850.0){\rule[-0.200pt]{0.400pt}{2.409pt}}
\put(794.0,123.0){\rule[-0.200pt]{0.400pt}{2.409pt}}
\put(794.0,850.0){\rule[-0.200pt]{0.400pt}{2.409pt}}
\put(811.0,123.0){\rule[-0.200pt]{0.400pt}{2.409pt}}
\put(811.0,850.0){\rule[-0.200pt]{0.400pt}{2.409pt}}
\put(826.0,123.0){\rule[-0.200pt]{0.400pt}{2.409pt}}
\put(826.0,850.0){\rule[-0.200pt]{0.400pt}{2.409pt}}
\put(840.0,123.0){\rule[-0.200pt]{0.400pt}{4.818pt}}
\put(840,82){\makebox(0,0){ 1000}}
\put(840.0,840.0){\rule[-0.200pt]{0.400pt}{4.818pt}}
\put(930.0,123.0){\rule[-0.200pt]{0.400pt}{2.409pt}}
\put(930.0,850.0){\rule[-0.200pt]{0.400pt}{2.409pt}}
\put(983.0,123.0){\rule[-0.200pt]{0.400pt}{2.409pt}}
\put(983.0,850.0){\rule[-0.200pt]{0.400pt}{2.409pt}}
\put(1020.0,123.0){\rule[-0.200pt]{0.400pt}{2.409pt}}
\put(1020.0,850.0){\rule[-0.200pt]{0.400pt}{2.409pt}}
\put(1049.0,123.0){\rule[-0.200pt]{0.400pt}{2.409pt}}
\put(1049.0,850.0){\rule[-0.200pt]{0.400pt}{2.409pt}}
\put(1073.0,123.0){\rule[-0.200pt]{0.400pt}{2.409pt}}
\put(1073.0,850.0){\rule[-0.200pt]{0.400pt}{2.409pt}}
\put(1093.0,123.0){\rule[-0.200pt]{0.400pt}{2.409pt}}
\put(1093.0,850.0){\rule[-0.200pt]{0.400pt}{2.409pt}}
\put(1110.0,123.0){\rule[-0.200pt]{0.400pt}{2.409pt}}
\put(1110.0,850.0){\rule[-0.200pt]{0.400pt}{2.409pt}}
\put(1126.0,123.0){\rule[-0.200pt]{0.400pt}{2.409pt}}
\put(1126.0,850.0){\rule[-0.200pt]{0.400pt}{2.409pt}}
\put(1140.0,123.0){\rule[-0.200pt]{0.400pt}{4.818pt}}
\put(1140,82){\makebox(0,0){ 10000}}
\put(1140.0,840.0){\rule[-0.200pt]{0.400pt}{4.818pt}}
\put(1230.0,123.0){\rule[-0.200pt]{0.400pt}{2.409pt}}
\put(1230.0,850.0){\rule[-0.200pt]{0.400pt}{2.409pt}}
\put(1282.0,123.0){\rule[-0.200pt]{0.400pt}{2.409pt}}
\put(1282.0,850.0){\rule[-0.200pt]{0.400pt}{2.409pt}}
\put(1320.0,123.0){\rule[-0.200pt]{0.400pt}{2.409pt}}
\put(1320.0,850.0){\rule[-0.200pt]{0.400pt}{2.409pt}}
\put(1349.0,123.0){\rule[-0.200pt]{0.400pt}{2.409pt}}
\put(1349.0,850.0){\rule[-0.200pt]{0.400pt}{2.409pt}}
\put(1373.0,123.0){\rule[-0.200pt]{0.400pt}{2.409pt}}
\put(1373.0,850.0){\rule[-0.200pt]{0.400pt}{2.409pt}}
\put(1393.0,123.0){\rule[-0.200pt]{0.400pt}{2.409pt}}
\put(1393.0,850.0){\rule[-0.200pt]{0.400pt}{2.409pt}}
\put(1410.0,123.0){\rule[-0.200pt]{0.400pt}{2.409pt}}
\put(1410.0,850.0){\rule[-0.200pt]{0.400pt}{2.409pt}}
\put(1425.0,123.0){\rule[-0.200pt]{0.400pt}{2.409pt}}
\put(1425.0,850.0){\rule[-0.200pt]{0.400pt}{2.409pt}}
\put(1439.0,123.0){\rule[-0.200pt]{0.400pt}{4.818pt}}
\put(1439,82){\makebox(0,0){ 100000}}
\put(1439.0,840.0){\rule[-0.200pt]{0.400pt}{4.818pt}}
\put(241.0,123.0){\rule[-0.200pt]{288.598pt}{0.400pt}}
\put(1439.0,123.0){\rule[-0.200pt]{0.400pt}{177.543pt}}
\put(241.0,860.0){\rule[-0.200pt]{288.598pt}{0.400pt}}
\put(241.0,123.0){\rule[-0.200pt]{0.400pt}{177.543pt}}
\put(40,491){\makebox(0,0){\begin{sideways}spectral condition number, $\kappa_2$\end{sideways}}}
\put(840,21){\makebox(0,0){number of vertices, $N$}}
\put(481,820){\makebox(0,0)[r]{isotropic}}
\put(501.0,820.0){\rule[-0.200pt]{24.090pt}{0.400pt}}
\put(519,360){\usebox{\plotpoint}}
\multiput(519.00,360.58)(0.923,0.499){111}{\rule{0.837pt}{0.120pt}}
\multiput(519.00,359.17)(103.263,57.000){2}{\rule{0.418pt}{0.400pt}}
\multiput(624.00,417.58)(0.914,0.498){99}{\rule{0.829pt}{0.120pt}}
\multiput(624.00,416.17)(91.279,51.000){2}{\rule{0.415pt}{0.400pt}}
\multiput(717.00,468.58)(0.923,0.499){111}{\rule{0.837pt}{0.120pt}}
\multiput(717.00,467.17)(103.263,57.000){2}{\rule{0.418pt}{0.400pt}}
\multiput(822.00,525.58)(0.956,0.498){83}{\rule{0.863pt}{0.120pt}}
\multiput(822.00,524.17)(80.209,43.000){2}{\rule{0.431pt}{0.400pt}}
\multiput(904.00,568.58)(0.973,0.497){61}{\rule{0.875pt}{0.120pt}}
\multiput(904.00,567.17)(60.184,32.000){2}{\rule{0.438pt}{0.400pt}}
\multiput(966.00,600.58)(0.959,0.498){89}{\rule{0.865pt}{0.120pt}}
\multiput(966.00,599.17)(86.204,46.000){2}{\rule{0.433pt}{0.400pt}}
\multiput(1054.00,646.58)(0.945,0.498){65}{\rule{0.853pt}{0.120pt}}
\multiput(1054.00,645.17)(62.230,34.000){2}{\rule{0.426pt}{0.400pt}}
\multiput(1118.00,680.58)(1.009,0.499){137}{\rule{0.906pt}{0.120pt}}
\multiput(1118.00,679.17)(139.120,70.000){2}{\rule{0.453pt}{0.400pt}}
\put(519,360){\raisebox{-.8pt}{\makebox(0,0){$\Diamond$}}}
\put(624,417){\raisebox{-.8pt}{\makebox(0,0){$\Diamond$}}}
\put(717,468){\raisebox{-.8pt}{\makebox(0,0){$\Diamond$}}}
\put(822,525){\raisebox{-.8pt}{\makebox(0,0){$\Diamond$}}}
\put(904,568){\raisebox{-.8pt}{\makebox(0,0){$\Diamond$}}}
\put(966,600){\raisebox{-.8pt}{\makebox(0,0){$\Diamond$}}}
\put(1054,646){\raisebox{-.8pt}{\makebox(0,0){$\Diamond$}}}
\put(1118,680){\raisebox{-.8pt}{\makebox(0,0){$\Diamond$}}}
\put(1259,750){\raisebox{-.8pt}{\makebox(0,0){$\Diamond$}}}
\put(551,820){\raisebox{-.8pt}{\makebox(0,0){$\Diamond$}}}
\put(481,779){\makebox(0,0)[r]{anisotropic}}
\put(501.0,779.0){\rule[-0.200pt]{24.090pt}{0.400pt}}
\put(542,238){\usebox{\plotpoint}}
\multiput(542.00,238.59)(2.552,0.485){11}{\rule{2.043pt}{0.117pt}}
\multiput(542.00,237.17)(29.760,7.000){2}{\rule{1.021pt}{0.400pt}}
\multiput(576.00,245.58)(0.608,0.498){99}{\rule{0.586pt}{0.120pt}}
\multiput(576.00,244.17)(60.783,51.000){2}{\rule{0.293pt}{0.400pt}}
\multiput(638.00,296.58)(1.530,0.492){21}{\rule{1.300pt}{0.119pt}}
\multiput(638.00,295.17)(33.302,12.000){2}{\rule{0.650pt}{0.400pt}}
\multiput(674.00,308.58)(1.175,0.496){45}{\rule{1.033pt}{0.120pt}}
\multiput(674.00,307.17)(53.855,24.000){2}{\rule{0.517pt}{0.400pt}}
\multiput(730.00,332.58)(1.464,0.495){35}{\rule{1.258pt}{0.119pt}}
\multiput(730.00,331.17)(52.389,19.000){2}{\rule{0.629pt}{0.400pt}}
\multiput(785.00,351.58)(0.858,0.495){31}{\rule{0.782pt}{0.119pt}}
\multiput(785.00,350.17)(27.376,17.000){2}{\rule{0.391pt}{0.400pt}}
\multiput(814.00,368.58)(1.222,0.497){53}{\rule{1.071pt}{0.120pt}}
\multiput(814.00,367.17)(65.776,28.000){2}{\rule{0.536pt}{0.400pt}}
\multiput(882.00,396.58)(0.861,0.498){75}{\rule{0.787pt}{0.120pt}}
\multiput(882.00,395.17)(65.366,39.000){2}{\rule{0.394pt}{0.400pt}}
\put(542,238){\raisebox{-.8pt}{\makebox(0,0){$\Box$}}}
\put(576,245){\raisebox{-.8pt}{\makebox(0,0){$\Box$}}}
\put(638,296){\raisebox{-.8pt}{\makebox(0,0){$\Box$}}}
\put(674,308){\raisebox{-.8pt}{\makebox(0,0){$\Box$}}}
\put(730,332){\raisebox{-.8pt}{\makebox(0,0){$\Box$}}}
\put(785,351){\raisebox{-.8pt}{\makebox(0,0){$\Box$}}}
\put(814,368){\raisebox{-.8pt}{\makebox(0,0){$\Box$}}}
\put(882,396){\raisebox{-.8pt}{\makebox(0,0){$\Box$}}}
\put(949,435){\raisebox{-.8pt}{\makebox(0,0){$\Box$}}}
\put(551,779){\raisebox{-.8pt}{\makebox(0,0){$\Box$}}}
\put(241.0,123.0){\rule[-0.200pt]{288.598pt}{0.400pt}}
\put(1439.0,123.0){\rule[-0.200pt]{0.400pt}{177.543pt}}
\put(241.0,860.0){\rule[-0.200pt]{288.598pt}{0.400pt}}
\put(241.0,123.0){\rule[-0.200pt]{0.400pt}{177.543pt}}
\end{picture}
  \end{center}
  \caption[Matrix conditioning of adapted meshes]{Matrix conditioning
    quality of adapted meshes measured by the spectral condition
    number $\kappa_2$ of the stiffness matrix.  The condition number
    grows as $\bigO(N)$ in both cases, but is offset by a factor of
    about 5 in the adapted anisotropic meshes. \label{f:comparecond}}
\end{figure}
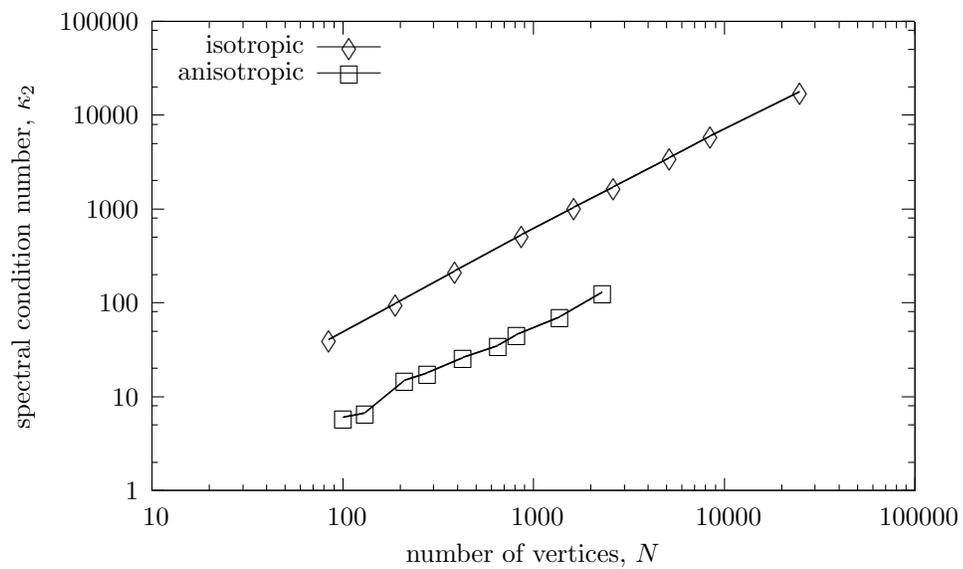

Figures~\ref{f:aniscompare} to~\ref{f:comparecond} give the results of a numerical experiment
illustrating the use of our algorithm for the case \begin{equation}
  Q = \begin{pmatrix} 0.1 & 0 \\ 0 & 10 \end{pmatrix}.
\end{equation} Consistent with theory, the quality measures of interpolation and matrix condition
number do not change at a greater rate than if an isotropic mesh is used with this conductivity.
However the constants in the performance metrics of the anisotropic meshes are less than those for the
isotropic meshes.

As a word of explanation, for a fixed number of points on the boundary of the domain, anisotropic
meshes tend to have fewer interior points than isotropic meshes.  Since the experimental meshes are
specified by their number of boundary points, this explains why the range of the total number of
vertices is greater for the isotropic meshes than the anisotropic meshes.

\begin{figure}
  \begin{center}
\setlength{\unitlength}{0.240900pt}
\ifx\plotpoint\undefined\newsavebox{\plotpoint}\fi
\begin{picture}(1500,900)(0,0)
\sbox{\plotpoint}{\rule[-0.200pt]{0.400pt}{0.400pt}}%
\put(241.0,123.0){\rule[-0.200pt]{4.818pt}{0.400pt}}
\put(221,123){\makebox(0,0)[r]{ 0.0001}}
\put(1419.0,123.0){\rule[-0.200pt]{4.818pt}{0.400pt}}
\put(241.0,197.0){\rule[-0.200pt]{2.409pt}{0.400pt}}
\put(1429.0,197.0){\rule[-0.200pt]{2.409pt}{0.400pt}}
\put(241.0,240.0){\rule[-0.200pt]{2.409pt}{0.400pt}}
\put(1429.0,240.0){\rule[-0.200pt]{2.409pt}{0.400pt}}
\put(241.0,271.0){\rule[-0.200pt]{2.409pt}{0.400pt}}
\put(1429.0,271.0){\rule[-0.200pt]{2.409pt}{0.400pt}}
\put(241.0,295.0){\rule[-0.200pt]{2.409pt}{0.400pt}}
\put(1429.0,295.0){\rule[-0.200pt]{2.409pt}{0.400pt}}
\put(241.0,314.0){\rule[-0.200pt]{2.409pt}{0.400pt}}
\put(1429.0,314.0){\rule[-0.200pt]{2.409pt}{0.400pt}}
\put(241.0,331.0){\rule[-0.200pt]{2.409pt}{0.400pt}}
\put(1429.0,331.0){\rule[-0.200pt]{2.409pt}{0.400pt}}
\put(241.0,345.0){\rule[-0.200pt]{2.409pt}{0.400pt}}
\put(1429.0,345.0){\rule[-0.200pt]{2.409pt}{0.400pt}}
\put(241.0,357.0){\rule[-0.200pt]{2.409pt}{0.400pt}}
\put(1429.0,357.0){\rule[-0.200pt]{2.409pt}{0.400pt}}
\put(241.0,369.0){\rule[-0.200pt]{4.818pt}{0.400pt}}
\put(221,369){\makebox(0,0)[r]{ 0.001}}
\put(1419.0,369.0){\rule[-0.200pt]{4.818pt}{0.400pt}}
\put(241.0,443.0){\rule[-0.200pt]{2.409pt}{0.400pt}}
\put(1429.0,443.0){\rule[-0.200pt]{2.409pt}{0.400pt}}
\put(241.0,486.0){\rule[-0.200pt]{2.409pt}{0.400pt}}
\put(1429.0,486.0){\rule[-0.200pt]{2.409pt}{0.400pt}}
\put(241.0,517.0){\rule[-0.200pt]{2.409pt}{0.400pt}}
\put(1429.0,517.0){\rule[-0.200pt]{2.409pt}{0.400pt}}
\put(241.0,540.0){\rule[-0.200pt]{2.409pt}{0.400pt}}
\put(1429.0,540.0){\rule[-0.200pt]{2.409pt}{0.400pt}}
\put(241.0,560.0){\rule[-0.200pt]{2.409pt}{0.400pt}}
\put(1429.0,560.0){\rule[-0.200pt]{2.409pt}{0.400pt}}
\put(241.0,576.0){\rule[-0.200pt]{2.409pt}{0.400pt}}
\put(1429.0,576.0){\rule[-0.200pt]{2.409pt}{0.400pt}}
\put(241.0,591.0){\rule[-0.200pt]{2.409pt}{0.400pt}}
\put(1429.0,591.0){\rule[-0.200pt]{2.409pt}{0.400pt}}
\put(241.0,603.0){\rule[-0.200pt]{2.409pt}{0.400pt}}
\put(1429.0,603.0){\rule[-0.200pt]{2.409pt}{0.400pt}}
\put(241.0,614.0){\rule[-0.200pt]{4.818pt}{0.400pt}}
\put(221,614){\makebox(0,0)[r]{ 0.01}}
\put(1419.0,614.0){\rule[-0.200pt]{4.818pt}{0.400pt}}
\put(241.0,688.0){\rule[-0.200pt]{2.409pt}{0.400pt}}
\put(1429.0,688.0){\rule[-0.200pt]{2.409pt}{0.400pt}}
\put(241.0,732.0){\rule[-0.200pt]{2.409pt}{0.400pt}}
\put(1429.0,732.0){\rule[-0.200pt]{2.409pt}{0.400pt}}
\put(241.0,762.0){\rule[-0.200pt]{2.409pt}{0.400pt}}
\put(1429.0,762.0){\rule[-0.200pt]{2.409pt}{0.400pt}}
\put(241.0,786.0){\rule[-0.200pt]{2.409pt}{0.400pt}}
\put(1429.0,786.0){\rule[-0.200pt]{2.409pt}{0.400pt}}
\put(241.0,805.0){\rule[-0.200pt]{2.409pt}{0.400pt}}
\put(1429.0,805.0){\rule[-0.200pt]{2.409pt}{0.400pt}}
\put(241.0,822.0){\rule[-0.200pt]{2.409pt}{0.400pt}}
\put(1429.0,822.0){\rule[-0.200pt]{2.409pt}{0.400pt}}
\put(241.0,836.0){\rule[-0.200pt]{2.409pt}{0.400pt}}
\put(1429.0,836.0){\rule[-0.200pt]{2.409pt}{0.400pt}}
\put(241.0,849.0){\rule[-0.200pt]{2.409pt}{0.400pt}}
\put(1429.0,849.0){\rule[-0.200pt]{2.409pt}{0.400pt}}
\put(241.0,860.0){\rule[-0.200pt]{4.818pt}{0.400pt}}
\put(221,860){\makebox(0,0)[r]{ 0.1}}
\put(1419.0,860.0){\rule[-0.200pt]{4.818pt}{0.400pt}}
\put(241.0,123.0){\rule[-0.200pt]{0.400pt}{4.818pt}}
\put(241,82){\makebox(0,0){ 10}}
\put(241.0,840.0){\rule[-0.200pt]{0.400pt}{4.818pt}}
\put(331.0,123.0){\rule[-0.200pt]{0.400pt}{2.409pt}}
\put(331.0,850.0){\rule[-0.200pt]{0.400pt}{2.409pt}}
\put(384.0,123.0){\rule[-0.200pt]{0.400pt}{2.409pt}}
\put(384.0,850.0){\rule[-0.200pt]{0.400pt}{2.409pt}}
\put(421.0,123.0){\rule[-0.200pt]{0.400pt}{2.409pt}}
\put(421.0,850.0){\rule[-0.200pt]{0.400pt}{2.409pt}}
\put(450.0,123.0){\rule[-0.200pt]{0.400pt}{2.409pt}}
\put(450.0,850.0){\rule[-0.200pt]{0.400pt}{2.409pt}}
\put(474.0,123.0){\rule[-0.200pt]{0.400pt}{2.409pt}}
\put(474.0,850.0){\rule[-0.200pt]{0.400pt}{2.409pt}}
\put(494.0,123.0){\rule[-0.200pt]{0.400pt}{2.409pt}}
\put(494.0,850.0){\rule[-0.200pt]{0.400pt}{2.409pt}}
\put(511.0,123.0){\rule[-0.200pt]{0.400pt}{2.409pt}}
\put(511.0,850.0){\rule[-0.200pt]{0.400pt}{2.409pt}}
\put(527.0,123.0){\rule[-0.200pt]{0.400pt}{2.409pt}}
\put(527.0,850.0){\rule[-0.200pt]{0.400pt}{2.409pt}}
\put(541.0,123.0){\rule[-0.200pt]{0.400pt}{4.818pt}}
\put(541,82){\makebox(0,0){ 100}}
\put(541.0,840.0){\rule[-0.200pt]{0.400pt}{4.818pt}}
\put(631.0,123.0){\rule[-0.200pt]{0.400pt}{2.409pt}}
\put(631.0,850.0){\rule[-0.200pt]{0.400pt}{2.409pt}}
\put(683.0,123.0){\rule[-0.200pt]{0.400pt}{2.409pt}}
\put(683.0,850.0){\rule[-0.200pt]{0.400pt}{2.409pt}}
\put(721.0,123.0){\rule[-0.200pt]{0.400pt}{2.409pt}}
\put(721.0,850.0){\rule[-0.200pt]{0.400pt}{2.409pt}}
\put(750.0,123.0){\rule[-0.200pt]{0.400pt}{2.409pt}}
\put(750.0,850.0){\rule[-0.200pt]{0.400pt}{2.409pt}}
\put(774.0,123.0){\rule[-0.200pt]{0.400pt}{2.409pt}}
\put(774.0,850.0){\rule[-0.200pt]{0.400pt}{2.409pt}}
\put(794.0,123.0){\rule[-0.200pt]{0.400pt}{2.409pt}}
\put(794.0,850.0){\rule[-0.200pt]{0.400pt}{2.409pt}}
\put(811.0,123.0){\rule[-0.200pt]{0.400pt}{2.409pt}}
\put(811.0,850.0){\rule[-0.200pt]{0.400pt}{2.409pt}}
\put(826.0,123.0){\rule[-0.200pt]{0.400pt}{2.409pt}}
\put(826.0,850.0){\rule[-0.200pt]{0.400pt}{2.409pt}}
\put(840.0,123.0){\rule[-0.200pt]{0.400pt}{4.818pt}}
\put(840,82){\makebox(0,0){ 1000}}
\put(840.0,840.0){\rule[-0.200pt]{0.400pt}{4.818pt}}
\put(930.0,123.0){\rule[-0.200pt]{0.400pt}{2.409pt}}
\put(930.0,850.0){\rule[-0.200pt]{0.400pt}{2.409pt}}
\put(983.0,123.0){\rule[-0.200pt]{0.400pt}{2.409pt}}
\put(983.0,850.0){\rule[-0.200pt]{0.400pt}{2.409pt}}
\put(1020.0,123.0){\rule[-0.200pt]{0.400pt}{2.409pt}}
\put(1020.0,850.0){\rule[-0.200pt]{0.400pt}{2.409pt}}
\put(1049.0,123.0){\rule[-0.200pt]{0.400pt}{2.409pt}}
\put(1049.0,850.0){\rule[-0.200pt]{0.400pt}{2.409pt}}
\put(1073.0,123.0){\rule[-0.200pt]{0.400pt}{2.409pt}}
\put(1073.0,850.0){\rule[-0.200pt]{0.400pt}{2.409pt}}
\put(1093.0,123.0){\rule[-0.200pt]{0.400pt}{2.409pt}}
\put(1093.0,850.0){\rule[-0.200pt]{0.400pt}{2.409pt}}
\put(1110.0,123.0){\rule[-0.200pt]{0.400pt}{2.409pt}}
\put(1110.0,850.0){\rule[-0.200pt]{0.400pt}{2.409pt}}
\put(1126.0,123.0){\rule[-0.200pt]{0.400pt}{2.409pt}}
\put(1126.0,850.0){\rule[-0.200pt]{0.400pt}{2.409pt}}
\put(1140.0,123.0){\rule[-0.200pt]{0.400pt}{4.818pt}}
\put(1140,82){\makebox(0,0){ 10000}}
\put(1140.0,840.0){\rule[-0.200pt]{0.400pt}{4.818pt}}
\put(1230.0,123.0){\rule[-0.200pt]{0.400pt}{2.409pt}}
\put(1230.0,850.0){\rule[-0.200pt]{0.400pt}{2.409pt}}
\put(1282.0,123.0){\rule[-0.200pt]{0.400pt}{2.409pt}}
\put(1282.0,850.0){\rule[-0.200pt]{0.400pt}{2.409pt}}
\put(1320.0,123.0){\rule[-0.200pt]{0.400pt}{2.409pt}}
\put(1320.0,850.0){\rule[-0.200pt]{0.400pt}{2.409pt}}
\put(1349.0,123.0){\rule[-0.200pt]{0.400pt}{2.409pt}}
\put(1349.0,850.0){\rule[-0.200pt]{0.400pt}{2.409pt}}
\put(1373.0,123.0){\rule[-0.200pt]{0.400pt}{2.409pt}}
\put(1373.0,850.0){\rule[-0.200pt]{0.400pt}{2.409pt}}
\put(1393.0,123.0){\rule[-0.200pt]{0.400pt}{2.409pt}}
\put(1393.0,850.0){\rule[-0.200pt]{0.400pt}{2.409pt}}
\put(1410.0,123.0){\rule[-0.200pt]{0.400pt}{2.409pt}}
\put(1410.0,850.0){\rule[-0.200pt]{0.400pt}{2.409pt}}
\put(1425.0,123.0){\rule[-0.200pt]{0.400pt}{2.409pt}}
\put(1425.0,850.0){\rule[-0.200pt]{0.400pt}{2.409pt}}
\put(1439.0,123.0){\rule[-0.200pt]{0.400pt}{4.818pt}}
\put(1439,82){\makebox(0,0){ 100000}}
\put(1439.0,840.0){\rule[-0.200pt]{0.400pt}{4.818pt}}
\put(241.0,123.0){\rule[-0.200pt]{288.598pt}{0.400pt}}
\put(1439.0,123.0){\rule[-0.200pt]{0.400pt}{177.543pt}}
\put(241.0,860.0){\rule[-0.200pt]{288.598pt}{0.400pt}}
\put(241.0,123.0){\rule[-0.200pt]{0.400pt}{177.543pt}}
\put(40,491){\makebox(0,0){\begin{sideways}$s$-interpolation error, $\norm{e_s}_{L_2}$\end{sideways}}}
\put(840,21){\makebox(0,0){number of vertices, $N$}}
\put(1279,820){\makebox(0,0)[r]{isotropic}}
\put(1299.0,820.0){\rule[-0.200pt]{24.090pt}{0.400pt}}
\put(519,810){\usebox{\plotpoint}}
\multiput(519.00,808.92)(0.500,-0.499){207}{\rule{0.500pt}{0.120pt}}
\multiput(519.00,809.17)(103.962,-105.000){2}{\rule{0.250pt}{0.400pt}}
\multiput(624.00,703.92)(0.560,-0.499){163}{\rule{0.548pt}{0.120pt}}
\multiput(624.00,704.17)(91.862,-83.000){2}{\rule{0.274pt}{0.400pt}}
\multiput(717.00,620.92)(0.577,-0.499){179}{\rule{0.562pt}{0.120pt}}
\multiput(717.00,621.17)(103.834,-91.000){2}{\rule{0.281pt}{0.400pt}}
\multiput(822.00,529.92)(0.641,-0.499){125}{\rule{0.613pt}{0.120pt}}
\multiput(822.00,530.17)(80.729,-64.000){2}{\rule{0.306pt}{0.400pt}}
\multiput(904.00,465.92)(0.596,-0.498){101}{\rule{0.577pt}{0.120pt}}
\multiput(904.00,466.17)(60.803,-52.000){2}{\rule{0.288pt}{0.400pt}}
\multiput(966.00,413.92)(0.595,-0.499){145}{\rule{0.576pt}{0.120pt}}
\multiput(966.00,414.17)(86.805,-74.000){2}{\rule{0.288pt}{0.400pt}}
\multiput(1054.00,339.92)(0.640,-0.498){97}{\rule{0.612pt}{0.120pt}}
\multiput(1054.00,340.17)(62.730,-50.000){2}{\rule{0.306pt}{0.400pt}}
\multiput(1118.00,289.92)(0.592,-0.499){235}{\rule{0.574pt}{0.120pt}}
\multiput(1118.00,290.17)(139.809,-119.000){2}{\rule{0.287pt}{0.400pt}}
\put(519,810){\raisebox{-.8pt}{\makebox(0,0){$\Diamond$}}}
\put(624,705){\raisebox{-.8pt}{\makebox(0,0){$\Diamond$}}}
\put(717,622){\raisebox{-.8pt}{\makebox(0,0){$\Diamond$}}}
\put(822,531){\raisebox{-.8pt}{\makebox(0,0){$\Diamond$}}}
\put(904,467){\raisebox{-.8pt}{\makebox(0,0){$\Diamond$}}}
\put(966,415){\raisebox{-.8pt}{\makebox(0,0){$\Diamond$}}}
\put(1054,341){\raisebox{-.8pt}{\makebox(0,0){$\Diamond$}}}
\put(1118,291){\raisebox{-.8pt}{\makebox(0,0){$\Diamond$}}}
\put(1259,172){\raisebox{-.8pt}{\makebox(0,0){$\Diamond$}}}
\put(1349,820){\raisebox{-.8pt}{\makebox(0,0){$\Diamond$}}}
\put(1279,779){\makebox(0,0)[r]{anisotropic}}
\put(1299.0,779.0){\rule[-0.200pt]{24.090pt}{0.400pt}}
\put(542,627){\usebox{\plotpoint}}
\multiput(542.00,625.92)(1.444,-0.492){21}{\rule{1.233pt}{0.119pt}}
\multiput(542.00,626.17)(31.440,-12.000){2}{\rule{0.617pt}{0.400pt}}
\multiput(576.58,612.39)(0.499,-0.662){121}{\rule{0.120pt}{0.629pt}}
\multiput(575.17,613.69)(62.000,-80.694){2}{\rule{0.400pt}{0.315pt}}
\multiput(638.00,531.92)(1.215,-0.494){27}{\rule{1.060pt}{0.119pt}}
\multiput(638.00,532.17)(33.800,-15.000){2}{\rule{0.530pt}{0.400pt}}
\multiput(674.00,516.92)(0.508,-0.499){107}{\rule{0.507pt}{0.120pt}}
\multiput(674.00,517.17)(54.947,-55.000){2}{\rule{0.254pt}{0.400pt}}
\multiput(730.00,461.92)(0.725,-0.498){73}{\rule{0.679pt}{0.120pt}}
\multiput(730.00,462.17)(53.591,-38.000){2}{\rule{0.339pt}{0.400pt}}
\multiput(785.00,423.92)(0.913,-0.494){29}{\rule{0.825pt}{0.119pt}}
\multiput(785.00,424.17)(27.288,-16.000){2}{\rule{0.413pt}{0.400pt}}
\multiput(814.00,407.92)(0.654,-0.498){101}{\rule{0.623pt}{0.120pt}}
\multiput(814.00,408.17)(66.707,-52.000){2}{\rule{0.312pt}{0.400pt}}
\multiput(882.00,355.92)(0.507,-0.499){129}{\rule{0.506pt}{0.120pt}}
\multiput(882.00,356.17)(65.950,-66.000){2}{\rule{0.253pt}{0.400pt}}
\put(542,627){\raisebox{-.8pt}{\makebox(0,0){$\Box$}}}
\put(576,615){\raisebox{-.8pt}{\makebox(0,0){$\Box$}}}
\put(638,533){\raisebox{-.8pt}{\makebox(0,0){$\Box$}}}
\put(674,518){\raisebox{-.8pt}{\makebox(0,0){$\Box$}}}
\put(730,463){\raisebox{-.8pt}{\makebox(0,0){$\Box$}}}
\put(785,425){\raisebox{-.8pt}{\makebox(0,0){$\Box$}}}
\put(814,409){\raisebox{-.8pt}{\makebox(0,0){$\Box$}}}
\put(882,357){\raisebox{-.8pt}{\makebox(0,0){$\Box$}}}
\put(949,291){\raisebox{-.8pt}{\makebox(0,0){$\Box$}}}
\put(1349,779){\raisebox{-.8pt}{\makebox(0,0){$\Box$}}}
\put(241.0,123.0){\rule[-0.200pt]{288.598pt}{0.400pt}}
\put(1439.0,123.0){\rule[-0.200pt]{0.400pt}{177.543pt}}
\put(241.0,860.0){\rule[-0.200pt]{288.598pt}{0.400pt}}
\put(241.0,123.0){\rule[-0.200pt]{0.400pt}{177.543pt}}
\end{picture}
  \end{center}
  \caption[$L_2$ interpolation quality of adapted
  meshes]{Interpolation quality of adapted meshes measured by the
  $L_2$-norm error in a linear interpolation of $s(x,y)$.  Error
  diminishes as $\bigO(N^{-1})$ in both cases, but is offset by a factor of
  about 4 in the adapted anisotropic meshes. \label{f:compareL2}}
\end{figure}

\begin{figure}
  \begin{center}
\setlength{\unitlength}{0.240900pt}
\ifx\plotpoint\undefined\newsavebox{\plotpoint}\fi
\begin{picture}(1500,900)(0,0)
\sbox{\plotpoint}{\rule[-0.200pt]{0.400pt}{0.400pt}}%
\put(201.0,123.0){\rule[-0.200pt]{4.818pt}{0.400pt}}
\put(181,123){\makebox(0,0)[r]{ 0.01}}
\put(1419.0,123.0){\rule[-0.200pt]{4.818pt}{0.400pt}}
\put(201.0,197.0){\rule[-0.200pt]{2.409pt}{0.400pt}}
\put(1429.0,197.0){\rule[-0.200pt]{2.409pt}{0.400pt}}
\put(201.0,240.0){\rule[-0.200pt]{2.409pt}{0.400pt}}
\put(1429.0,240.0){\rule[-0.200pt]{2.409pt}{0.400pt}}
\put(201.0,271.0){\rule[-0.200pt]{2.409pt}{0.400pt}}
\put(1429.0,271.0){\rule[-0.200pt]{2.409pt}{0.400pt}}
\put(201.0,295.0){\rule[-0.200pt]{2.409pt}{0.400pt}}
\put(1429.0,295.0){\rule[-0.200pt]{2.409pt}{0.400pt}}
\put(201.0,314.0){\rule[-0.200pt]{2.409pt}{0.400pt}}
\put(1429.0,314.0){\rule[-0.200pt]{2.409pt}{0.400pt}}
\put(201.0,331.0){\rule[-0.200pt]{2.409pt}{0.400pt}}
\put(1429.0,331.0){\rule[-0.200pt]{2.409pt}{0.400pt}}
\put(201.0,345.0){\rule[-0.200pt]{2.409pt}{0.400pt}}
\put(1429.0,345.0){\rule[-0.200pt]{2.409pt}{0.400pt}}
\put(201.0,357.0){\rule[-0.200pt]{2.409pt}{0.400pt}}
\put(1429.0,357.0){\rule[-0.200pt]{2.409pt}{0.400pt}}
\put(201.0,369.0){\rule[-0.200pt]{4.818pt}{0.400pt}}
\put(181,369){\makebox(0,0)[r]{ 0.1}}
\put(1419.0,369.0){\rule[-0.200pt]{4.818pt}{0.400pt}}
\put(201.0,443.0){\rule[-0.200pt]{2.409pt}{0.400pt}}
\put(1429.0,443.0){\rule[-0.200pt]{2.409pt}{0.400pt}}
\put(201.0,486.0){\rule[-0.200pt]{2.409pt}{0.400pt}}
\put(1429.0,486.0){\rule[-0.200pt]{2.409pt}{0.400pt}}
\put(201.0,517.0){\rule[-0.200pt]{2.409pt}{0.400pt}}
\put(1429.0,517.0){\rule[-0.200pt]{2.409pt}{0.400pt}}
\put(201.0,540.0){\rule[-0.200pt]{2.409pt}{0.400pt}}
\put(1429.0,540.0){\rule[-0.200pt]{2.409pt}{0.400pt}}
\put(201.0,560.0){\rule[-0.200pt]{2.409pt}{0.400pt}}
\put(1429.0,560.0){\rule[-0.200pt]{2.409pt}{0.400pt}}
\put(201.0,576.0){\rule[-0.200pt]{2.409pt}{0.400pt}}
\put(1429.0,576.0){\rule[-0.200pt]{2.409pt}{0.400pt}}
\put(201.0,591.0){\rule[-0.200pt]{2.409pt}{0.400pt}}
\put(1429.0,591.0){\rule[-0.200pt]{2.409pt}{0.400pt}}
\put(201.0,603.0){\rule[-0.200pt]{2.409pt}{0.400pt}}
\put(1429.0,603.0){\rule[-0.200pt]{2.409pt}{0.400pt}}
\put(201.0,614.0){\rule[-0.200pt]{4.818pt}{0.400pt}}
\put(181,614){\makebox(0,0)[r]{ 1}}
\put(1419.0,614.0){\rule[-0.200pt]{4.818pt}{0.400pt}}
\put(201.0,688.0){\rule[-0.200pt]{2.409pt}{0.400pt}}
\put(1429.0,688.0){\rule[-0.200pt]{2.409pt}{0.400pt}}
\put(201.0,732.0){\rule[-0.200pt]{2.409pt}{0.400pt}}
\put(1429.0,732.0){\rule[-0.200pt]{2.409pt}{0.400pt}}
\put(201.0,762.0){\rule[-0.200pt]{2.409pt}{0.400pt}}
\put(1429.0,762.0){\rule[-0.200pt]{2.409pt}{0.400pt}}
\put(201.0,786.0){\rule[-0.200pt]{2.409pt}{0.400pt}}
\put(1429.0,786.0){\rule[-0.200pt]{2.409pt}{0.400pt}}
\put(201.0,805.0){\rule[-0.200pt]{2.409pt}{0.400pt}}
\put(1429.0,805.0){\rule[-0.200pt]{2.409pt}{0.400pt}}
\put(201.0,822.0){\rule[-0.200pt]{2.409pt}{0.400pt}}
\put(1429.0,822.0){\rule[-0.200pt]{2.409pt}{0.400pt}}
\put(201.0,836.0){\rule[-0.200pt]{2.409pt}{0.400pt}}
\put(1429.0,836.0){\rule[-0.200pt]{2.409pt}{0.400pt}}
\put(201.0,849.0){\rule[-0.200pt]{2.409pt}{0.400pt}}
\put(1429.0,849.0){\rule[-0.200pt]{2.409pt}{0.400pt}}
\put(201.0,860.0){\rule[-0.200pt]{4.818pt}{0.400pt}}
\put(181,860){\makebox(0,0)[r]{ 10}}
\put(1419.0,860.0){\rule[-0.200pt]{4.818pt}{0.400pt}}
\put(201.0,123.0){\rule[-0.200pt]{0.400pt}{4.818pt}}
\put(201,82){\makebox(0,0){ 10}}
\put(201.0,840.0){\rule[-0.200pt]{0.400pt}{4.818pt}}
\put(294.0,123.0){\rule[-0.200pt]{0.400pt}{2.409pt}}
\put(294.0,850.0){\rule[-0.200pt]{0.400pt}{2.409pt}}
\put(349.0,123.0){\rule[-0.200pt]{0.400pt}{2.409pt}}
\put(349.0,850.0){\rule[-0.200pt]{0.400pt}{2.409pt}}
\put(387.0,123.0){\rule[-0.200pt]{0.400pt}{2.409pt}}
\put(387.0,850.0){\rule[-0.200pt]{0.400pt}{2.409pt}}
\put(417.0,123.0){\rule[-0.200pt]{0.400pt}{2.409pt}}
\put(417.0,850.0){\rule[-0.200pt]{0.400pt}{2.409pt}}
\put(442.0,123.0){\rule[-0.200pt]{0.400pt}{2.409pt}}
\put(442.0,850.0){\rule[-0.200pt]{0.400pt}{2.409pt}}
\put(463.0,123.0){\rule[-0.200pt]{0.400pt}{2.409pt}}
\put(463.0,850.0){\rule[-0.200pt]{0.400pt}{2.409pt}}
\put(481.0,123.0){\rule[-0.200pt]{0.400pt}{2.409pt}}
\put(481.0,850.0){\rule[-0.200pt]{0.400pt}{2.409pt}}
\put(496.0,123.0){\rule[-0.200pt]{0.400pt}{2.409pt}}
\put(496.0,850.0){\rule[-0.200pt]{0.400pt}{2.409pt}}
\put(511.0,123.0){\rule[-0.200pt]{0.400pt}{4.818pt}}
\put(511,82){\makebox(0,0){ 100}}
\put(511.0,840.0){\rule[-0.200pt]{0.400pt}{4.818pt}}
\put(604.0,123.0){\rule[-0.200pt]{0.400pt}{2.409pt}}
\put(604.0,850.0){\rule[-0.200pt]{0.400pt}{2.409pt}}
\put(658.0,123.0){\rule[-0.200pt]{0.400pt}{2.409pt}}
\put(658.0,850.0){\rule[-0.200pt]{0.400pt}{2.409pt}}
\put(697.0,123.0){\rule[-0.200pt]{0.400pt}{2.409pt}}
\put(697.0,850.0){\rule[-0.200pt]{0.400pt}{2.409pt}}
\put(727.0,123.0){\rule[-0.200pt]{0.400pt}{2.409pt}}
\put(727.0,850.0){\rule[-0.200pt]{0.400pt}{2.409pt}}
\put(751.0,123.0){\rule[-0.200pt]{0.400pt}{2.409pt}}
\put(751.0,850.0){\rule[-0.200pt]{0.400pt}{2.409pt}}
\put(772.0,123.0){\rule[-0.200pt]{0.400pt}{2.409pt}}
\put(772.0,850.0){\rule[-0.200pt]{0.400pt}{2.409pt}}
\put(790.0,123.0){\rule[-0.200pt]{0.400pt}{2.409pt}}
\put(790.0,850.0){\rule[-0.200pt]{0.400pt}{2.409pt}}
\put(806.0,123.0){\rule[-0.200pt]{0.400pt}{2.409pt}}
\put(806.0,850.0){\rule[-0.200pt]{0.400pt}{2.409pt}}
\put(820.0,123.0){\rule[-0.200pt]{0.400pt}{4.818pt}}
\put(820,82){\makebox(0,0){ 1000}}
\put(820.0,840.0){\rule[-0.200pt]{0.400pt}{4.818pt}}
\put(913.0,123.0){\rule[-0.200pt]{0.400pt}{2.409pt}}
\put(913.0,850.0){\rule[-0.200pt]{0.400pt}{2.409pt}}
\put(968.0,123.0){\rule[-0.200pt]{0.400pt}{2.409pt}}
\put(968.0,850.0){\rule[-0.200pt]{0.400pt}{2.409pt}}
\put(1006.0,123.0){\rule[-0.200pt]{0.400pt}{2.409pt}}
\put(1006.0,850.0){\rule[-0.200pt]{0.400pt}{2.409pt}}
\put(1036.0,123.0){\rule[-0.200pt]{0.400pt}{2.409pt}}
\put(1036.0,850.0){\rule[-0.200pt]{0.400pt}{2.409pt}}
\put(1061.0,123.0){\rule[-0.200pt]{0.400pt}{2.409pt}}
\put(1061.0,850.0){\rule[-0.200pt]{0.400pt}{2.409pt}}
\put(1082.0,123.0){\rule[-0.200pt]{0.400pt}{2.409pt}}
\put(1082.0,850.0){\rule[-0.200pt]{0.400pt}{2.409pt}}
\put(1100.0,123.0){\rule[-0.200pt]{0.400pt}{2.409pt}}
\put(1100.0,850.0){\rule[-0.200pt]{0.400pt}{2.409pt}}
\put(1115.0,123.0){\rule[-0.200pt]{0.400pt}{2.409pt}}
\put(1115.0,850.0){\rule[-0.200pt]{0.400pt}{2.409pt}}
\put(1130.0,123.0){\rule[-0.200pt]{0.400pt}{4.818pt}}
\put(1130,82){\makebox(0,0){ 10000}}
\put(1130.0,840.0){\rule[-0.200pt]{0.400pt}{4.818pt}}
\put(1223.0,123.0){\rule[-0.200pt]{0.400pt}{2.409pt}}
\put(1223.0,850.0){\rule[-0.200pt]{0.400pt}{2.409pt}}
\put(1277.0,123.0){\rule[-0.200pt]{0.400pt}{2.409pt}}
\put(1277.0,850.0){\rule[-0.200pt]{0.400pt}{2.409pt}}
\put(1316.0,123.0){\rule[-0.200pt]{0.400pt}{2.409pt}}
\put(1316.0,850.0){\rule[-0.200pt]{0.400pt}{2.409pt}}
\put(1346.0,123.0){\rule[-0.200pt]{0.400pt}{2.409pt}}
\put(1346.0,850.0){\rule[-0.200pt]{0.400pt}{2.409pt}}
\put(1370.0,123.0){\rule[-0.200pt]{0.400pt}{2.409pt}}
\put(1370.0,850.0){\rule[-0.200pt]{0.400pt}{2.409pt}}
\put(1391.0,123.0){\rule[-0.200pt]{0.400pt}{2.409pt}}
\put(1391.0,850.0){\rule[-0.200pt]{0.400pt}{2.409pt}}
\put(1409.0,123.0){\rule[-0.200pt]{0.400pt}{2.409pt}}
\put(1409.0,850.0){\rule[-0.200pt]{0.400pt}{2.409pt}}
\put(1425.0,123.0){\rule[-0.200pt]{0.400pt}{2.409pt}}
\put(1425.0,850.0){\rule[-0.200pt]{0.400pt}{2.409pt}}
\put(1439.0,123.0){\rule[-0.200pt]{0.400pt}{4.818pt}}
\put(1439,82){\makebox(0,0){ 100000}}
\put(1439.0,840.0){\rule[-0.200pt]{0.400pt}{4.818pt}}
\put(201.0,123.0){\rule[-0.200pt]{298.234pt}{0.400pt}}
\put(1439.0,123.0){\rule[-0.200pt]{0.400pt}{177.543pt}}
\put(201.0,860.0){\rule[-0.200pt]{298.234pt}{0.400pt}}
\put(201.0,123.0){\rule[-0.200pt]{0.400pt}{177.543pt}}
\put(40,491){\makebox(0,0){\begin{sideways}$s$-interpolation error, $\abs{e_s}_{h^1}$\end{sideways}}}
\put(820,21){\makebox(0,0){number of vertices, $N$}}
\put(1279,820){\makebox(0,0)[r]{isotropic}}
\put(1299.0,820.0){\rule[-0.200pt]{24.090pt}{0.400pt}}
\put(489,638){\usebox{\plotpoint}}
\multiput(489.00,636.92)(0.967,-0.499){109}{\rule{0.871pt}{0.120pt}}
\multiput(489.00,637.17)(106.191,-56.000){2}{\rule{0.436pt}{0.400pt}}
\multiput(597.00,580.92)(1.176,-0.498){79}{\rule{1.037pt}{0.120pt}}
\multiput(597.00,581.17)(93.849,-41.000){2}{\rule{0.518pt}{0.400pt}}
\multiput(693.00,539.92)(1.189,-0.498){89}{\rule{1.048pt}{0.120pt}}
\multiput(693.00,540.17)(106.825,-46.000){2}{\rule{0.524pt}{0.400pt}}
\multiput(802.00,493.92)(1.320,-0.497){61}{\rule{1.150pt}{0.120pt}}
\multiput(802.00,494.17)(81.613,-32.000){2}{\rule{0.575pt}{0.400pt}}
\multiput(886.00,461.92)(1.290,-0.497){47}{\rule{1.124pt}{0.120pt}}
\multiput(886.00,462.17)(61.667,-25.000){2}{\rule{0.562pt}{0.400pt}}
\multiput(950.00,436.92)(1.249,-0.498){71}{\rule{1.095pt}{0.120pt}}
\multiput(950.00,437.17)(89.728,-37.000){2}{\rule{0.547pt}{0.400pt}}
\multiput(1042.00,399.92)(1.310,-0.497){47}{\rule{1.140pt}{0.120pt}}
\multiput(1042.00,400.17)(62.634,-25.000){2}{\rule{0.570pt}{0.400pt}}
\multiput(1107.00,374.92)(1.220,-0.499){117}{\rule{1.073pt}{0.120pt}}
\multiput(1107.00,375.17)(143.772,-60.000){2}{\rule{0.537pt}{0.400pt}}
\put(489,638){\raisebox{-.8pt}{\makebox(0,0){$\Diamond$}}}
\put(597,582){\raisebox{-.8pt}{\makebox(0,0){$\Diamond$}}}
\put(693,541){\raisebox{-.8pt}{\makebox(0,0){$\Diamond$}}}
\put(802,495){\raisebox{-.8pt}{\makebox(0,0){$\Diamond$}}}
\put(886,463){\raisebox{-.8pt}{\makebox(0,0){$\Diamond$}}}
\put(950,438){\raisebox{-.8pt}{\makebox(0,0){$\Diamond$}}}
\put(1042,401){\raisebox{-.8pt}{\makebox(0,0){$\Diamond$}}}
\put(1107,376){\raisebox{-.8pt}{\makebox(0,0){$\Diamond$}}}
\put(1253,316){\raisebox{-.8pt}{\makebox(0,0){$\Diamond$}}}
\put(1349,820){\raisebox{-.8pt}{\makebox(0,0){$\Diamond$}}}
\put(1279,779){\makebox(0,0)[r]{anisotropic}}
\put(1299.0,779.0){\rule[-0.200pt]{24.090pt}{0.400pt}}
\put(512,500){\usebox{\plotpoint}}
\multiput(512.00,498.92)(1.798,-0.491){17}{\rule{1.500pt}{0.118pt}}
\multiput(512.00,499.17)(31.887,-10.000){2}{\rule{0.750pt}{0.400pt}}
\multiput(547.00,488.92)(0.835,-0.498){75}{\rule{0.767pt}{0.120pt}}
\multiput(547.00,489.17)(63.409,-39.000){2}{\rule{0.383pt}{0.400pt}}
\multiput(612.00,449.92)(1.573,-0.492){21}{\rule{1.333pt}{0.119pt}}
\multiput(612.00,450.17)(34.233,-12.000){2}{\rule{0.667pt}{0.400pt}}
\multiput(649.00,437.92)(1.103,-0.497){49}{\rule{0.977pt}{0.120pt}}
\multiput(649.00,438.17)(54.972,-26.000){2}{\rule{0.488pt}{0.400pt}}
\multiput(706.00,411.92)(1.196,-0.496){45}{\rule{1.050pt}{0.120pt}}
\multiput(706.00,412.17)(54.821,-24.000){2}{\rule{0.525pt}{0.400pt}}
\multiput(763.00,387.93)(2.247,-0.485){11}{\rule{1.814pt}{0.117pt}}
\multiput(763.00,388.17)(26.234,-7.000){2}{\rule{0.907pt}{0.400pt}}
\multiput(793.00,380.92)(1.357,-0.497){49}{\rule{1.177pt}{0.120pt}}
\multiput(793.00,381.17)(67.557,-26.000){2}{\rule{0.588pt}{0.400pt}}
\multiput(863.00,354.92)(1.099,-0.497){61}{\rule{0.975pt}{0.120pt}}
\multiput(863.00,355.17)(67.976,-32.000){2}{\rule{0.488pt}{0.400pt}}
\put(512,500){\raisebox{-.8pt}{\makebox(0,0){$\Box$}}}
\put(547,490){\raisebox{-.8pt}{\makebox(0,0){$\Box$}}}
\put(612,451){\raisebox{-.8pt}{\makebox(0,0){$\Box$}}}
\put(649,439){\raisebox{-.8pt}{\makebox(0,0){$\Box$}}}
\put(706,413){\raisebox{-.8pt}{\makebox(0,0){$\Box$}}}
\put(763,389){\raisebox{-.8pt}{\makebox(0,0){$\Box$}}}
\put(793,382){\raisebox{-.8pt}{\makebox(0,0){$\Box$}}}
\put(863,356){\raisebox{-.8pt}{\makebox(0,0){$\Box$}}}
\put(933,324){\raisebox{-.8pt}{\makebox(0,0){$\Box$}}}
\put(1349,779){\raisebox{-.8pt}{\makebox(0,0){$\Box$}}}
\put(201.0,123.0){\rule[-0.200pt]{298.234pt}{0.400pt}}
\put(1439.0,123.0){\rule[-0.200pt]{0.400pt}{177.543pt}}
\put(201.0,860.0){\rule[-0.200pt]{298.234pt}{0.400pt}}
\put(201.0,123.0){\rule[-0.200pt]{0.400pt}{177.543pt}}
\end{picture}
  \end{center}
  \caption[$h_1$ interpolation quality of adapted
  meshes]{Interpolation quality of adapted meshes measured by the
  $h^1$-semi-norm error in a linear interpolation of $s(x,y)$.  Error
  diminishes as $\bigO(N^{-\onehalf})$ in both cases, but is offset by
  a factor of about 3 in the adapted anisotropic meshes. \label{f:compareh1}}
\end{figure}

\section{Relation with inverse homogenization}\label{lshkdgshdgjjhdsg}
Consider the following sequence of PDEs indexed by $\epsilon$.
\begin{equation}\label{sdsddsefrfrc}
  \left\{
  \begin{aligned}
    -\Div\left(\sigma(\frac{x}{\epsilon}) \Grad u_\epsilon\right) &= f, \enspace &x\in\Omega, \\
    u &= 0, \enspace &x\in\partial\Omega.
  \end{aligned}
  \right.
\end{equation}
Assuming that $y\rightarrow \sigma(y)$ is periodic (with in $L^\infty(\mathbb{T}^d)$, where $\mathbb{T}^d$ is the $d$-dimensional unit torus) we know from classical homogenization theory \cite{BeLiPa78} that
$u_\epsilon$ converges towards $u_0$ as $\epsilon \downarrow 0$ where $u_0$ is the solution of the following PDE
\begin{equation}\label{sdsdddfdfersefrfrc}
  \left\{
  \begin{aligned}
    -\Div\left(\sigma_e \Grad u_0\right) &= f, \enspace &x\in\Omega, \\
    u &= 0, \enspace &x\in\partial\Omega.
  \end{aligned}
  \right.
\end{equation}
Moreover $\sigma_e$ is constant positive definite $d\times d$ matrix defined by
\begin{equation}\label{kdsjdhskj}
\sigma_e:=\int_{\mathbb{T}^d} \sigma(y)(I_d+\nabla \chi(y))
\end{equation}
where the entries of the vector field $\chi:=(\chi_1,\ldots,\chi_d)$ are solutions of the cell problems
\begin{equation}\label{sdssdssddsfrc}
  \left\{
  \begin{aligned}
    -\Div\left(\sigma(y) \Grad (y_i+\chi_i(y))\right) = 0, \enspace &y\in\mathbb{T}^d, \\
    \chi_i \in H^1(\mathbb{T}^d)\quad \int_{\mathbb{T}^d}\chi_i(y)=0
  \end{aligned}
  \right.
\end{equation}

Consider the following problem:

\paragraph{Inverse homogenization problem:} Given the effective matrix $\sigma_e$ find $\sigma$.\\

This problem belongs to a class of problems in engineering called inverse homogenization, structural or shape optimization corresponding to the computation of the
microstructure of a material from its effective or homogenized properties  or the optimization of effective properties with respect  to microstructures belonging to an ``admissible set''. These problems are known to be ill posed in the sense that they don't admit a solution but a minimizing sequence of designs. It is possible to characterize the limits of the sequences by following the theory of G-convergence \cite{MR0477444} as observed in \cite{MR885713, MR844873}. For non-symmetric matrices the notion of H-convergence has been introduced \cite{Mur78, MR844873}. The modern theory for the optimal design of materials is the relaxation method through homogenization
\cite{ MR844873, MR885713, MR1763123, MR1211415, MR1804685}. This theory has lead to numerical  methods allowing for the design of nearly optimal micro-structures \cite{MR2008524, MR1859696, MR2270127}. We also refer to \cite{MR1899805, MR1862782, MR1493036} for the related theory of composite materials.

\begin{figure} \begin{center}
\xymatrix{ *-------{} & \txt{\footnotesize \it Physical \\\footnotesize
\it conductivity space} && \txt{\footnotesize \it Divergence-free \\ \footnotesize \it matrix space }
&& \txt{\footnotesize \it Convex functions\\ \footnotesize  \it space }\\
 *-------{\txt{\Large $\Omega$}} & *+[F-,]{\txt{{\huge $\sigma$}}} \ar@{~>} @<1ex>
 [rr]^{\txt{\footnotesize non-linear bijection}} \ar@{~>}[rrddd]|-{\txt{\footnotesize  non-linear, \\ \footnotesize  non-injective }}
  & &
*+[F-,]{\txt{ {\huge $Q$}}} \ar@{~>}@<1ex>[ll]^{\txt{ \footnotesize \eref{kdjhlsjdhlekhlkhcd},
\eref{sdssdsd32frc}} } \ar@{->}@<1ex>[rr]^{\txt{\footnotesize linear bijection}}
\ar@{->}[ddd]|-{\txt{\footnotesize  \eref{ksjdhslkhkkee}\\ \footnotesize linear\\ \footnotesize non injective}} & & *+[F-,]{\txt{{\huge
$s$}}}\ar@{->}@<1ex>[ll]^{\txt{}}\ar@{->}[dddll]\\ *-------{} & && &&\\ *-------{} & && &&\\
*-------{\txt{effective\\ conductivity}} & && *+[F-,]{\txt{ {\huge $\sigma_e$}}}
 &&
  }
\end{center} \caption[Relationships in inverse homogenization]{Illustration of theorem \ref{kdjshkjdhlshldhlk}. \label{fsdseedfofdsfdgspaces}}
\end{figure}
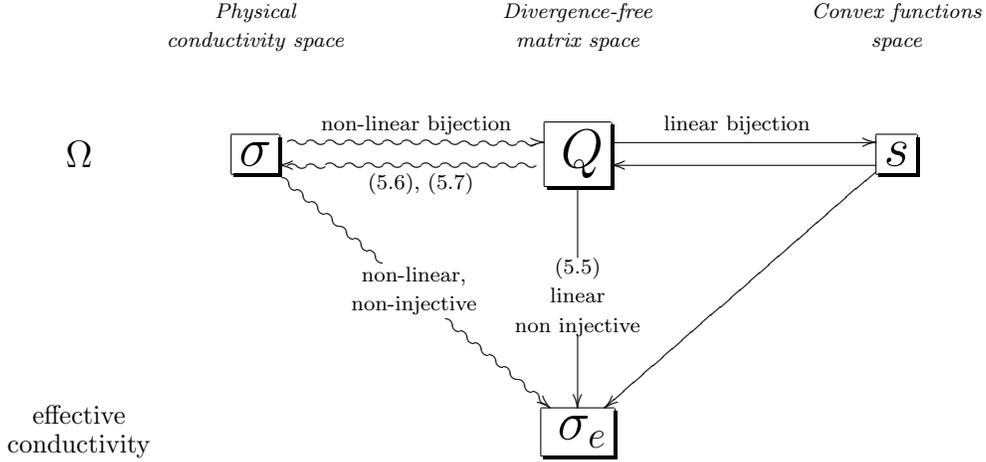

In this paper we observe that at least for the conductivity problem in dimension it is possible to transform the problem of looking for an optimal solution within a highly non-linear into the problem of looking for an optimal solution within a
linear space, as illustrated by theorem \ref{kdjshkjdhlshldhlk} and  Figure~\ref{fsdseedfofdsfdgspaces} for which efficient
optimization algorithms could be developed.

Define $F(y):=y+\chi(y)$ where $\chi$ is defined by \eref{sdssdssddsfrc}.

\begin{Theorem}\label{kdjshkjdhlshldhlk}
Let $Q$ be defined by \eref{edgeconddfdsuct}, then
\begin{enumerate}
\item $Q$ is divergence-free, periodic and associated with a convex function $s$ on $[0,1]^2$ through \eref{kjshwjhgkjwhg}.
\item
\begin{equation}\label{ksjdhslkhkkee}
\sigma_e=\int_{\mathbb{T}^d} Q(y)\,dy
\end{equation}
\item If $\sigma$ is isotropic then
\begin{equation}\label{kdjhlsjdhlekhlkhcd}
\sigma=\sqrt{\det(Q)\circ G^{-1}} I_d
\end{equation}
where $G(y):=y+\bar{\chi}$, $\bar{\chi}:=(\bar{\chi}_1,\bar{\chi}_2)$ and
\begin{equation}\label{sdssdsd32frc}
  \left\{
  \begin{aligned}
    -\Div\left(\frac{Q}{\sqrt{\det(Q)}}(y) \Grad (y_i+\bar{\chi}_i(y))\right) = 0, \enspace &y\in\mathbb{T}^d, \\
    \bar{\chi}_i \in H^1(\mathbb{T}^d)\quad \int_{\mathbb{T}^d}\bar{\chi}_i(y)=0
  \end{aligned}
  \right.
\end{equation}
\end{enumerate}
\end{Theorem}

\begin{figure} \begin{center} \xymatrix{ *-------{} & \txt{\footnotesize \it Physical \\\footnotesize
\it conductivity space} && \txt{\footnotesize \it Divergence-free \\ \footnotesize \it matrix space }
&& \txt{\footnotesize \it Convex functions\\ \footnotesize  \it space }\\
 *-------{\txt{\Large $\Omega$}} & *+[F-,]{\txt{{\huge $\sigma$}}} \ar@{~>} @<1ex>
 [rr]^{\txt{\footnotesize non-linear bijection}} \ar@{~>}[rrddd]
 \ar@{~>}@<1ex>[rrddddd]|-{\txt{\footnotesize  non-linear, \\ \footnotesize  non-injective }}
  & &
*+[F-,]{\txt{ {\huge $Q$}}} \ar@{~>}@<1ex>[ll]^{\txt{ \footnotesize \eref{eQharmonic},
\eref{jwhgdjhwjh}} } \ar@{->}@<1ex>[rr]^{\txt{\footnotesize linear bijection}}
\ar@{->}[ddd]|-{\txt{\footnotesize  volume-averaging\\ \footnotesize linear, non-injective\\
\eref{klyiuyoiyu}}} & & *+[F-,]{\txt{{\huge
$s$}}}\ar@{->}@<1ex>[ll]^{\txt{}}\ar@{->}[ddd]|-{\txt{\footnotesize  linear-interpolation\\
\footnotesize linear, non-injective\\ \eref{klkglkhg}}}\\ *-------{} & && &&\\ *-------{} & && &&\\
*-------{\txt{\Large $\Omega_h$}} & && *+[F-,]{\txt{ {\huge $q^h$}}}
\ar@{->}@<1ex>[rr]^{\txt{\footnotesize linear bijection}} && *+[F-,]{\txt{ {\huge
$s^h$}}}\ar@{->}@<1ex>[ll]^{\txt{}}\ar@{~>}@<1ex>[lldd]|-{\txt{\footnotesize non-linear bijection}}\\
*-------{} & && &&\\ *-------{} & && *+[F-,]{\txt{{\huge DATA}}} \ar@{~~>}'[rr]|-{\txt{\footnotesize
non-linear,\\\footnotesize well posed}}'[rruu]  \ar@{~~>}'[ll]|-{\txt{\footnotesize
non-linear,\\\footnotesize ill posed}}'[lluuuuu] &&
  }
\end{center} \caption[Relationships in inverse homogenization]{Relationships between spaces
  in inverse homogenization. \label{f:invhomogspaces}}
\end{figure}
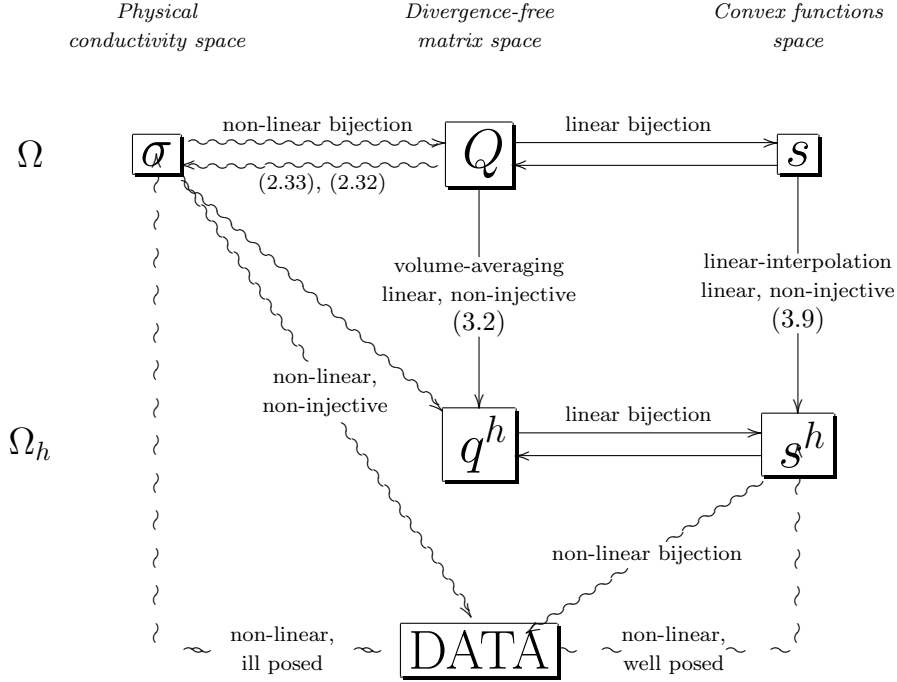

The problem of the computation of the
microstructure of a material from macroscopic information is not limited to inverse homogenization. Indeed in many ill posed inverse problems one can choose a scale coarse enough at which the problem admits a unique solution. Hence these problems can be formulated as the composition of a well posed (eventually non-linear) problem with an inverse homogenization problem.
The approach proposed here can also used to transform these problems (looking for an optimal solution within
a highly non-linear, non-convex space) into the problem of looking for an optimal solution within a
linear space, as illustrated in  Figure~\ref{f:invhomogspaces} for which efficient
 optimization algorithms can be used.  As an example, we examine Electric Impedance Tomography (EIT) in subsection \ref{s:inversenumeric}.

\section{Electric Impedance Tomography}\label{klcolhdckhcdeh}

We now apply our new approach to the inverse problem referred to as Electric Impedance Tomography
(EIT), which considers the electrical interpretation of~\eqref{e:basic}.  The goal is to determine
electrical conductivity from boundary voltage and current measurements, whereupon $\sigma(x)$ is an
image of the materials comprising the domain.  Boundary data is given as the Dirichlet-to-Neumann
(DtN) map, $\Lambda_{\sigma}:H^\frac{1}{2}(\partial\Omega) \rightarrow
H^{-\frac{1}{2}}(\partial\Omega)$, where this operator returns the electrical current pattern at the
boundary for a given boundary potential.

$\Lambda_{\sigma}$ can be sampled by solving the Dirichlet problem \begin{equation}
  \left\{
  \begin{aligned}
    -\Div(\sigma \Grad u) &= 0, \enspace &x\in\Omega, \\
    u &= g, \enspace &x\in\partial\Omega,
  \end{aligned}
  \right.
  \label{e:bvp}
\end{equation} and measuring the resulting Neumann data $f = \sigma \frac{\partial u}{\partial
  n}, x \in\partial\Omega$.  (In EIT, Neumann data is interpreted as electric
current.)  Although boundary value problem~\eqref{e:bvp} is not identical to the basic
problem~\eqref{e:basic}, we can still appeal to the homogenization results in~\cite{owha2005} provided
we restrict $g \in W^{2-\frac{1}{p}, p}$, in which case, we again have regular homogenization
solutions $\uhat \in W^{2,p}$, that is, those obtained by applying conductivity $Q(x)$ or $s(x)$.  If
$p>2$, a Sobolev embedding theorem gives $\uhat \in C^{1,\alpha}, \alpha>0$, already seen in
Section~\ref{jhskdgkjsghdkj}, although this restriction is not necessary for this section.

The EIT problem was first identified in the mathematics literature in the seminal 1980 paper by
Calder\'{o}n~\cite{cald1980}, although the technique had been known in geophysics since the 1930s.
We refer to \cite{MR2375322} and references therein for  simulated and experimental implementations of the method proposed by Calder\'{o}n.
From the work of Uhlmann, Sylvester, Kohn, Vogelius, Isakov and more recently, Alessandrini and
Vessella, we know that complete knowledge of $\Lambda_{\sigma}$ uniquely determines an {\em isotropic}
$\sigma(x) \in \Linf(\Omega)$, where $\Omega \subset \real^d, d \geq2$
~\cite{ales2005,kohn1984,sylv1987,isak1993}.

For a given diffeomorphism $F$ from $\Omega$ onto $\Omega$, write \begin{equation}\label{lkjjshjdhk}
  F_{\ast}\sigma:=\frac{(\nabla F)^T \sigma \nabla F}{\det (\nabla F)}\circ F^{-1}
\end{equation} It is known~\cite{gree2003} (see also \cite{gree2003a, sylv1990,nach1995, asta2005})
that for any diffeomorphism $F:\Omega \rightarrow \Omega,\quad F \in H^1(\Omega)$, the transformed
conductivity $\tilde{\sigma}(x) = F_{\ast} \sigma(x)$ has the same DtN map as the original
conductivity.  If $\sigma(x)$ is not isotropic, then $\tilde{\sigma} \neq \sigma$.  However, as shown
in~\cite{asta2005}, this is the {\em only} manner in which $\sigma(x) \in \Linf(\Omega)$ can be
non-unique.

Let $\Sigma(\Omega)$ be set of uniformly elliptic and bounded conductivities on $\Omega \in \real^2$,
that is, \begin{equation}
  \Sigma(\Omega) = \{\sigma \in \Linf(\Omega; \real^{2\times 2}) \, \mid \, \sigma
    = \sigma^T, 0< \lambda_{\min}(\sigma) <\lambda_{\max}(\sigma) < \infty \}.
\end{equation} The main result of~\cite{asta2005} is that $\Lambda_{\sigma}$ uniquely determines the
equivalence class of conductivities \begin{equation}
  \begin{aligned}
    E_{\sigma} = &\{ \sigma_1 \in \Sigma(\Omega) \mid \sigma_1 =
    F_{\ast}\sigma,\\
    & F : \Omega \rightarrow \Omega \quad \text{is an
    $H^1$-diffeomorphism and $F\mid_{\partial\Omega}=x$} \}.
  \end{aligned}
\end{equation} It has also been shown that there exists at most $\gamma \in E_{\sigma}$ such that
$\gamma$ is isotropic~\cite{asta2005}.

Our contributions in this section are as follows: \begin{itemize} \item Proposition~\ref{hasshljhedje}
gives an alternate and very simple proof of
  the uniqueness of an isotropic $\sigma \in E_\sigma$.  This is in contrast to
  Lemma 3.1 of~\cite{asta2005}, which appeals to quasi-conformal
  mappings. Moreover, Proposition~\eref{hasshljhedje} identifies isotropic
  $\gamma$ by explicit construction from an abribitrary $M\in E_{\sigma}$.
\item Proposition~\ref{kjdlskjhedjhe} shows that there exists equivalent classes
  $E_\sigma$ admitting no isotropic conductivities.
\end{itemize}
Proposition~\ref{kjhlkjhelkhwlkeh3e} shows that a given $\sigma \in
  \Sigma(\Omega)$ there exists a unique divergence-free matrix $Q$ such that
  $\Lambda_Q=\Lambda_\sigma$. It has been brought to our attention that proposition~\ref{kjhlkjhelkhwlkeh3e} has also been proven in \cite{MR2308972}. Hence, although for a given DtN map there may not exist an isotropic
  $\sigma$ such that $\Lambda=\Lambda_\sigma$,  there always exists a
  \emph{unique} divergence-free $Q$ such that $\Lambda=\Lambda_Q$. This is of
  practical importance since the medium to be recovered in a real application
  may not be isotropic and the associated EIT problem may not admit an isotropic
  solution.
Although the inverse of the map $\sigma\rightarrow \Lambda_\sigma$ is not continuous with respect to the topology of $G$-convergence when $\sigma$ is restricted to the set of isotropic matrices, it has also been shown in section 3 of \cite{MR2308972} that this inverse is continuous with respect to the topology of $G$-convergence when $\sigma$ is restricted to the set of divergence-free matrices.

We suggest from the previous observations and from theorem \ref{kjshgdkdjsgdhjg} that the space of convex functions on $\Omega$ is
  \emph{a natural space} in which to look for a parametrization of solutions of the EIT problem. In particular if an isotropic solution does exist, proposition \eref{hasshljhedje} allows for its recovery through the resolution the PDE \eref{eQharmsdsonic} involving the hessian of that convex function.

\begin{proposition}\label{hasshljhedje}
  Let  $\gamma \in E_{\sigma}$ such that $\gamma$ is isotropic. Then in dimension $d=2$,
  \begin{enumerate}
  \item $\gamma$ is unique.
  \item For any $M\in E_{\sigma}$
    \begin{equation}\label{jwhgdsdsdsdrrdsjhwjh}
      \gamma=\sqrt{\det(M)\circ G^{-1}}\, I_d,
    \end{equation}
    where $G=(G_1,G_2)$ are the harmonic coordinates associated to
    $\frac{M}{\sqrt{\det(M)}}$, that is, $G$ is the solution of
    \begin{equation}\label{eQharmsdsonic}
      \left\{
        \begin{aligned}
          \diiv \Big(\frac{M}{\sqrt{\det(M)}} \nabla G_i\Big)&=0, \enspace
          &x\in\Omega, \\
          G_i(x)&=x_i, \enspace &x\in\partial\Omega.
        \end{aligned}
      \right.
    \end{equation}
  \item $G=F^{-1}$ where $G$ is the transformation given
    by~\eqref{eQharmsdsonic}, and $F$ is the diffeomorphism mapping $\gamma$
    onto $M$ through equation~\eref{lkjjshjdhk}.
  \end{enumerate}
\end{proposition} \begin{proof}
  The proof is identical to the proof of Theorem~\ref{kjlaehlwhelkjwh}.
\end{proof}

\begin{proposition}\label{kjdlskjhedjhe}
  If $\sigma$ is a non isotropic, symmetric, definite positive, constant
  $2\times 2$ matrix, then there exists no isotropic $\gamma \in E_{\sigma}$.
\end{proposition} \begin{proof}
  Let us prove the proposition by contradiction.  Assume that $\gamma$
  exists. Then, it follows from Proposition~\ref{hasshljhedje} that $\gamma$ is
  constant and equal to $\sqrt{\det(\sigma)}I_d$. Moreover, it follows
  from~\eref{eQharmsdsonic} that $F^{-1}(x)=x$.  Using
  \begin{equation}
    \sigma=\frac{(\nabla F)^T \gamma \nabla F}{\det(\nabla F)}\circ F^{-1}
  \end{equation}
  we obtain that $\sigma$ is isotropic which is a contradiction.
\end{proof}

\begin{proposition}\label{kjhlkjhelkhwlkeh3e}
  Let $\sigma\in E_\sigma$.  Then in dimension $d=2$,
  \begin{enumerate}
  \item There exists a unique $Q$ such that $Q$ is a positive,
    symmetric divergence-free and $Q=F_{\ast}\sigma$. Moreover, $F$ are the
    harmonic coordinates associated to $\sigma$ given by~\eref{e:harmonic}.
  \item $Q$ is bounded and uniformly elliptic if and only if the non-degeneracy
    condition~\eref{sjdsgkjdksjhdg} is satisfied for an arbitrary $M \in
    E_\sigma$.
  \item $Q$ is the unique positive, symmetric, and divergence-free matrix
    such that $\Lambda_Q=\Lambda_\sigma$.
\end{enumerate} \end{proposition} \begin{proof}
  The existence of $Q$ follows from Proposition~\ref{kljssdlhkshdlk}. Let us
  prove the uniqueness of $Q$. If
  \begin{equation}
    Q=\frac{(\nabla F)^T \sigma \nabla F}{\det(\nabla F)}\circ F^{-1}
  \end{equation}
  is divergence free, then for all $l\in \R^d$ and $\varphi \in C_0^\infty(\Omega)$
  \begin{equation}
    \int_{\Omega} (\nabla \varphi)^T Q\cdot l=0.
  \end{equation}
  Using the change of variables $y=F(x)$ we obtain that
  \begin{equation}
    \int_{\Omega} (\nabla \varphi)^T Q\cdot l=\int_{\Omega} (\nabla (\varphi\circ F))^T \sigma \nabla
    F\cdot l.
  \end{equation}
  It follows that $F$ are the harmonic coordinates associated to $\sigma$ which
  proves the uniqueness of $Q$. The second part of the proposition follows from
  Proposition~\ref{kljssdlhkshdlk}.
\end{proof}

\subsection{Numerical reconstructions with incomplete boundary measurements using geometric
homogenization} \label{s:inversenumeric}

We close by examining two numerical methods for recovering conductivities from incomplete boundary
data using the ideas of geometric homogenization.  By incomplete we mean that potentials and currents
are measured at only a finite number of points of the boundary of the domain.  (For example, we have
data at $8$ points in Figure~\ref{f:Fiteration} for the medium shown in Figure~\ref{f:sigma2orig}.)

The first method is an iteration between the harmonic coordinates $F(x)$ and the conductivity
$\sigma(x)$.  The second recovers $s^h(x)$ from incomplete boundary data, and from $s^h(x)$ we compute
$q^h(x)$, then $Q$.  Both methods regularize the reconstruction in a natural way as to provide {\em
super-resolution} of the conductivity in a sense we now make precise.

The inverse conductivity problem with an imperfectly known
              boundary has also been considered in  \cite{MR2203860}. We refer to \cite{MR2300316}
           and reference therein for an analysis of the reconstruction of realistic conductivities
from noisy EIT data (using the D-bar method by studying its application to piecewise
smooth conductivities).

Even with complete boundary data this inverse problem is ill-posed with respect to the resolution of
$\sigma(x)$.  The Lipschitz stability estimate in~\cite{ales2005} states \begin{multline}
  \norm{\sigma^{(1,N)}(x)-\sigma^{(2,N)}(x)}_{\Linf}  \\ \leq C(N)
  \norm{\Lambda_{\sigma^{(1,N)}}-\Lambda_{\sigma^{(2,N)}}}_{\LambdaOpNorm},
\end{multline} where $\LambdaOpNorm$ is the natural operator norm for the DtN map. $\sigma^{(j,N)}(x)$
are scalar conductivities satisfying the ellipticity condition $0<\lambda\leq
\sigma(x)\leq\lambda^{-1}$ almost everywhere in $\Omega$ and belonging to a finite-dimensional space
such that \begin{equation}
  \sigma^{(j,N)}(x) = \sum_{i=1}^N \sigma^{(j,N)}_i z_i^{(N)}(x)
  \label{e:isoparamsigma}
\end{equation} for known basis functions $z_i^{(N)}(x)$.  Thus, the inverse problem in this setting is
to determine the real numbers $\sigma^{(j,N)}_i$ from the given DtN map $\Lambda_{\sigma^{(j,N)}}$.

The Lipschitz constant $C(N)$ depends on $\lambda, \Omega$ and $z_i^{(N)}$.  As shown by construction
in~\cite{rond2006}, when $z_i^{(N)}$ are characteristic functions of $N$ disjoint sets covering
$\Omega \subset \real^d$, the bound \begin{equation}
  C(N) \geq A \exp \left( B N^{\frac{1}{2d-1}} \right)
  \label{e:Cresolution}
\end{equation} for absolute constants $A, B>0$ is sharp.  That is, the amplification of error in the
recovered conductivity with respect to boundary data error increases exponentially with $N$.

From~\eqref{e:Cresolution}, we infer a resolution limit on the identification of $\sigma(x)$.  Setting
$\bar{C}$ our upper tolerance for the amplification of error in recovering $\sigma(x)$ with respect to
boundary data error, and introducing resolution $\bar{r} = N^{-1/d}$, which scales with length, we
have \begin{equation}
  \bar{r} \geq \left( \frac{1}{B} \log \frac{\bar{C}}{A} \right)^{-\frac{2d-1}{d}}.
\end{equation} We refer to any features of $\sigma(x)$ resolved at scales greater than this limit as
{\em stably-resolved} and knowledge of features below this limit as {\em super-resolved}.

\subsubsection{Harmonic coordinate iteration.} The first method provides super-resolution in two steps.
First, we stably-resolve conductivity using a resistor-network interpretation.  From this stable
resolution, we super-resolve conductivity by computing a function $\sigma(x)$ and its harmonic
coordinates $F(x)$ consistent with the stable resolution.

To solve for the conductivity at a stable resolution, we consider a coarse triangulation of $\Omega$.
Assigning a piecewise linear basis over the triangulation gives the edge-wise conductivies
\begin{oldequation}{edgesddsnzdduct}
  q_{ij}^h:=-\int_{\Omega} (\nabla \varphi_i )^T Q(x) \nabla \varphi_j \,dx.
\end{oldequation} As we have already examined, when $\sigma(x)$ (hence $Q(x)$) is known, so too are
the $q_{ij}^h$.  The discretized inverse problem is, given data at boundary vertices, determine an
appropriate triangulation of the domain and the $q_{ij}^h$ over the edges of the triangulation.  We
next specify our discrete model of conductivity in order to define what we mean by ``boundary data.''

Let $\calV_I$ be the set of interior vertices of a triangulation of $\Omega$, and let $\calV_B$ be the
boundary vertices, namely, the set of vertices on $\partial\Omega$.  Let the cardinality of $\calV_B$
be $V_B$.  Suppose vector $u^{(k)}$ solves the matrix equation \begin{equation}
  \left\{
  \begin{aligned}
    \sum_{j\sim i} \qh_{ij} (u^{(k)}_i-u^{(k)}_j) &= 0, \enspace& i \in \calV_I, \\
    u^{(k)}_i &= g^{(k)}_i, \enspace& i \in \calV_B,
  \end{aligned}
  \right.
  \label{e:discreteqharmonic}
\end{equation} where $g^{(k)}$ is given {\em discrete Dirichlet data}.  Then we define
\begin{equation}
  f^{(k)}_i = \sum_{j \sim i} \qh_{ij} (u^{(k)}_i-u^{(k)}_j), \enspace i \in
  \calV_B
\end{equation} as the {\em discrete Neumann data}.  The $V_B$ linearly independent $g^{(k)}$ and their
associated $f^{(k)}$ together determine the matrix $\Lambda^{\calV_B}_{\qh}$, which we call the {\em
discrete Dirichlet-to-Neumann map}. $\Lambda^{\calV_B}_{\qh}$ is linear, symmetric, and has the vector
$g = (1,1,\ldots,1)$ as its nullspace.  Hence, $\Lambda^{\calV_B}_{\qh}$ has $V_B(V_B-1)/2$ degrees of
freedom.

In practise, the discrete DtN map is provided as problem data without a triangulation specified: only
the boundary points where the Dirichlet and Neumann data are experimentally collected are given.  We
refer to this experimentally-determined discrete DtN map as $\Lambda^{\calV_B}_{\sigma}$.

We are also aware that to make sense in the homogenization setting, the $\qh_{ij}$ must be discretely
divergence-free.  That is, we require that \begin{equation}
  \sum_{j\sim i} \qh_{ij} (x^{(p)}_i-x^{(p)}_j) = 0, \enspace i \in \calV_I,
  \enspace p=1,2,
\end{equation} where $(x^{(1)}_i,x^{(2)}_i)$ is the $xy$-location of vertex $i$.

Set $\calT^{\calV_B}$ the set of triangulations having boundary vertices $\calV_B$ specified by
$\Lambda^{\calV_B}_{\sigma}$, and $\{\qh_{ij}\}$ the edge-values over ${\cal T}^{\calV_B}$.  The
complete problem is \begin{equation}
  \left\{
    \begin{aligned}
      &\mathop{\minimise}_{\calT^{\calV_B}, \{\qh_{ij}\}} \enspace
      \norm{\Lambda^{\calV_B}_{\qh}-\Lambda^{\calV_B}_{\sigma}}_{\ast}, \\
      &\subjectto \enspace \{\qh_{ij}\} \, \text{discretely divergence-free}.
    \end{aligned}
  \right.
\end{equation} The norm $\norm{\cdot}_{\ast}$ is a suitable matrix norm---as a form of regularization,
we use a thresholded spectral norm which under-weights error in the modes associated to the smallest
eigenvalues.  We solve this non-convex constrained problem using Constrained Simulated Annealing
(CSA). See~\cite{wah1999}, for example, for details on the CSA method.

EIT has already been cast in a similar form in~\cite{borc2007}, where edge-based data was solved for
using a finite-volume treatment, interpreting edges of the graph which connects adjacent cells as
electrical conductances.  They determine the edge values using a direct calculation provided by the
inverse theory for resistor networks~\cite{curt1990,curt2000}.  Although our work shares some
similarities with this prior art, we do not assume that a connectivity is known a priori.

An inversion algorithm for tomographic imaging of high contrast media based on a resistor network
theory has also been introduced in \cite{MR1421651}. The algorithm of \cite{MR1421651} is based on the
results of an asymptotic analysis of the forward problem  showing that when the contrast of the
conductivity is high, the current flow can be roughly approximated by that of a resistor network. Here
our algorithm is based on geometric structures hidden in homogenization of divergence form elliptic
equations with rough coefficients.

Given an optimal triangulation $\calT^{\ast}$ and its associated stably-resolved $\{\qh_{ij}\}$
representing conductivity, we now compute a fine-scale conductivity $\sigma^f(x)$ consistent with our
edge values, as well as its harmonic coordinates $F(x)$.  To help us super-resolve the conductivity,
we also regularize $\sigma^f(x)$.

Set $\calT^f$ a triangulation which is a refinement of triangulation $\calT^{\ast}$ from the solution
to the stably-resolved problem.  Let $\sigma^f(x)$ be constant on triangles of $T^f$.  Suppose
coordinates $F(x)$ are given, and solve \begin{equation}
  \left\{
    \begin{aligned}
      &\mathop{\minimise} \enspace \norm{\sigma^f(x)}_{\ast}, \\
      &\subjectto \enspace -\int_{\Omega} (\Grad(\myphi_i \circ F))^T \sigma^f(x)
      \Grad(\myphi_j \circ F) = \qh_{ij}, \enspace i,j \in \calT^{\ast}.
    \end{aligned}
  \right.
  \label{e:sigmaTV}
\end{equation} Here, $\norm{\cdot}_{\ast}$ is some smoothness measure of $\sigma^f(x)$. Following the
success of regularization by total variation norms in other contexts, see~\cite{koen2004,acar1994} for
example (in particular we refer to \cite{MR2424823} and references therein for convergence results the regularization of the inverse conductivity problem with discontinuous conductivities using total variation and the Mumford-Shah functional.), we choose \begin{equation}
  \norm{z(x)}_{\ast} = \norm{z(x)}_{\TV} := \int_{\Omega} \abs{\Grad z(x)}.
\end{equation}
 This norm makes particular sense for typical test cases, where the conductivity takes
on a small number of constant values.  This ``cartoon-like'' scenario is common where a small blob of
unusual material is included within a constant background material.  The constraints
in~\eqref{e:sigmaTV} are linear in the values of $\sigma^f(x)$ on triangles of $\calT^f$, and the norm
is convex, so~\eqref{e:sigmaTV} is a convex optimization problem.  In particular, it is possible to
recast~\eqref{e:sigmaTV} as a linear program, see~\cite{boyd2004}, for example. We then use the GNU
Linear Programing Kit to do the optimization of this resulting linear program~\cite{glpk}.  Note that
we build our refined triangulation using Shewchuk's Triangle program~\cite{tria}.

The harmonic coordinates $F(x)$ are not in general known.  We set $F(x)=x$ initially, and following
the solution of~\eqref{e:sigmaTV}, we compute \begin{equation}
  \left\{
    \begin{aligned}
      -\Div(\sigma^f\Grad F) &=0, \enspace &x \in \Omega, \\
      F&=x, \enspace &x \in\partial\Omega,
    \end{aligned}
  \right.
\end{equation} using $\sigma^f(x)$ from the previous step.  We can now iterate, returning to
solve~\eqref{e:sigmaTV} with these new harmonic coordinates.

\begin{figure}
  \begin{center}
    \includegraphics[width=0.34\textwidth]{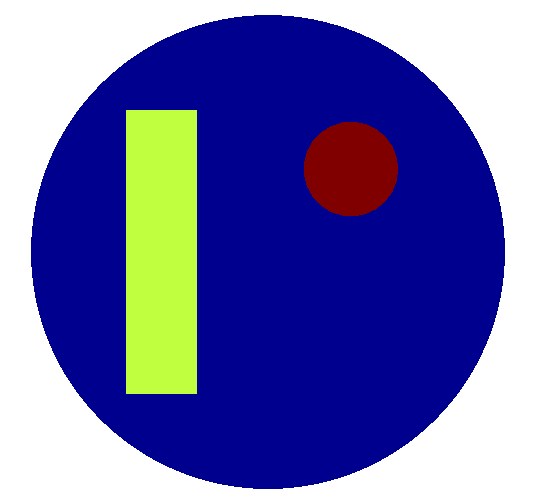}
    \hspace{0.1\textwidth}
    \includegraphics[width=0.35\textwidth]{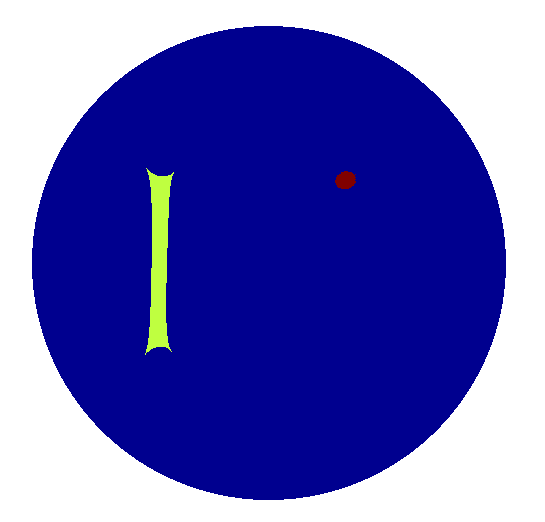}
  \end{center}
  \caption[A sample isotropic conductivity for testing
    reconstruction]{A sample isotropic conductivity for testing
    reconstruction.  The image on the left is $\sigma$, while the
    image on the right is $\sqrt{\det Q} = \sigma \circ F^{-1}$.  The
    dark blue background has conductivity 1.0, the red circle has
    conductivity 10.0 and the yellow bar has conductivity 5.0.  In
    this case, all of the features shrink in harmonic coordinates.
    \label{f:sigma2orig}}
\end{figure}

\begin{figure}
  \begin{center}
    \includegraphics[width=0.49\textwidth]{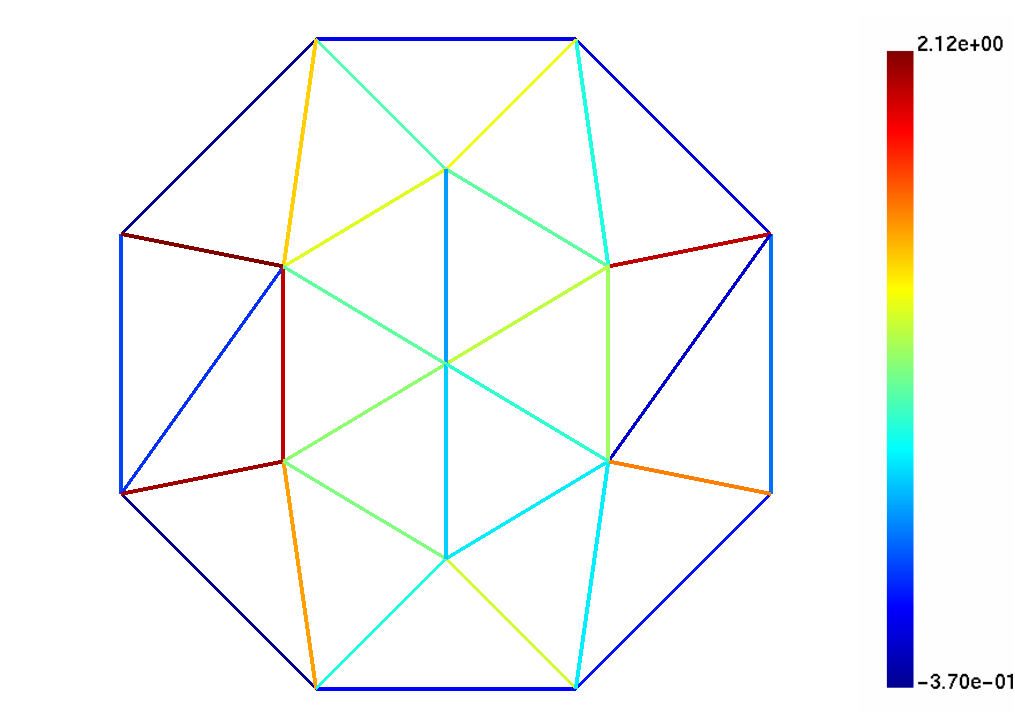}
    \includegraphics[width=0.49\textwidth]{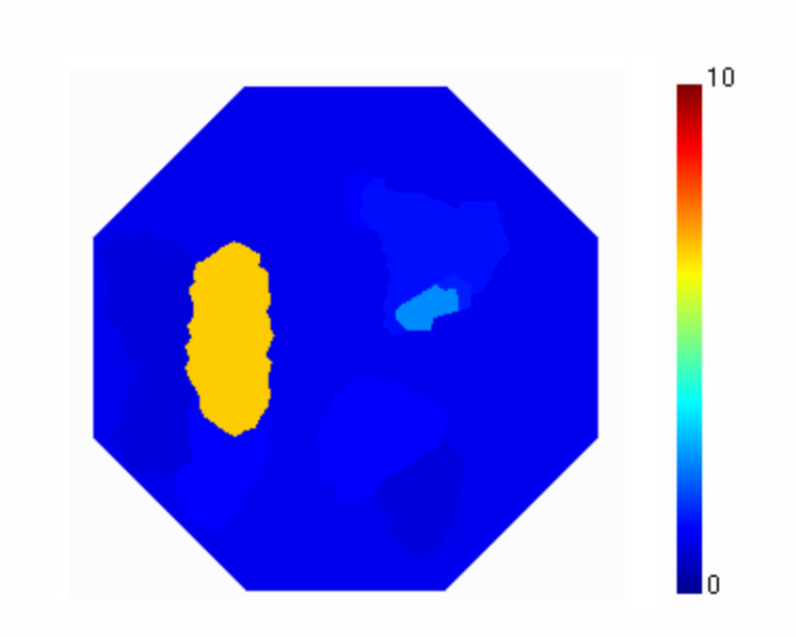} \\
    \includegraphics[width=0.49\textwidth]{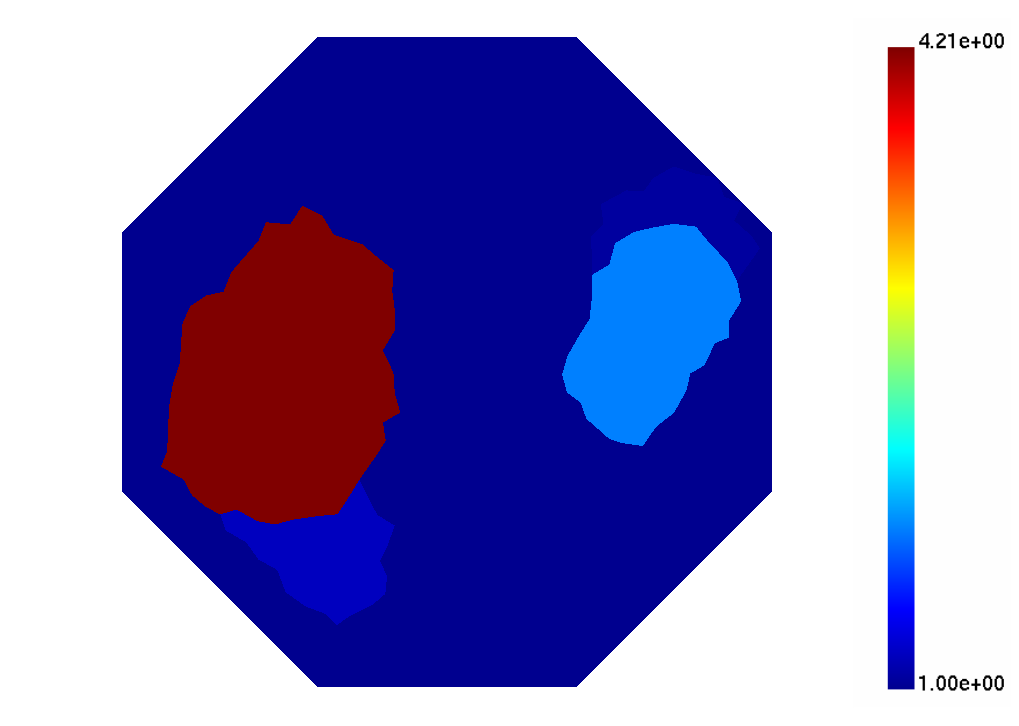}
    \includegraphics[width=0.49\textwidth]{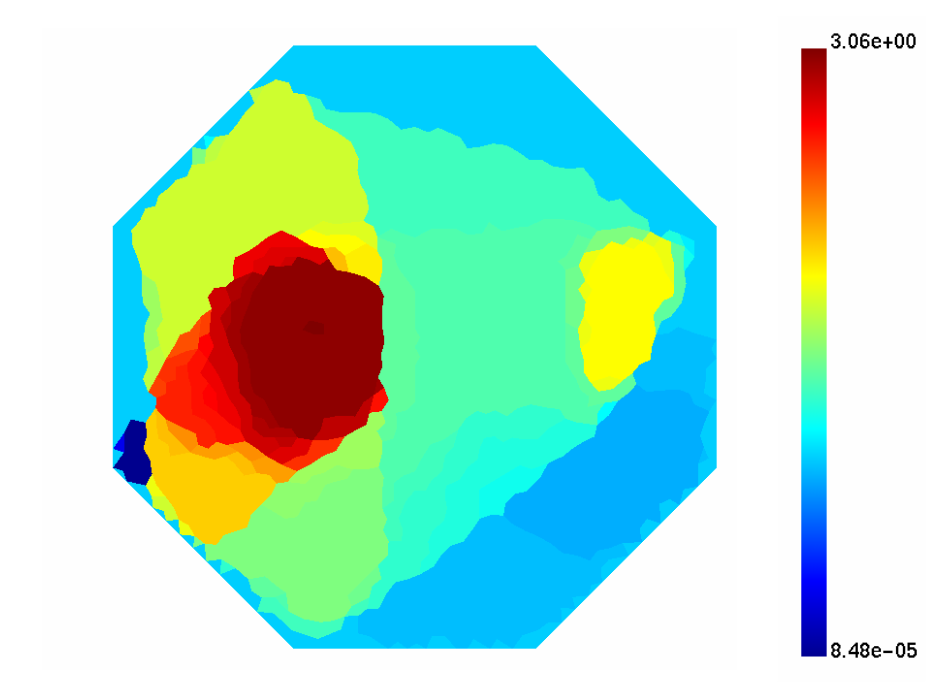}
  \end{center}
  \caption[Output of the harmonic coordinate iteration]{Output of the
    harmonic coordinate iteration. The figure on the top left is the
    coarse mesh produced by simulated annealing, the input to the
    harmonic coordinate iteration.  Left to right, top to bottom, the
    remaining three images show the progression of the iteration at 1,
    10, and 20 steps, showing its instability.  The true conductivity
    is that of Figure~\ref{f:sigma2orig}. \label{f:Fiteration}}
\end{figure}

Figures~\ref{f:sigma2orig} and~\ref{f:Fiteration} show the results of a numerical experiment
illustrating the method.  In particular, the harmonic coordinate iteration resolves details of the
true conductivity at scales below that of the coarse mesh used to resolve $\{\qh_{ij}\}$.  We observe
numerically that this iteration can become unstable and fail to converge.  However, before becoming
unstable the algorithm indeed super-resolves the conductivity.  We believe that this algorithm can be
stabilized and plan to investigate its regularization in a future paper.

\subsubsection{Divergence-free parameterization recovery.} Our second numerical method computes $s(x)$
from boundary data in one step.  In essence, we recover the divergence-free conductivity consistent
with the boundary data, without concern for the fine-scale conductivity that gives rise to the
coarse-scale anisotropy.

We begin by tessellating $\Omega$ by a fine-scale Delaunay triangulation, and we parameterize
conductivity by $\sh_i$, the piecewise linear interpolants of $s(x)$ over vertices of the
triangulation.  From the $\sh_i$, we can compute $\qh_{ij}$ using the hinge formula in order to solve
the discretized problem~\eqref{e:discreteqharmonic}.  This determines the discrete DtN map
$\Lambda_{\sh}$ in this setting.

\begin{figure}
  \begin{center}
\begin{picture}(0,0)%
\includegraphics{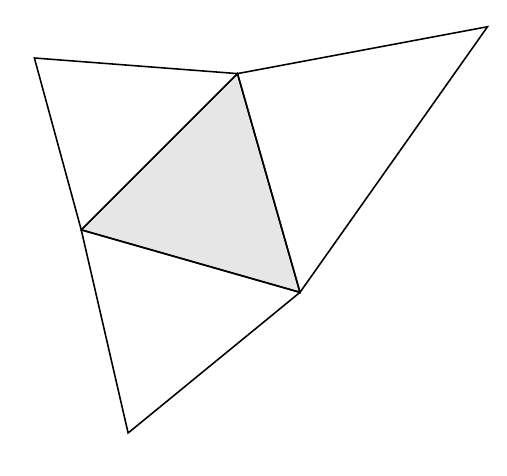}%
\end{picture}%
\setlength{\unitlength}{3947sp}%
\begingroup\makeatletter\ifx\SetFigFont\undefined%
\gdef\SetFigFont#1#2#3#4#5{%
  \reset@font\fontsize{#1}{#2pt}%
  \fontfamily{#3}\fontseries{#4}\fontshape{#5}%
  \selectfont}%
\fi\endgroup%
\begin{picture}(2430,2278)(2161,-2441)
\put(2701,-2386){\makebox(0,0)[lb]{\smash{{\SetFigFont{10}{12.0}{\familydefault}{\mddefault}{\updefault}{\color[rgb]{0,0,0}$m$}%
}}}}
\put(4576,-286){\makebox(0,0)[lb]{\smash{{\SetFigFont{10}{12.0}{\familydefault}{\mddefault}{\updefault}{\color[rgb]{0,0,0}$n$}%
}}}}
\put(2176,-436){\makebox(0,0)[lb]{\smash{{\SetFigFont{10}{12.0}{\familydefault}{\mddefault}{\updefault}{\color[rgb]{0,0,0}$l$}%
}}}}
\put(3226,-436){\makebox(0,0)[lb]{\smash{{\SetFigFont{10}{12.0}{\familydefault}{\mddefault}{\updefault}{\color[rgb]{0,0,0}$i$}%
}}}}
\put(3676,-1561){\makebox(0,0)[lb]{\smash{{\SetFigFont{10}{12.0}{\familydefault}{\mddefault}{\updefault}{\color[rgb]{0,0,0}$k$}%
}}}}
\put(2401,-1261){\makebox(0,0)[lb]{\smash{{\SetFigFont{10}{12.0}{\familydefault}{\mddefault}{\updefault}{\color[rgb]{0,0,0}$j$}%
}}}}
\end{picture}%
  \end{center}
  \caption[Stencil for approximating $Q$]{Stencil for approximating
    $Q(x)$ on triangle $ijk$ from the interpolants $\sh_i$ at nearby
    vertices. \label{f:quadratic}}
\end{figure}
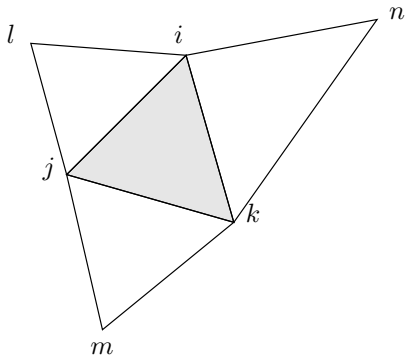

We shall also need a relationship between the $\sh_i$ and $\Qh_{ijk}$, an approximation of $Q(x)$
constant on triangles.  One choice is to presume that $s(x)$ can be locally interpolated by a
quadratic polynomial at the vertices of each triangle, and the opposite vertices of its three
neighbours, see Figure~\ref{f:quadratic}.  Taking second partial derivatives of this quadratic
interpolant gives a linear relationship between $\Qh_{ijk}$ and the six nearest $\sh_i$. This
quadratic interpolation presents a small difficulty, as triangles at the edge of $\Omega$ have at most
two neighbours.  Our solution is to place ghost vertices outside the domain near each boundary edge,
thus extending the domain of $s(x)$ and adding points where $\sh_i$ must be determined.

We solve for the $\sh_i$ using optimization by an interior point method. Although the algorithm we
choose is intended for convex optimization---the non-linear relationship between $\Lambda_{\sh}$ and
the $\sh_i$ makes the resulting problem non-convex---we follow the practise in the EIT literature of
relying on regularization to make the algorithm stable~\cite{borc2003,dobs1992a}.  We thus solve
\begin{equation}
  \left\{
    \begin{aligned}
      &\minimise \enspace \onehalf \sum_{k=1}^K \norm{\Lambda_{\sh} g^{(k)} -
        f^{(k)}}^2_{L^2(\partial\Omega)} + \alpha \norm{\tr \Qh}_{\TV}, \\
      &\subjectto \enspace \qh_{ij} \geq 0.
    \end{aligned}
  \right.
\end{equation} We use IpOpt software package for the optimization~\cite{ipopt}.

The data are provided as $K$ measured Dirichlet-Neumann pairs of data, $\{(g^{(k)},f^{(k)})\}$, and
the Tikhonov parameter $\alpha$ is determined experimentally (a common method is the L-curve method).
Again, the total variation norm is used to evaluate the smoothness of the conductivity.  We could just
as well regularize using $\det Q$ rather than $\tr Q$.  Using the trace makes the problem more
computationally tractable (the Jacobian is easier to compute), and our experience with such
optimizations shows that regularizing with respect to the determinant does not improve our results.
We compute the Jacobian of the objective's ``quadratic'' term using a primal-adjoint method,
see~\cite{dobs1992a}, for example.

We constrain $\qh_{ij}\geq 0$ on all edges, despite the possibility that our choice of triangulation
may require that some edges should have negative values. Our reasons for this choice are practical:
edge-flipping in this case de-stabilizes the interior point method.  Moreover, numerical experiments
using triangulations well-adapted to $\sigma(x)$ do not give qualitatively better results.

\begin{figure}
  \begin{center}
    \includegraphics[width=0.49\textwidth]{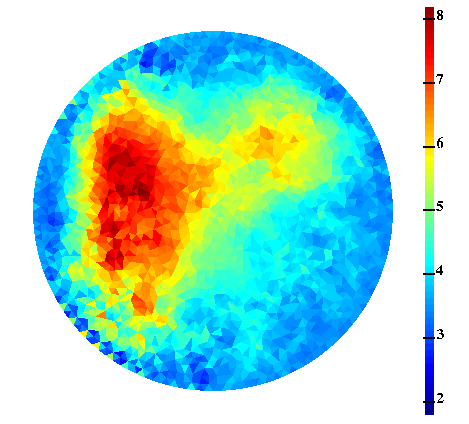}
    \includegraphics[width=0.49\textwidth]{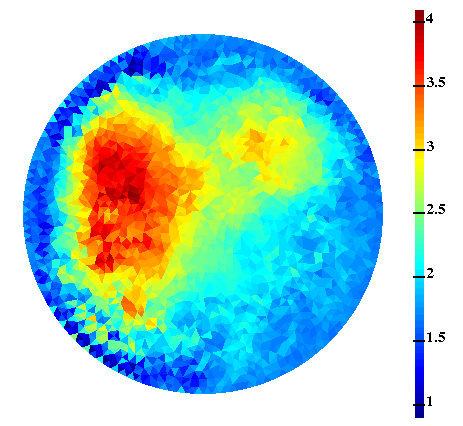}
  \end{center}
  \caption[Reconstruction of an isotropic conductivity]{Reconstruction
    of the isotropic conductivity in Figure~\ref{f:sigma2orig}.  The
    left-hand figure shows $\tr Q$, while the right-hand figure shows
    $\sqrt{\det Q}$.  The reconstruction blurs the original $\sigma$,
    similar to other methods in the literature, but does not
    underestimate the dynamic range of the large rectangle.
    \label{f:sigma2soln}}
\end{figure}

\begin{figure}
  \begin{center}
    \includegraphics[width=0.49\textwidth]{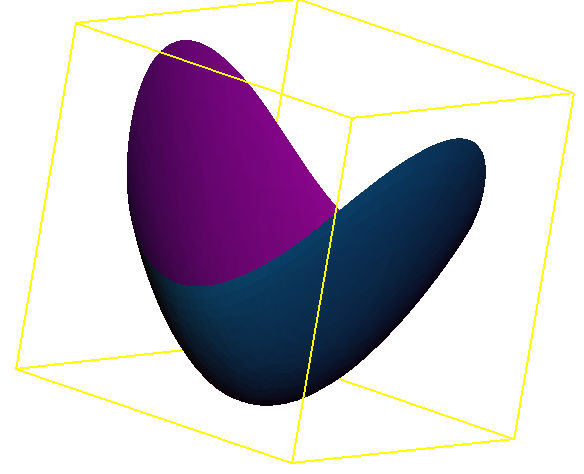}
    \includegraphics[width=0.49\textwidth]{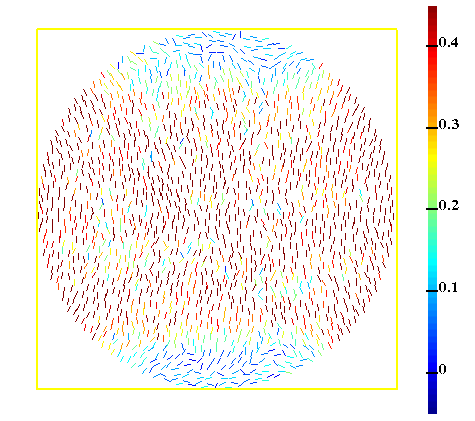}
  \end{center}
  \caption[Anisotropic reconstruction I]{Anisotropic reconstruction
    showing the parameterisation $s(x)$ we recover using EIT (left
    image), and the pattern of anisotropy we see by rendering the
    orientation of the maximal eigenvalue of $Q$ (right image).  The
    colour-bar in the right image indicates the strength of the
    anisotropy as $\abs{\lambda_{\max}-\lambda_{\min}}/\tr Q$.
    \label{f:anis5solnI}}
\end{figure}

Figures~\ref{f:sigma2soln} and~\ref{f:anis5solnI} show reconstructions of the conductivities in
Figures~\ref{f:sigma2orig} and~\ref{f:anis4sigma}, respectively.  We include the reconstruction of the
conductivity in Figure~\ref{f:sigma2orig} only to show that our parameterization can resolve this test
case, a typical one in the EIT literature.  For such tests recovering ``cartoon blobs,'' our method
does not compete with existing methods such as variational approaches~\cite{borc2003}, or those based
on quasi-conformal mappings~\cite{isaa2006b,knud2004}.  Our recovery of the conductivity in
Figure~\ref{f:anis4sigma}, however, achieves a reconstruction, to our knowledge, which has not
previously been realized.  The pitch of the laminations in this test case are below a reasonable limit
of stable resolution.  Hencee do not aim to recover the laminations themselves, but we do recover
their up-scaled representation.  The anisotropy of this up-scaled representation is apparent in
Figure~\ref{f:anis5solnI}.  Admitting the possibility of recovering anisotropic, though
divergence-free, conductivities by parameterizing conductivity by $s(x)$ gives a sensical recovered
conductivity.  Owing to the stable resolution limit, parameterizing $\sigma(x)$ directly, by a usual
parameterization---such as the linear combination in~\eqref{e:isoparamsigma}, choosing $z_i(x)$ as
characteristic functions---is not successful. 

\begin{flushleft} 
\def\cprime{$'$}
\providecommand{\bysame}{\leavevmode\hbox to3em{\hrulefill}\thinspace}
\providecommand{\MR}{\relax\ifhmode\unskip\space\fi MR }
\providecommand{\MRhref}[2]{%
  \href{http://www.ams.org/mathscinet-getitem?mr=#1}{#2}
}
\providecommand{\href}[2]{#2}

\end{flushleft}

\end{document}